\documentclass[11pt]{elsarticle}
\usepackage{graphicx}
\usepackage{float}
\usepackage{subfig}
\usepackage{latexsym,amsmath,amsfonts,amscd, amsthm, dsfont}
\usepackage{nameref}

\usepackage[percent]{overpic}
\usepackage{cleveref}
\usepackage{bm,color}
\usepackage{epsfig,verbatim,epstopdf,graphics}
\usepackage{changebar}
\usepackage{multirow}
\usepackage{enumerate}
\usepackage{yhmath}
 \usepackage{booktabs} 
 \usepackage{tikz}
\usepackage{verbatim}
\usepackage{appendix}
\usetikzlibrary{arrows,backgrounds,snakes,shapes}
 \numberwithin{equation}{section}

\graphicspath{{./}{./figure/}}
\allowdisplaybreaks

\newtheoremstyle{plainNoItalics}{}{}{\normalfont}{}{\bfseries}{.}{ }{}

\theoremstyle{plain}
\newtheorem{thm}{Theorem}[section]
\theoremstyle{plainNoItalics}

\newtheorem{rem}[thm]{Remark}

\newtheorem{exa}[thm]{Example}

\newcommand{\be}{\begin{eqnarray}}
\newcommand{\ee}{\end{eqnarray}}
\newcommand{\beno}{\begin{eqnarray*}}
\newcommand{\eeno}{\end{eqnarray*}}

\makeatletter

\newcommand{\Rmnum}[1]{\expandafter\@slowromancap\romannumeral #1@}
\makeatother

\begin{document}

 \baselineskip=1.8pc

\begin{frontmatter}
\title{Fourth-order conservative non-splitting semi-Lagrangian Hermite WENO schemes for kinetic and fluid simulations}

\author[1]{Nanyi Zheng}
\ead{nyzheng@stu.xmu.edu.cn}
\author[2,3]{Xiaofeng Cai\corref{cor1}}
\ead{xfcai@bnu.edu.cn}
\author[4]{Jing-Mei Qiu}
\ead{jingqiu@udel.edu}
\author[5]{Jianxian Qiu}
\ead{jxqiu@xmu.edu.cn}
\cortext[cor1]{Corresponding authors}
\address[1]{School of Mathematical Sciences, Xiamen University, Xiamen, Fujian 361005, China}
\address[2]{Research Center for Mathematics, Beijing Normal University,
	Zhuhai 519087, China}
 \address[3]{ BNU-HKBU United International College, Zhuhai 519087, China}
\address[4]{Department of Mathematical Sciences, University of Delaware, Newark, DE, 19716, USA}
\address[5]{School of Mathematical Sciences and Fujian Provincial Key Laboratory of Mathematical Modeling and High-Performance Scientific Computing, Xiamen University, Xiamen, Fujian 361005, China}

\begin{abstract}
We present fourth-order conservative non-splitting semi-Lagrangian (SL) Hermite essentially non-oscillatory (HWENO) schemes for linear transport equations with applications for nonlinear problems including the Vlasov-Poisson system, the guiding center Vlasov model, and the incompressible Euler equations in the vorticity-stream function formulation. The proposed SL HWENO schemes combine a weak formulation of the characteristic Galerkin method with two newly constructed HWENO reconstruction methods. Fourth-order accuracy is accomplished in both space and time under a non-splitting setting. Mass conservation naturally holds due to the weak formulation of the characteristic Galerkin method and the design of the HWENO reconstructions. We apply a positive-preserving limiter to maintain the positivity of numerical solutions when needed. Although the proposed SL framework allows us to take large time steps for improving computational efficiency, it also brings challenges to the spatial reconstruction technique; we construct two kind of novel HWENO reconstructions to fit the need for the proposed SL framework. Abundant benchmark tests are performed to verify the effectiveness of the proposed SL HWENO schemes.

\vfill

\end{abstract}

\end{frontmatter}



\section{Introduction}

Semi-Lagrangian (SL) schemes have been developed in areas of applications such as numerical weather prediction \cite{Staniforth1991,Guo2014A,Bosler2019}, kinetic description of plasma \cite{FILBET2001166,Nicolas2010,ROSSMANITH20116203}, and fluid simulations \cite{Pironneau1982,XIU2001658}. Most SL schemes are designed to solve the transport equation in the form of
\begin{equation}\label{eq:transport_R_D}
	u_t + \nabla_{\mathbf{x}}\cdot\left(\mathbf{a}(u,\mathbf{x},t)u\right) = 0,
\end{equation}
where $u( \mathbf{x},t)$ usually represents a density function of a conservative quantity in a velocity field $\mathbf{a}(u,\mathbf{x},t)$ with $\mathbf{x}\in\mathbb{R}^d$. Comparing with the Eulerian and the Lagrangian approaches, the SL approach naturally holds its advantages in terms of accuracy and efficiency for certain applicable problems. On one hand, the same with the Lagrangian approach, the SL approach evaluates the solution along the convection characteristics. Hence, it allows large numerical time steps comparing with the Eulerian approach. On the other hand, as the Eulerian approach, the SL approach adopts a fixed spatial mesh equipped with a wide range of different solution spaces for high-order spatial accuracy. As a comparison, the Lagrangian approach suffers from statistical noises and only achieves a low order of $O(1/\sqrt{N})$ with $N$ representing the number of sampling points.

 The geometry structure of transport dynamics can be very complicated. For instance, the Vlasov-Poisson system has drastic high-frequency filamentation structures \cite{ROSSMANITH20116203,QIU20118386} and the guiding center Vlasov model has steep structures \cite{2019ConservativeXiong}. To handle such complicated structures, the SL approach has been successfully coupled with the discontinuous Galerkin (DG) method \cite{RESTELLI2006195,QIU20118386,lee2016a}, the  weighted essentially non-oscillatory (WENO) schemes \cite{QIU2011863,huang2016semi,SIRAJUDDIN2019619}, and the Hermite WENO (HWENO) schemes \cite{YANG201418,CAI201695,ZhengSLHWENO2021}. However, there are pros and cons for each type of spatial discretization. For the DG method, it needs many degrees of freedom (DOF) per element for its high-order version, especially in a high-dimensional setting.  For high-order WENO schemes, they always require very wide reconstruction stencils, compared with the DG method. The HWENO methods can be regarded as an intermediate transition from the DG method to a WENO method. It requires significantly lower degrees of freedom per element compared with the DG method, and it uses a more compact stencil for high-order reconstruction compared with a WENO method.

 The HWENO method was first introduced by Qiu and Shu in \cite{qiu2004hermite} and was further developed in \cite{QIU2005642,Cai2016Positivity,DU2018385,zhao2020hermite}. Notice that the HWENO schemes in \cite{qiu2004hermite,QIU2005642,Cai2016Positivity,DU2018385} use point-wise positive linear weights to present the information of a high-degree polynomial as a convex combination of the information of several low-degree polynomials. In \cite{qiu2004hermite,QIU2005642,Cai2016Positivity,DU2018385}, such positive linear weights do exist for the specific Guassian points they require. However, one can find that such positive linear weights do not exist for some special locations (even regions) in an Eulerian grid for each HWENO reconstruction in \cite{qiu2004hermite,QIU2005642,Cai2016Positivity,DU2018385}. Notice that characteristic feet for \eqref{eq:transport_R_D} can be located  anywhere on an Eulerian grid; and the reconstruction designs in \cite{qiu2004hermite,QIU2005642,Cai2016Positivity,DU2018385} are not suitable for an SL method. Based on the observation above, in our previous work \cite{ZhengSLHWENO2021}, we adopted a 1-D hybrid HWENO reconstruction method proposed in \cite{zhao2020hermite}, and developed a dimensional-splitting SL hybrid HWENO scheme. Although the numerical results of the splitting-based SL HWENO scheme are satisfying for most cases, we observe that the scheme is very dissipative for large-gradient extreme points. The reason for the dissipation problem is that both the troubled cell indicator and the HWENO reconstruction are not able to distinguish large-gradient extreme points from discontinuities. Hence, the hybrid HWENO reconstruction decays to a first-degree polynomial reconstruction near large-gradient extreme points.

 To overcome the dissipation issue, we propose two new two-dimensional (2-D) HWENO reconstruction methods, denoted by HWENO-1 and HWENO-2. Comparing with the HWENO reconstruction technique in \cite{zhao2020hermite}, a key difference of the newly constructed reconstructions is that we rule out the participation of any first-degree polynomial, which leads to large numerical dissipation. The two newly proposed HWENO methods gather  information from central or one-sided constructed polynomials, which are at least quadratic.
 Theoretically, a good HWENO reconstruction approximates the highest-degree polynomial where the solution is continuous and reduces to a one-sided lower-degree polynomial when a  discontinuity is involved.
 To achieve such a principle, we attempt two different strategies. For the HWENO-1, we follow the same technique in \cite{zhao2020hermite}, but base on newly chosen polynomials. The resulting HWENO-1 method has a significant improvement in reducing numerical dissipation. On the other hand, the HWENO-2 directly selects one polynomial from all candidates. This is equivalent to assigning weights of one or zero to candidate polynomials. The stencil selection strategy of the HWENO-2 is motivated by the targeted essentially non-oscillatory (TENO) schemes developed by Fu et al. \cite{FU2016333,FU2018724,FU2019117}.
 Recently, a Hermite TENO (HTENO) is developed in \cite{Wibisono2021HTENO}, which shows very good numerical results.
 There are two reasons not to apply the HTENO method proposed in \cite{Wibisono2021HTENO}. Firstly, it depends on point-wise positive linear weights as we mentioned before.
 Secondly, it is only a one-dimensional version method. Both the HWENO-1 and HWENO-2 have a huge improvement in terms of dissipation comparing with the HWENO reconstruction in \cite{zhao2020hermite}. But, there are subtle differences in effectiveness. On one hand, The HWENO-1 is ``smoother" due to its nature of assigning non-zero weight to each candidate polynomial. As a comparison, the HWENO-2 (or HTENO) only chooses one candidate, which means the piecewise polynomial constructed by the HWENO-2 can be fragmented sometimes. On the other hand, The HWENO-2 is more robust for capturing discontinuity since the highest-degree polynomial is not involved at all as long as the stencil selecting strategy works correctly. In this paper, we will describe the construction of the HWENO-1 and HWENO-2 in details and compare these two strategies by abundant numerical tests.

 The proposed non-splitting SL HWENO schemes adopt the weak formulation of the SL DG  method \cite{Guo2014A,CAI2018529}, which comes from the characteristic Galerkin method \cite{dahle1995eulerian,russell2002overview}, to update the zeroth-order and first-order moments. In other words, we define the solution space as a $P^1$ DG space. However, through the HWENO reconstructions, the $P^1$ DG solution is replaced with a piecewise $P^3$ polynomial for the solution evaluation. Such a procedure can be viewed as a $P_NP_M$ method introduced in \cite{dumbser2008unified}. The combination of the weak formulation of SL DG method and the HWENO reconstruction can also be regarded as a one-step evolution Galerkin scheme introduced in \cite{Morton1998}. The proposed SL scheme requires a remapping procedure between a fixed Eulerian mesh and a characteristic upstream twisted mesh. To accomplish a fourth-order accuracy, the remapping procedure can be summarized in the following three steps. Firstly, we define a cubic-curved numerical upstream mesh to approximate the real upstream mesh. Secondly, a clipping technique is involved to gather the curved polygons generated from the Eulerian mesh and the cubic-curved mesh. Finally, apply piecewise integration over each cubic-curved upstream cell. For nonlinear models, we couple the proposed SL HWENO schemes with a fourth-order Runge-Kutta exponential integrator (RKEI) \cite{CELLEDONI2003341}, which freezes the velocity field for each stage, for high-order temporal accuracy. Notice that the evaluation of zeroth-order moment (cell average) is equivalent to the formulation of the SL finite volume (FV) method in \cite{ZHENG2022114973}. The mass conservation, positivity preservation (PP), and fourth-order accuracy can be proved similar to \cite{ZHENG2022114973}. For stability, under a linearized setting, we numerically prove the unconditionally stable property of the proposed schemes by the Fourier analysis.

 An outline of this paper is as follows. In \Cref{sec:SL_FV_WENO_scheme}, we introduce the construction of the SL HWENO schemes. In \Cref{sec:nonlinear}, we describe how to couple the proposed SL HWENO scheme with the fourth-order RKEI method for nonlinear models. A variety of numerical tests are provided in \Cref{sec:numerical_tests}. Finally, a conclusion is presented in \Cref{sec:conclusion}.

\section{Two-dimensional SL HWENO schemes}\label{sec:SL_FV_WENO_scheme}
Consider the following linear transport equation
\begin{equation}
u_t + (a(x,y,t)u )_x + ( b(x,y,t)u )_y = 0,
\label{2d_linear}
\end{equation}
where $\left(a(x,y,t),b(x,y,t)\right)$ represents a known velocity field and $u(x,y,t)$ is a density function. We assume a rectangle computational domain $\Omega := [x_L,x_R]\times[y_B,y_T]$ and a corresponding discretization such that $x_L = x_{\frac12}<x_{\frac32}<\cdots<x_{N_x+\frac12}=x_R$, $y_B = y_{\frac12}<y_{\frac32}<\cdots<y_{N_y+\frac12}=y_T$, with $I^x_i:=\left[x_{i-\frac12},x_{i+\frac12}\right]$, $I^y_j:=\left[y_{j-\frac12},y_{j+\frac12}\right]$, $I_{i,j} :=I^x_i \times I^y_j$, $x_i:=\frac{x_{i-\frac12}+x_{i+\frac12}}2$, $y_j:=\frac{y_{j-\frac12}+y_{j+\frac12}}{2}$, $\Delta x_i := x_{i+\frac12}-x_{i-\frac12}$ and $\Delta y_j := y_{j+\frac12} - y_{j-\frac12}$. We define a $P^1$ DG solution space $V^1_h := \{ u_h|u_h(x,y)|_{I_{i,j}}\in P^1(I_{i,j})~\forall i,j \}$. Consider an Eulerian cell $I_{i,j}$ at $t = t^{n+1}$ and define a dynamic characteristic region
$I_{i,j}(t):=\{\left(x^*,y^*\right)|\left(x^*,y^*\right)=\left(X(x,y;t),Y(x,y;t)\right),~~(x,y)\in I_{i,j}\}$, where $\left(X(x,y;t),Y(x,y;t)\right)$ represents the characteristic curve emanating from $(x,y,t^{n+1})$, i.e., the solution of the ordinary differential equations (ODEs)
\begin{equation}\label{eq:ODEs}
  \begin{cases}
     d X(t)/dt = a(X(t) ,Y(t) ,t ),\\
      d Y(t)/dt = b(X(t) ,Y(t) ,t ), \\
     X(t^{n+1} ) = x,\\
     Y(t^{n+1} ) = y.
  \end{cases}
\end{equation}

We define an adjoint problem of \eqref{2d_linear} as in \cite{dahle1995eulerian,russell2002overview} on $I_{i,j}(t)\times[t^n,t^{n+1}]$: for a given test function $W(x,y)\in P^1(I_{i,j})$,
\begin{equation}\label{eq:adjoint}
\begin{cases}
\begin{split}
&w_t+a(x,y,t)w_x+b(x,y,t)w_y=0,\quad t\in [t^n,t^{n+1})\\
&w(t=t^{n+1})=W(x,y).
\end{split}
\end{cases}
\end{equation}
By the Reynolds transport Theorem, we have,
\begin{equation}\label{SL_2}
\begin{split}
&\frac{d}{dt}\iint_{I_{i,j}(t)}u(x,y,t)w(x,y,t)dxdy = 0.
\end{split}
\end{equation}

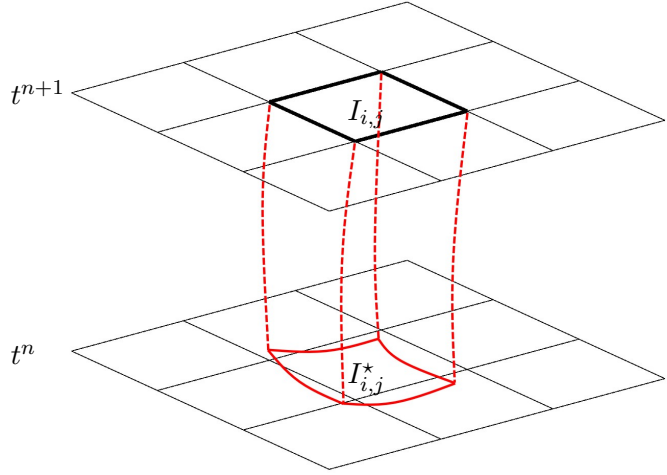
\begin{figure}[htb]
\begin{center}
\begin{overpic}[scale=0.4]{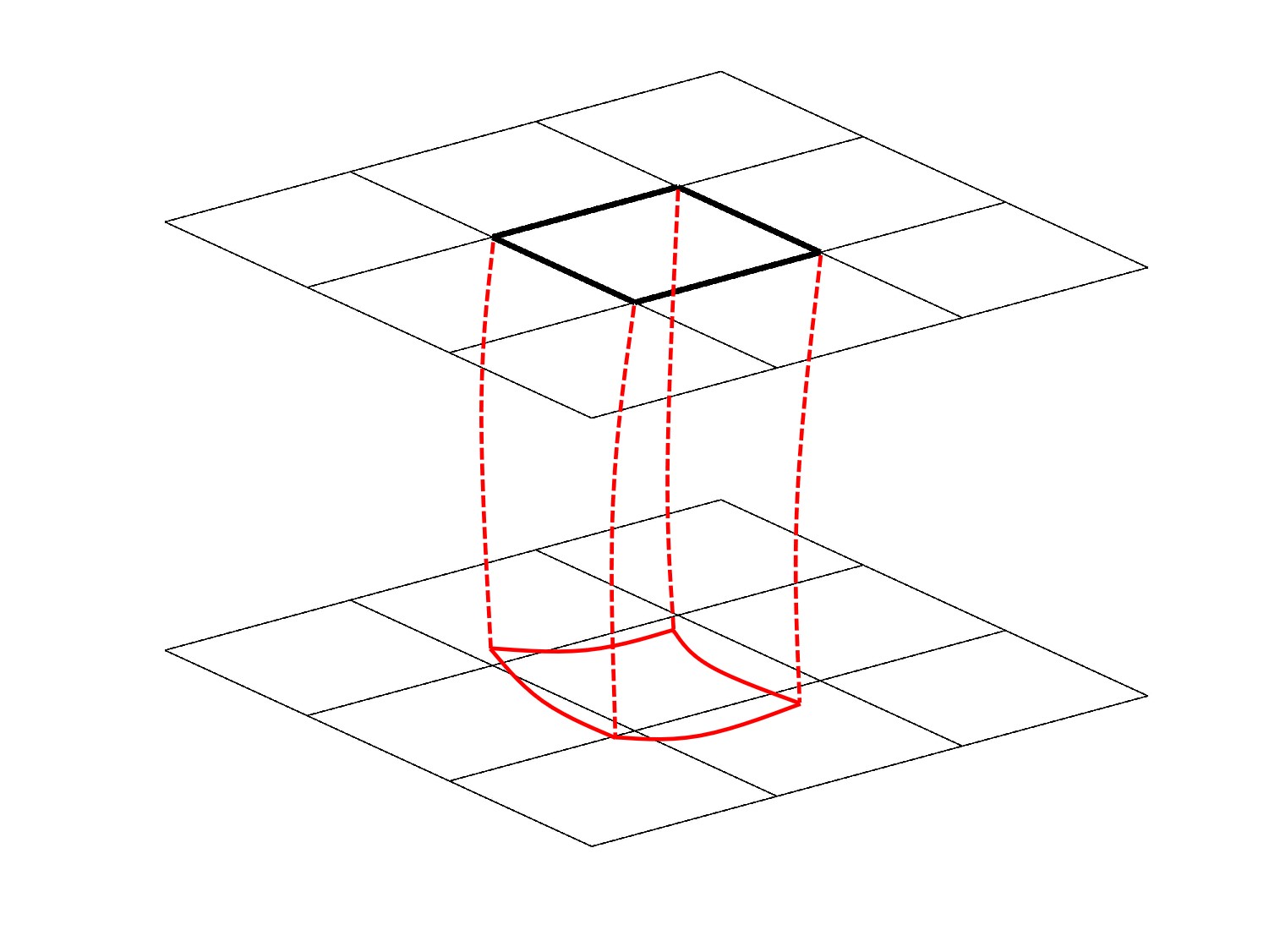}
\put(49,54){$I_{i,j}$}
\put(49,19.5){$I_{i,j}^\star$}
\put(5,22) {$t^n$}
\put(5,56){$t^{n+1}$}
\end{overpic}
\end{center}
\caption{Schematic  illustration for the characteristic upstream cell $I_{i,j}^\star$.}
\label{charac_region}
\end{figure}

From \eqref{SL_2}, an SL scheme is naturally formulated:
\begin{equation}
\frac1{\Delta x_i \Delta y_j}\iint_{I_{i,j}}u(x,y,t^{n+1})W(x,y)dxdy = \frac1{\Delta x_i \Delta y_j}
\iint_{I_{i,j}^\star} u(x,y,t^n)w(x,y,t^n)dxdy,\\
\label{SL_formulation}
\end{equation}
where $I^\star_{i,j}=I_{i,j}(t^n)$ (see \Cref{charac_region}).
For given time level $t^n$, we denote the first three moments of the solution by $\{\overline{u}_{i,j}\}$, $\{\overline{v}_{i,j}\}$, and $\{\overline{w}_{i,j}\}$. Then, we denote the numerical solution on $t^n$ by $u^n$ with
\begin{equation}
u^n(x,y) = \overline{u}_{i,j} + 12\overline{v}_{i,j}\left(\frac{x-x_i}{\Delta x_i}\right) + 12\overline{w}_{i,j}\left(\frac{y-y_j}{\Delta y_j}\right),\quad (x,y)\in I_{i,j}.
\end{equation}
For constructing an SL HWENO scheme, it is sufficient to take $W(x,y) = 1$, $(x-x_i)/\Delta x_i$, $(y-y_j)/\Delta y_j$ and evaluate the right-hand side of \eqref{SL_formulation} accurately. To approximate the right-hand side of \eqref{SL_formulation}, in \Cref{2dHWENO}, we first introduce the two newly constructed HWENO reconstructions to recover a piece-wise cubic polynomial, denoted by $H^n(x,y)$, to approximate $u(x,y,t^n)$ in \eqref{SL_formulation}. Then, in \Cref{cubic_upstream}, cubic-curved quadrilaterals, denoted by $\{\widetilde{I}^\star_{i,j}\}$, are defined for approximating $\{I_{i,j}^\star\}$ and a cubic polynomial $\widetilde{w}(x,y)$ is constructed over each $\widetilde{I}^\star_{i,j}$ by a least square procedure to approximate the test function $w(x,y)$. Finally, we briefly summarize the integration strategy on $\widetilde{I}^\star_{i,j}$ in \Cref{sec:numerical_integration}.

\subsection{Two-dimensional HWENO reconstruction methods}\label{2dHWENO}

For convenience, we require $\Delta x_i \equiv \Delta x,\quad \Delta y_j \equiv \Delta y\quad \forall i,j$. Based on the $P^1$ DG solution, $u^n(x,y)$, we recover a piecewise $P^3$ polynomial,
\begin{equation}
	H^n(x,y) = H^{(i,j)}(x,y),~~(x,y)\in I_{i,j},~~\forall (i,j),
\end{equation}
where $H^{(i,j)}(x,y)\in P^3(I_{i,j})$. We define a set of local orthogonal basis of the high order polynomial denoted as $\{P_l(x,y)\}$ with:
\begin{equation}\label{basis}
\begin{split}
  &P_1(x,y)=1,~~P_2(x,y)=\mu_i(x):=\frac{x-x_i}{\Delta x},~~P_3(x,y)=\nu_j(y):=\frac{y-y_j}{\Delta y},\\
  &P_4(x,y)=\mu_i^2-\frac1{12},~~P_5(x,y)=\mu_i\nu_j,~~P_6(x,y)= \nu_j^2 - \frac1{12},\\
  &P_7(x,y)=\mu_i^3-\frac3{20}\mu_i,~~P_8(x,y)=\left(\mu_i^2-\frac1{12}\right)\nu_j,~~P_9(x,y)=\mu_i\left(\nu_j^2 - \frac1{12}\right),~~\\
  &P_{10}(x,y)=\nu_j^3 - \frac3{20}\nu_j,~~P_{11}(x,y)=\left(\mu_i^2-\frac1{12}\right)\left(\nu_j^2 - \frac1{12}\right).
\end{split}
\end{equation}
We also define that $\overline{u}_5^n:=\overline{u}_{i,j}^n$, $\overline{v}_5^n:=\overline{v}_{i,j}^n$, $\overline{w}_5^n:=\overline{w}_{i,j}^n$, $I_{5}:=I_{i,j}$ and other $\{\overline{u}_{s}^n\}$, $\{\overline{v}_{s}^n\}$, $\{\overline{w}_{s}^n\}$, $\{I_s\}$ represent corresponding moments and Eulerian cells based on the serial numbers in \Cref{schematic_2d}. Then, the 2-D  HWENO reconstruction methods over $I_{i,j}$ are summarized as follows.

\begin{figure}[htb]
\centering
\begin{tikzpicture}

		\draw[black,thin] (0.5,0.5) node[left] {} -- (5,0.5) node[right]{};
		\draw[black,thin] (0.5,2) node[left] {} -- (5,2) node[right]{};
		\draw[black,thin] (0.5,3.5) node[left] {} -- (5,3.5) node[right]{};
    \draw[black,thin] (0.5,5) node[left] {} -- (5,5) node[right]{};

    \draw[black,thin] (0.5,0.5) node[left] {} -- (0.5,5) node[right]{};
    \draw[black,thin] (2,0.5) node[left] {} -- (2,5) node[right]{};
    \draw[black,thin] (3.5,0.5) node[left] {} -- (3.5,5) node[right]{};
    \draw[black,thin] (5,0.5) node[left] {} -- (5,5) node[right]{};

    \draw(0.5,0.5) node[above right=14pt] {$1$} -- (2,0.5) node[above right=14pt] {$2$};
    \draw(3.5,0.5) node[above right=14pt] {$3$} -- (3.5,0.5) node[above right=14pt] {};

    \draw(0.5,2) node[above right=14pt] {} -- (0.5,2) node[above right=14pt] {$4$};
    \draw(2,2) node[above right=14pt] {$5$} -- (3.5,2) node[above right=14pt] {$6$};

    \draw(0.5,3.5) node[above right=14pt] {$7$} -- (2,3.5) node[above right=14pt] {$8$};
    \draw(3.5,3.5) node[above right=14pt] {$9$} -- (3.5,3.5) node[above right=14pt] {};

    \draw(1.25,0) node[below] {$i-1$} -- (1.25,0) node[above right=14pt] {};
    \draw(1.25+1.5,0) node[below] {$i$} -- (1.25+1.5,0) node[above right=14pt] {};
    \draw(1.25+3,0) node[below] {$i+1$} -- (1.25+3,0) node[above right=14pt] {};

    \draw(0,-0.25+1.5) node[left] {$j-1$} -- (0,-0.25+1.5) node[above right=14pt] {};
    \draw(0,-0.25+3) node[left] {$j$} -- (0,-0.25+3) node[above right=14pt] {};
    \draw(0,-0.25+4.5) node[left] {$j+1$} -- (0,-0.25+4.5) node[above right=14pt] {};

\end{tikzpicture}

\caption{Stencil for the HWENO reconstructions on 2-D Cartesian mesh.}
\label{schematic_2d}
\end{figure}

\textbf{Step 1.} Reconstruct the first-order moments.

The first-order moments, $\{\overline{v}^n_5,\overline{w}^n_5\}$, can be extremely large when $u(x,y,t^n)$ is discontinuous in $I_{i,j}$ since $\overline{v}^n_5 \sim \frac{\Delta x}{12}\frac{\partial u}{\partial x}|_{(x_i,y_j)}$ and $\overline{w}^n_5 \sim \frac{\Delta y}{12}\frac{\partial u}{\partial y}|_{(x_i,y_j)}$. Hence, before recovering a cubic polynomial, we reconstruct the first-order moments with the following two goals: it will provide high-order approximations of first-order moments when $u(x,y,t^n)$ is smooth in $I_{i,j}$; when $u(x,y,t^n)$ is discontinuous in $I_{i,j}$, first-order moments will be reduced to a reasonable level.

The two first-order moments can be regarded as local indicators of the changing rates for the $x$- and $y$-dimensions. Each of them is highly independent of the other dimension. Hence, the reconstruction is performed in a dimension-by-dimension manner. Below, we take the $x$ direction as an example to illustrate the procedure for reconstructing its first-order moment.

\textbf{Step 1.1.} Compute approximations to the first-order moment and smoothness indicators from 1-D polynomial reconstructions.

\begin{enumerate}
	\item Construct a quartic polynomial $p_0(x)$, and three quadratic polynomials $\{p_k(x)\}_{k=1}^3$ satisfying
	\begin{equation}
		\begin{split}
		&\frac1{\Delta x}\int_{I^x_{i+l}}p_0(x)dx=\overline{u}_{i+l,j}^n,~~~l=-1,0,1,\\
		&\frac1{\Delta x}\int_{I^x_{i+l}}p_0(x)\left(\frac{x-x_{i+l}}{\Delta x}\right)dx=\overline{v}_{i+l,j}^n,~~~l=-1,1,\\
		\end{split}
	\end{equation}
	and
	\begin{equation}
		\begin{split}
			&\frac{1}{\Delta x}\int_{I^x_{i+l}}p_1(x)dx=\overline{u}_{i+l,j}^n~~l=-1,0,~~~\frac{1}{\Delta x}\int_{I^x_{i-1}}p_1(x)\left(\frac{x-x_{i-1}}{\Delta x}\right)dx=\overline{v}_{i-1,j}^n;\\
			&\frac{1}{\Delta x}\int_{I^x_{i+l}}p_2(x)dx=\overline{u}_{i+l,j}^n~~l=-1,0,1;\\
			&\frac{1}{\Delta x}\int_{I^x_{i+l}}p_3(x)dx=\overline{u}_{i+l,j}^n~~l=0,1,~~~\frac{1}{\Delta x}\int_{I^x_{i+1}}p_1(x)\left(\frac{x-x_{i+1}}{\Delta x}\right)dx=\overline{v}_{i+1,j}^n.\\
		\end{split}
	\end{equation}
	\item Compute the first-order moments of $\{p_k(x)\}_{k=0}^3$:
	\begin{equation}
		\begin{split}
			&\widetilde{\overline{v}}_{i,j,0}^n := \frac1{\Delta x}\int_{I^x_i}p_0(x)\left(\frac{x-x_i}{\Delta x}\right)dx = -\frac{5}{76}\overline{u}^n_{i-1,j}-\frac{11}{38}\overline{v}^n_{i-1,j}-\frac{11}{38}\overline{v}^n_{i+1,j}+\frac{5}{76}\overline{u}^n_{i+1,j},\\
			&\widetilde{\overline{v}}_{i,j,1}^n := \frac1{\Delta x}\int_{I^x_i}p_1(x)\left(\frac{x-x_i}{\Delta x}\right)dx=\frac16\overline{u}^n_{i,j}-\frac16\overline{u}^n_{i-1,j}-\overline{v}^n_{i-1,j},\\
			&\widetilde{\overline{v}}_{i,j,2}^n := \frac1{\Delta x}\int_{I^x_i}p_2(x)\left(\frac{x-x_i}{\Delta x}\right)dx=\frac1{24}\overline{u}^n_{i+1,j}-\frac1{24}\overline{u}^n_{i-1,j},\\
			&\widetilde{\overline{v}}_{i,j,3}^n := \frac1{\Delta x}\int_{I^x_i}p_2(x)\left(\frac{x-x_i}{\Delta x}\right)dx=\frac16\overline{u}^n_{i+1,j}-\frac16\overline{u}^n_{i,j}-\overline{v}^n_{i+1,j}.\\
		\end{split}
	\end{equation}

	\item Compute the smoothness indicators $\{\beta_k\}_{k=0}^3$ of $\{p_k(x)\}_{k=0}^3$ \cite{liu1994weighted,jiang1996efficient,shu1998essentially}:
	\begin{equation}
		\beta_k=\sum_{l=1}^r\frac1{\Delta x}\int_{I_i^x}\left(\Delta x^l\frac{\partial^l}{\partial x^l}p_k(x)\right)^2dx~~~k=0,1,2,3
	\end{equation}
	with $r$ representing the degree of the corresponding polynomial. Here, the smoothness indicators $\{\beta_k\}_{k=1}^3$ can be explicitly expressed by
	\begin{equation}
		\begin{split}
			&\beta_1 = \left(12\widetilde{\overline{v}}_{i,j,1}^n\right)^2+156\left(\widetilde{\overline{v}}_{i,j,1}^n-\overline{v}^n_{i-1,j}\right)^2,\\
			&\beta_2 = \left(12\widetilde{\overline{v}}_{i,j,2}^n\right)^2+\frac{13}{12}\left(\overline{u}^n_{i+1,j}-2\overline{u}^n_{i,j}+\overline{u}^n_{i-1,j}\right)^2,\\
			&\beta_3 = \left(12\widetilde{\overline{v}}_{i,j,3}^n\right)^2+156\left(\overline{v}^n_{i+1,j}-\widetilde{\overline{v}}_{i,j,3}^n\right)^2.
		\end{split}
	\end{equation}
	We refer to \cite{zhao2020hermite} for the explicit expression of $\beta_0$.

	\item Compute a full stencil global reference smoothness indicator \cite{Wibisono2021HTENO}:
	\begin{equation}
		\tau := \left(\frac{|\beta_0-\beta_1|+|\beta_0-\beta_3|}{2}\right)^2.
	\end{equation}
	By the Taylor expansion, we can find that $\tau = O(\Delta x^6)$, if there is no discontinuity involved.

\end{enumerate}

\textbf{Step 1.2} Weight the collected first-order moments.

Below, we first introduce the weighting strategy for the HWENO-1.
\begin{enumerate}
	\item Compute the nonlinear weights by
	\begin{equation}
		\omega_k = \frac{\overline{\omega}_k}{\sum_l{\overline{\omega}_l}}~~~ \text{with}~~~ \overline{\omega}_k = \gamma_k\left(1+\frac{\tau}{\beta_k+\epsilon}\right)~~~ k = 0,1,3,
	\end{equation}
	where $\epsilon = 10^{-40}$ is set to avoid the denominator being zero. The linear weights, $\{\gamma_0, \gamma_1, \gamma_3\}$, are chosen as $\{0.6,0.2,0.2\}$ in this paper. We refer to \cite{ZHU2016110,zhao2020hermite} for more details on linear and nonlinear weights.
	
	\item Reconstruct the $x$-dimension first-order moment, $\overline{v}_{5}^n$, by:
	\begin{equation}
		\widetilde{\overline{v}}_{5}^n = \frac{\omega_0}{\gamma_0}\left(\widetilde{\overline{v}}_{i,j,0}^n - \gamma_1\widetilde{\overline{v}}_{i,j,1}^n - \gamma_3\widetilde{\overline{v}}_{i,j,3}^n\right) + \omega_1\widetilde{\overline{v}}_{i,j,1}^n + \omega_3\widetilde{\overline{v}}_{i,j,3}^n.
	\end{equation}
	
\end{enumerate}

The HWENO-2 reconstructs the first-order moment $\overline{v}^n_5$ as follows.

\begin{enumerate}
	\item Separate the discontinuities from broad-band smooth fluctuations as illustrated in \cite{FU2016333,FU2019117} by first taking
\begin{equation}
	\eta_k = \left( 1 + \frac{\tau}{\beta_k+\epsilon} \right)^6~~~k = 1,2,3,
\end{equation}
where $\epsilon = 10^{-40} $ is used to avoid the denominator being zero as in \cite{BORGES20083191}. If these is no discontinuity involved, we can find that $\eta_k \approx 1$ for all $k$. If there is discontinuity involved for the global three-cells stencil, the $\eta$ value for a smooth small stencil can be greatly enlarge with a magnitude of $O(\Delta x^{-12})$. Then, we normalize $\{\eta_k\}_{k=1}^3$ with
\begin{equation}
	\kappa_k = \frac{\eta_k}{\sum_{l=1}^3\eta_l}~~~k = 1,2,3.
\end{equation}

\item Reconstruct the $x$-dimension first-order moment $\overline{v}_{5}^n$ by:
\begin{equation}
	\widetilde{\overline{v}}_{5}^n =\widetilde{\overline{v}}_{i,j,0}^n,~~~\text{if}~~~\min_{k}\kappa_k > C_T,
\end{equation}
otherwise
\begin{equation}
	\widetilde{\overline{v}}_{5}^n = \begin{cases}
		\begin{split}
			&\widetilde{\overline{v}}_{i,j,1}^n,~~~~~~~~~~~~\text{if}~~~\kappa_1>\kappa_3,\\
			&\widetilde{\overline{v}}_{i,j,3}^n,~~~~~~~~~~~~\text{if}~~~\kappa_3>\kappa_1.
		\end{split}
	\end{cases}
\end{equation}
where $C_T$ is a parameter deciding whether a corresponding polynomial should be involved \cite{FU2016333}. We empirically choose $C_T = 10^{-3}$ as in \cite{Wibisono2021HTENO}. For treating discontinuity, unlike \cite{FU2016333,Wibisono2021HTENO}, we only choose the smoothest polynomial for reconstruction, since the weighted summation does not increase the order of accuracy for our case. When $\{\kappa_k\}$ matches the smoothness of the three-cells stencil, we directly use $\widetilde{\overline{v}}_{i,j,0}^n$ for optimal accuracy.
\end{enumerate}

The first-order moment in $y$ direction can be reconstructed in a similar way; we the reconstructed  first-order moment in this direction by $\widetilde{\overline{w}}_{5}^n$.

\textbf{Step 2.} Recover the reconstructed $H^n(x,y)$ on $I_{i,j}$.

\textbf{Step 2.1.} Collect information from different 2-D polynomials.
\begin{enumerate}

\item Construct a polynomial $\widetilde{q}_0(x,y):=\sum_{l=1}^{11}a^{q_0}_lP_l(x,y)$ such that
\begin{equation}
	\begin{split}
		&\frac{1}{\Delta x\Delta y}\iint_{I_s}\widetilde{q}_0(x,y)dxdy = \overline{u}_s^n,~~~s=1,2,\ldots,9,\\
		&\frac{1}{\Delta x\Delta y}\iint_{I_{5}}\widetilde{q}_0(x,y)\left(\frac{x-x_i}{\Delta x}\right)dxdy = \widetilde{\overline{v}}_5^n,\\
		&\frac{1}{\Delta x\Delta y}\iint_{I_{5}}\widetilde{q}_0(x,y)\left(\frac{y-y_j}{\Delta y}\right)dxdy = \widetilde{\overline{w}}_5^n.\\
	\end{split}
\end{equation}
Let $q_0(x,y) = \sum_{l=1}^{10}a^{q_0}_lP_l(x,y)$, which is the orthogonal projection of $\widetilde{q}_0(x,y)$ to $P^3(I_{i,j})$. We provide the explicit expressions of $\{a^{q_0}_{l}\}_{l=1}^{10}$ in \Cref{appendix:coefficients}.

\item Construct four quadratic polynomial $\{q_k(x,y)\}_{k=1}^4:=\{\sum_{l=1}^6a^{q_k}_lP_l(x,y)\}_{k=1}^4$ satisfying
\begin{equation}
\begin{split}
&\frac{1}{\Delta x\Delta y}\iint_{I_s}q_k(x,y)dxdy = \overline{u}^n_s,\\
&\frac{1}{\Delta x\Delta y}\iint_{I_5}q_k(x,y)\frac{x-x_i}{\Delta x}dxdy = \widetilde{\overline{v}}^n_5,\\
&\frac{1}{\Delta x\Delta y}\iint_{I_5}q_k(x,y)\frac{y-y_i}{\Delta y}dxdy = \widetilde{\overline{w}}^n_5,
\end{split}
\end{equation}
where
\begin{equation}
\begin{split}
&s=1,2,4,5,~~\text{for}~~ k=1;~~s=2,3,5,6,~~\text{for}~~ k=2;\\
&s=4,5,7,8,~~\text{for}~~ k=3;~~s=5,6,8,9,~~\text{for}~~ k=4.
\end{split}
\end{equation}
The explicit expressions of $\{a^{q_k}_l\}$ are given in \Cref{appendix:coefficients}.

\item Compute the smoothness indicator \cite{liu1994weighted,jiang1996efficient,shu1998essentially} $\{\hat{\beta}_k\}_{k=0}^4$ of $\{q_k(x,y)\}_{k=0}^4$:
\begin{equation}\label{detector_beta}
	\begin{split}
		\hat{\beta}_0 = &\frac1{\Delta x\Delta y}\sum_{l_1+l_2<=3}\iint_{I_{i,j}}\left(\Delta x^{l_1}\Delta y^{l_2}\frac{\partial^{|l_1+l_2|}}{\partial x^{l_1}\partial y^{l_2}}q_0(x,y)\right)^2dxdy;\\
		\hat{\beta}_k=&\frac1{\Delta x\Delta y}\sum_{l_1+l_2<=2}\iint_{I_{i,j}}\left(\Delta x^{l_1}\Delta y^{l_2}\frac{\partial^{|l_1+l_2|}}{\partial x^{l_1}\partial y^{l_2}}q_k(x,y)\right)^2dxdy,~~k=1,2,3,4.
	\end{split}
\end{equation}
The explicit expression of $\{\hat{\beta}_k\}_{k=0}^4$ can be given by
\begin{equation}
	\begin{split}
		\hat{\beta}_0 = &\left(a^{q_0}_2+\frac1{10}a^{q_0}_7\right)^2+\left(a^{q_0}_3+\frac1{10}a^{q_0}_{10}\right)^2+\frac{13}3\left(a^{q_0}_4\right)^2+\frac76\left(a^{q_0}_5\right)^2+\frac{13}3\left(a^{q_0}_6\right)^2\\
		&+\frac{781}{20}\left(a^{q_0}_7\right)^2+\frac{47}{10}\left(a^{q_0}_8\right)^2+\frac{47}{10}\left(a^{q_0}_9\right)^2+\frac{781}{20}\left(a^{q_0}_{10}\right)^2;\\
		\hat{\beta}_k =&\left(a^{q_k}_2\right)^2+\left(a^{q_k}_3\right)^2+\frac{13}3\left(a^{q_k}_4\right)^2+\frac76\left(a^{q_k}_5\right)^2+\frac{13}3\left(a^{q_k}_6\right)^2~~~\text{for}~~k=1,2,3,4.
	\end{split}
\end{equation}

\item Compute a full stencil global reference smoothness indicator:
\begin{equation}
	\hat{\tau}=\left(\frac{|\hat{\beta}_0-\hat{\beta_1}|+|\hat{\beta}_0-\hat{\beta_2}|+|\hat{\beta}_0-\hat{\beta_3}|+|\hat{\beta}_0-\hat{\beta_4}|}{4}\right)^2.
\end{equation}
Similarly, by the Taylor expansion, we can easily check that $\hat{\tau}=O(\Delta x^6)$, if there is no discontinuity involved.

\end{enumerate}

\textbf{Step 2.2} Weight the collected 2-D polynomials.

The HWENO-1 weights the polynomials as follows.
\begin{enumerate}
	\item Compute the nonlinear weights by
\begin{equation}
	\hat{\omega}_k = \frac{\widetilde{\omega}_k}{\sum_l{\widetilde{\omega}_l}}~~~ \text{with}~~~ \widetilde{\omega}_k = \gamma_k\left(1+\frac{\tau}{\beta_k+\epsilon}\right)~~~ k = 0,1,2,3,4,
\end{equation}
where $\epsilon = 10^{-40}$ is set to avoid the denominator being zero. The linear weights, $\{\gamma_k\}_{k=0}^4$ are chosen as $\{0.6, 0.1, 0.1, 0.1, 0.1\}$ in this paper.

	\item Construct $H^{(i,j)}(x,y):=\sum_{l=1}^{10}a_lP_l(x,y)$:
\begin{equation}
	\begin{split}
		H^{(i,j)}(x,y)=\frac{\hat{\omega}_0}{\hat{\gamma}_0}\left(q_0(x)-\sum_{k=1}^4\hat{\gamma}_kq_k(x,y)\right)+\sum_{k=1}^4\hat{\omega}_kq_k(x,y),~~~(x,y)\in I_{i,j}.
	\end{split}
\end{equation}
Here, the coefficients $\{a_l\}$ can be explicitly given by
\begin{equation}
	\begin{split}
		&a_1 = \overline{u}^n_5;~~ a_2 = 12\widetilde{\overline{v}}^n_5;~~ a_3 =  12\widetilde{\overline{w}}^n_5;\\
		&a_l = \frac{\hat{\omega}_0}{\hat{\gamma}_0}a^{q_0}_l + \sum_{k=1}^4\left(\hat{\omega}_k-\frac{\hat{\omega}_0}{\hat{\gamma}_0}\hat{\gamma}_k\right)a^{q_k}_l~~ \text{for}~~ l= 4,5,6;\\
		&a_l = \frac{\hat{\omega}_0}{\hat{\gamma}_0}a^{q_0}_l~~\text{for}~~l=7,8,\cdots,10.
	\end{split}
\end{equation}
\end{enumerate}

The HWENO-2 weights the polynomials as follows.
\begin{enumerate}
\item Compute the separation parameters $\{\hat{\eta}_k\}_{k=1}^4$ and the corresponding normalized parameters $\{\hat{\kappa}_k\}_{k=1}^4$ similarly:
\begin{equation}
	\hat{\kappa}_k=\frac{\hat{\eta}_k}{\sum_{l=1}^3\hat{\eta}_l}~~~\text{with $\hat{\eta}_k=\left(1+\frac{\hat{\tau}_6}{\hat{\beta}_k+\epsilon}\right)^6$ for $k=1,2,3,4,$}
\end{equation}
where $\epsilon = 10^{-40}$.

\item Construct $H^{(i,j)}(x,y)$ by
\begin{equation}
	\begin{split}
		H^{(i,j)}(x,y)=q_0(x,y),~~~\text{if}~~~\min_k\hat{\kappa}_k > C_T,
	\end{split}
\end{equation}
otherwise
\begin{equation}
	\begin{split}
		H^{(i,j)}(x,y)=q_K(x,y)~~~\text{with $K$ being the index such that $\hat{\kappa}_K=\max_k\hat{\kappa}_k$}.
	\end{split}
\end{equation}
Here, $C_T$ is chosen as $10^{-3}$.
\end{enumerate}

In particular, if $\widetilde{\overline{v}}^n_{i,j} = \widetilde{\overline{v}}^n_{i,j,0}$, $\widetilde{\overline{w}}^n_{i,j} = \widetilde{\overline{w}}^n_{i,j,0}$ (with similar notation), and $H^{(i,j)}(x,y) = q_0(x,y)\quad \forall i,j,$ we call such a reconstruction linear reconstruction.


\subsection{Constructing cubic-curved upstream cells and approximating $w(x,y,t^n)$}\label{cubic_upstream}

To evaluate the right-hand side of \eqref{SL_formulation}, it remains to provide the approximations of $\{I_{i,j}^\star\}$ and $w(x,y,t^n)$. For a given index $(i,j)$, the cubic-curved quadrilateral upstream cell, $\widetilde{I}^\star_{i,j}$, and the cubic polynomial $\widetilde{w}(x,y)$ on $\widetilde{I}^\star_{i,j}$ are constructed by the following procedure.

\begin{enumerate}
\item Tracing characteristics backward in time.

We locate $4\times 4$ Gauss-Legendre-Lobatto (GLL) points on $I_{i,j}$. We determine the characteristic feet of these GLL points by solving \eqref{eq:ODEs} at $t = t^n$ (see \Cref{charac_trace}). In practice, we solve the ODEs \eqref{eq:ODEs} by a fourth-order Runge-Kutta (RK) method.

\begin{figure}[htb]
  \centering
  \subfloat[]{
  \begin{minipage}[t]{7cm}
\begin{overpic}[scale=0.07]{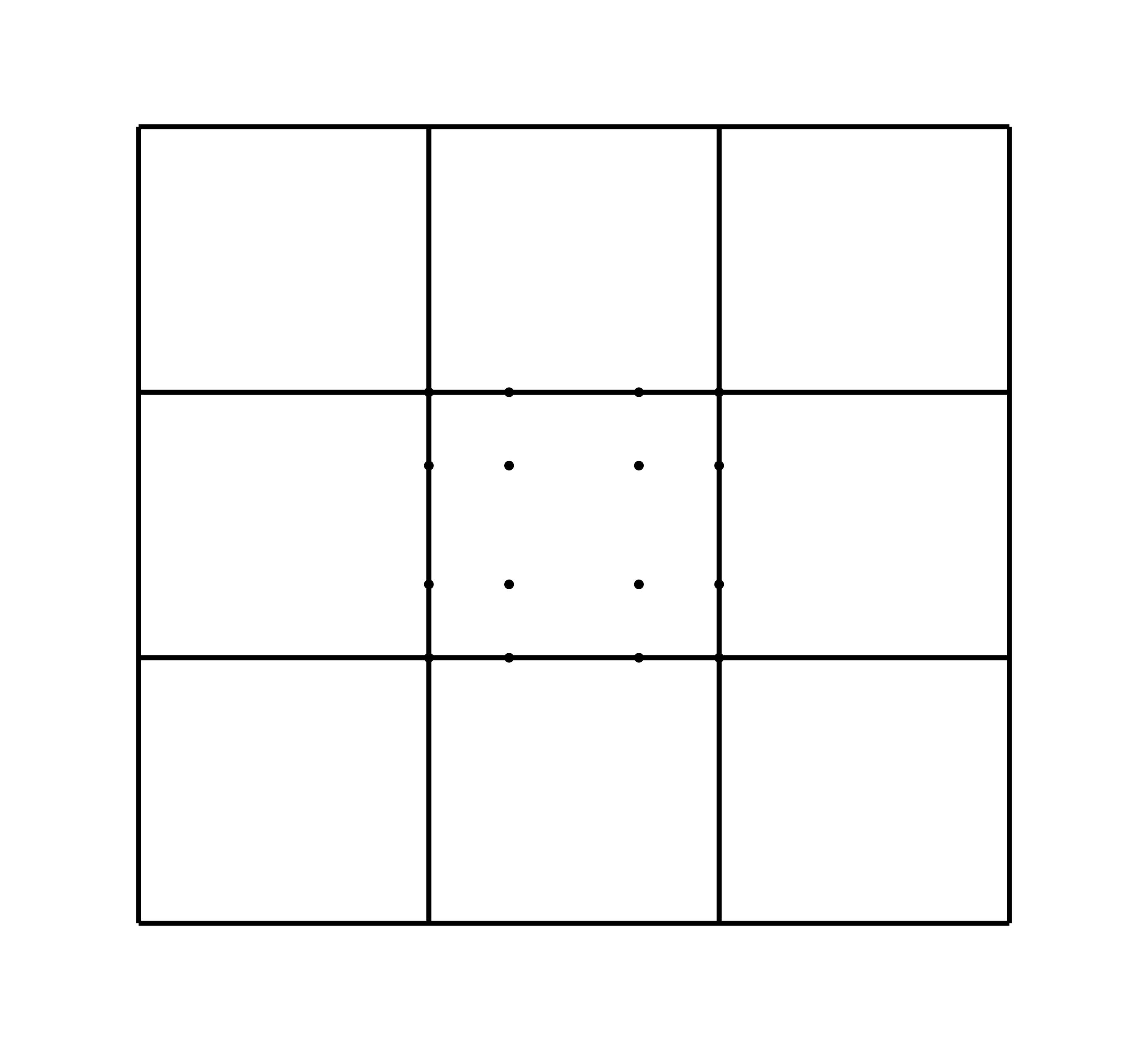}
\put(49,4){$i$}
\put(4,44.5){$j$}
\end{overpic}
  \end{minipage}
  }
  \subfloat[]{
  \begin{minipage}[t]{7cm}
\begin{overpic}[scale=0.07]{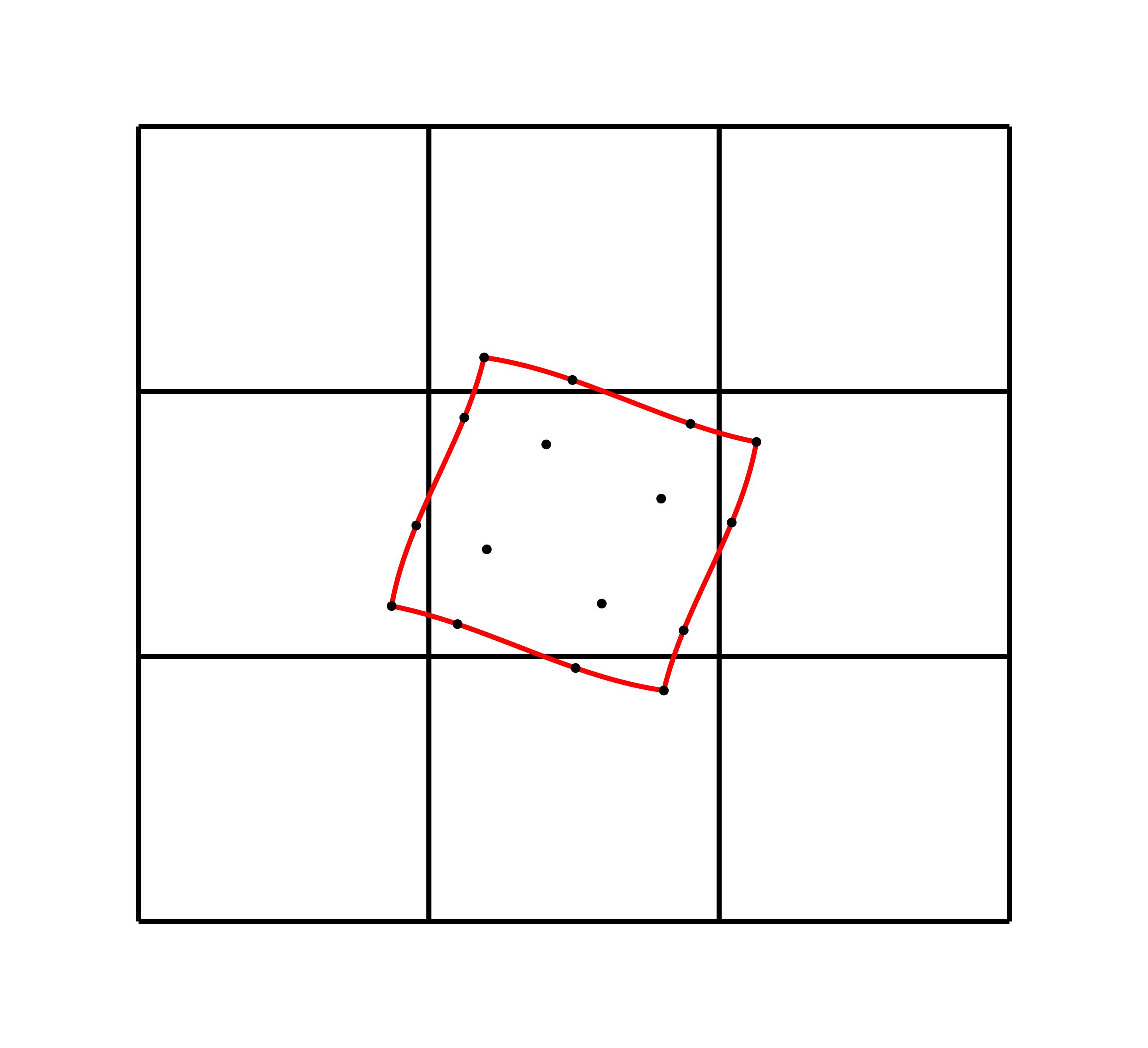}
\put(49,4){$i$}
\put(4,44.5){$j$}
\end{overpic}
  \end{minipage}
  }
  \caption{Left: the black solid lines represent the Eulerian mesh; the black dots are the GLL points located on the Eulerian cell $I_{i,j}$. Right: the black solid lines represent the Eulerian mesh; the red solid lines represent the boundary of $I_{i,j}^\star$; the black dots are the characteristic feet obtained by solving \eqref{eq:ODEs}.}
  \label{charac_trace}
\end{figure}

\item Constructing the edges of the cubic-curved quadrilateral upstream cells.

The edges of the cubic-curved quadrilateral upstream cells are constructed by a cubic interpolation procedure (see \cite{ZHENG2022114973}). In particular, we prefer to use a parametric form to present any cubic-curved edge:
\begin{equation}\label{parame_cubic}
     \begin{cases}
     x(\xi) = x_a\xi^3+x_b\xi^2+x_c\xi+x_d, \\
     y(\xi) = y_a\xi^3+y_b\xi^2+y_c\xi+y_d,~~~\xi\in[-1,1].
     \end{cases}
\end{equation}

\item Constructing $\widetilde{w}(x,y)$.

We denote the sixteen GLL points in \Cref{charac_trace} (a) by $\{v_k\}$ and denote the corresponding characteristic feet by $\{v_k^\star\}$. By the adjoint problem \eqref{eq:adjoint}, we have
\begin{equation}
	w(x(v_k^\star),y(v_k^\star))=W(x(v_k),y(v_k))~~~\text{for}~~k=1,2,\ldots,16.
\end{equation}

Hence, by a standard least square procedure, we can find a cubic polynomial $\widetilde{w}(x,y)\in P^3$.

\end{enumerate}

\subsection{Numerical integration}\label{sec:numerical_integration}

\begin{figure}[htb]
	\centering
	\subfloat[]{
		\begin{minipage}[t]{7cm}
			\begin{overpic}[scale=0.35]{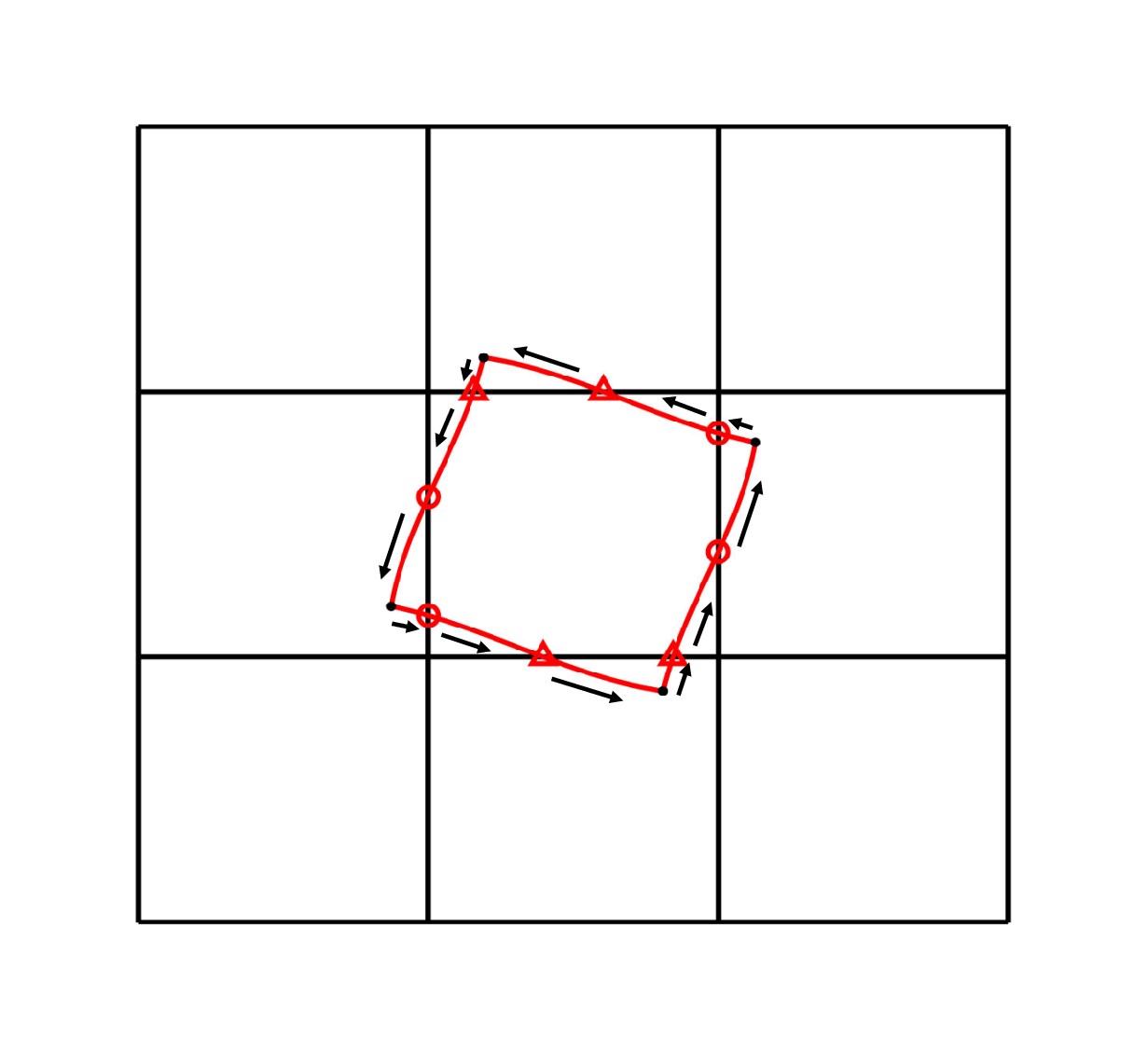}
				\put(48,44){$\widetilde{I}_{i,j}^\star$}
				\put(68,47){$\{\mathcal{L}_{i,j;p,q}^{k}\}$}
			\end{overpic}
		\end{minipage}
	}
	\subfloat[]{
		\begin{minipage}[t]{7cm}
			\begin{overpic}[scale=0.35]{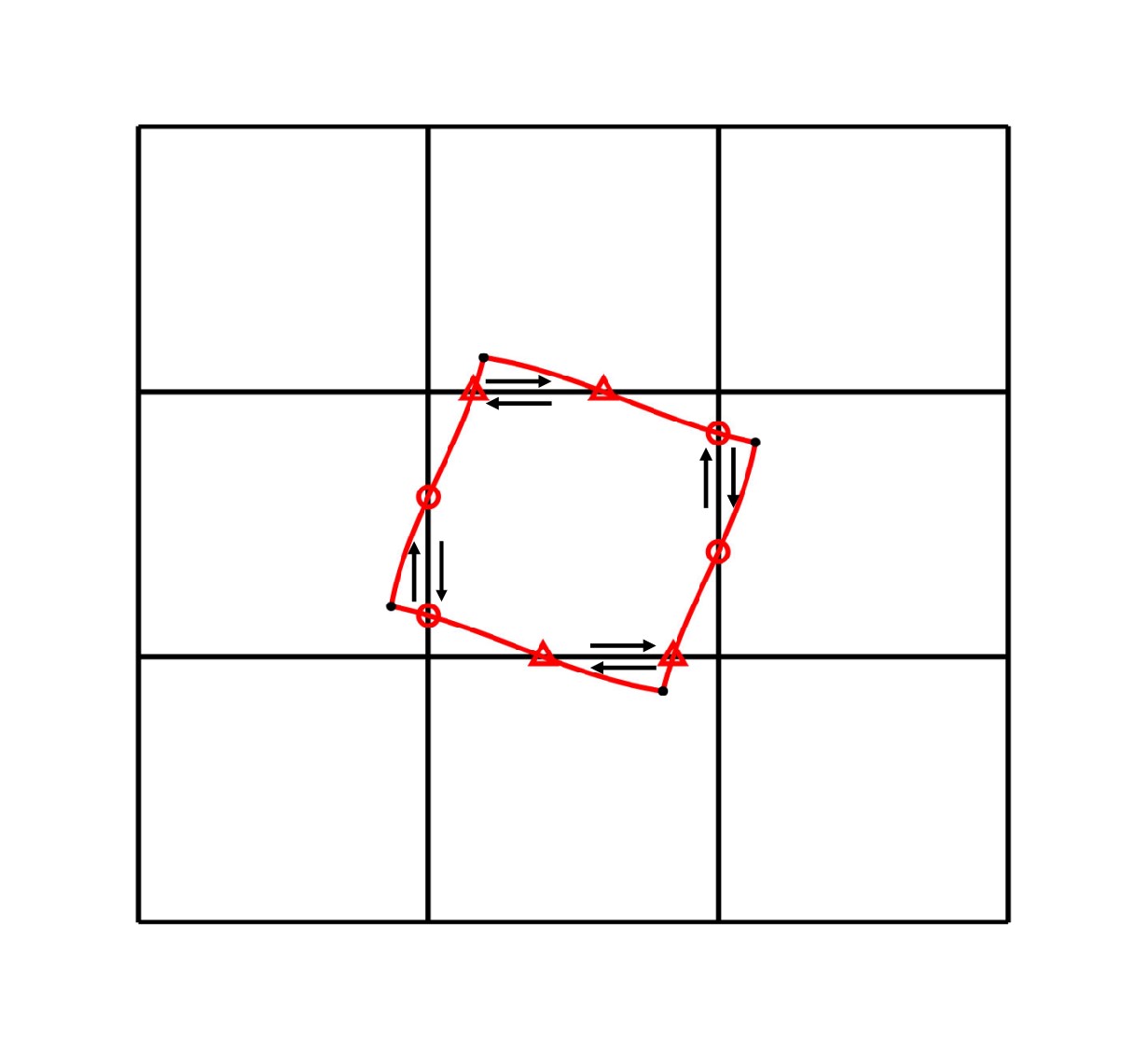}
				\put(48,44){$\widetilde{I}_{i,j}^\star$}
				\put(68,47){$\{\mathcal{S}_{i,j;p,q}^{k}\}$}
			\end{overpic}
		\end{minipage}
	}
	\caption{Schematic illustration for the outer integral segments (a) and the inner integral segments (b). The red circles and triangles are the intersections of $\widetilde{I}_{i,j}^\star$ and the Eulerian mesh.}
	\label{fig:innerseg_outerseg_def}
\end{figure}

For numerically evaluating the right-hand side of \eqref{SL_formulation}, we integrate the piecewise polynomial $H^n(x,y)\widetilde{w}(x,y)$ over each $\widetilde{I}_{i,j}^\star$, which may cross different Eulerian background cells. Hence, $\widetilde{I}_{i,j}^\star$ is first clipped into curved polygons such that $\widetilde{I}_{i,j}^\star=\cup_{(p,q)}(\widetilde{I}^\star_{i,j}\cap I_{p,q})$ and the integrand is smooth in each polygon. For conciseness, the curved polygons $\left\{\widetilde{I}^\star_{i,j}\cap I_{p,q}\right\}$  are denoted by $\left\{\widetilde{I}_{i,j;p,q}^\star\right\}$. In particular, we are concerned about the information along the edges of $\left\{\widetilde{I}_{i,j;p,q}^\star\right\}$. The edges of $\left\{\widetilde{I}_{i,j;p,q}^\star\right\}$, overlapping $\partial \widetilde{I}_{i,j}^\star$, with counterclockwise direction with respect to $\widetilde{I}_{i,j}^\star$ are denoted by $\{\mathcal{L}_{i,j;p,q}^{k}\}$ (see \Cref{fig:innerseg_outerseg_def} (a)). We call $\{\mathcal{L}_{i,j;p,q}^{k}\}$ the {\em outer integral segments} of the upstream cell $\widetilde{I}_{i,j}^\star$. Similarly, the edges of $\left\{\widetilde{I}_{i,j;p,q}^\star\right\}$, overlapping mesh lines,  with counterclockwise direction with respect to the corresponding curved polygons are denoted by $\{\mathcal{S}_{i,j;p,q}^{k}\}$ (see \Cref{fig:innerseg_outerseg_def} (b)). We call $\{\mathcal{S}_{i,j;p,q}^{k}\}$ the {\em inner integral segments} of $\widetilde{I}_{i,j}^\star$. For implementation, we refer to \cite{ZHENG2022114973} for more details.

With the clipped outer integral segments, $\{\mathcal{L}_{i,j;p,q}^k\}$, as well as the inner integral segments, $\{\mathcal{S}_{i,j;p,q}^k\}$, we evaluate the right-hand side of \eqref{SL_formulation} as follows
\begin{equation}
\begin{split}
&\frac1{\Delta x \Delta y}\iint_{I_{i,j}^\star} u(x,y,t^n)w(x,y)dxdy\\
\approx&\frac1{\Delta x \Delta y}\iint_{\widetilde{I}_{i,j}^\star} H^n(x,y)\widetilde{w}(x,y)dxdy\\
=&\frac1{\Delta x \Delta y}\sum_{(p,q)}\iint_{\widetilde{I}_{i,j;p,q}^\star}H^{(p,q)}(x,y)\widetilde{w}(x,y)dxdy\\
=&  \frac1{\Delta x \Delta y}\sum_{(p,q)}\int_{\partial(\widetilde{I}_{i,j;p,q}^\star)}\left[Pdx+Qdy\right]\\
=&\frac1{\Delta x \Delta y}\sum_{(p,q)}\Big\{\sum_{k}\int_{\mathcal{L}_{i,j;p,q}^k}\left[Pdx+Qdy\right]+\sum_{k}\int_{\mathcal{S}_{i,j;p,q}^k}\left[Pdx+Qdy\right]\Big\},
\end{split}
\label{SL_formulation_numer}
\end{equation}
where $P(x,y)$ and $Q(x,y)$ are piecewise smooth auxiliary functions such that
\begin{equation}
-\frac{\partial P}{\partial y} + \frac{\partial Q}{\partial x} = H^{(p,q)}(x,y)\widetilde{w}(x,y).
\end{equation}
It is straightforward for evaluating the line integral on $\{\mathcal{S}_{i,j;p,q}^k\}$. For the integral on an outer integral segment, say $\mathcal{L}_{i,j;p,q}^k$,
\begin{equation}
\begin{split}
&\int_{\mathcal{L}_{i,j;p,q}^k}[Pdx+Qdy]=\int_{\xi_k}^{\xi_{k+1}}\left[P\left(x(\xi),y(\xi)\right)x^\prime(\xi)+Q\left(x(\xi),y(\xi)\right)\nu_q^\prime(\xi)\right]d\xi,
\end{split}
\end{equation}
where $\xi_k$ and $\xi_{k+1}$ represent the $\xi$ value of the start point and end point of $\mathcal{L}_{i,j;p,q}^k$ \eqref{parame_cubic}.

The proposed SL HWENO schemes also equip a PP limiter \cite{Zhang2010On1}, when the analytical solution of \eqref{2d_linear} stays positive. For the implementation of this PP limiter, we refer to \cite{ZHENG2022114973} for a detailed description.

\begin{rem}\label{prop:stability}
	We can numerically prove that the numerical update provided by \eqref{SL_formulation_numer} is unconditionally stable for linear transport equations with constant coefficients and periodic boundary condition if $H^n(x,y)$ is reconstructed by the linear reconstruction defined in \Cref{2dHWENO}. The proof is accomplished by the standard von Neumann analysis. We arrange this proof in \Cref{appendix:stability}.
\end{rem}

\section{SL HWENO schemes for nonlinear models}\label{sec:nonlinear}

 The non-splitting SL HWENO schemes are coupled with a fourth-order RKEI in the same framework as in \cite{CAI2021110036} to solve nonlinear models such as the Vlasov-Poisson system, the guiding center Vlasov model, and the incompressible Euler equations in the vorticity-stream function. Below, we briefly describe these three models.

Arising from collisionless plasma, the Vlasov-Poisson system reads
\begin{equation}\label{eq:VP_1}
	f_t+vf_x+E(x,t)f_v=0,
\end{equation}
\begin{equation}\label{eq:VP_2}
	E(x,t)=-\phi_x,~~~-\phi_{xx}(x,t)=\rho(x,t),
\end{equation}
where $x$ represents the spatial position, $v$ is the velocity, $f(x,v,t)$ is the probability distribution function describing the probability of a particle arises at position $x$ with velocity $v$ at time $t$. The electric field $E$ is determined by the Poisson's equation \eqref{eq:VP_2}. $\phi$ is the self-consistent electrostatic potential. $\rho=\int_{\mathbb{R}}f(x,v,t)dv-\rho_0$ is the charge density with $\rho_0 = \frac{1}{|\Omega_x|}\int_{\Omega_x}\int_{\mathbb{R}}f(x,v,0)dvdx$.

The 2-D guiding center Vlasov model describes highly magnetized plasma in the transverse plane of a tokamak \cite{YANG201418,frenod2015}, which can be written as:
\begin{equation}\label{eq:GC_1}
\rho_t + \nabla\cdot\left(\mathbf{E}^{\bot}\rho\right)=0,
\end{equation}
\begin{equation}\label{eq:GC_2}
-\Delta \Phi = \rho,~~\mathbf{E}^{\bot}=(-\Phi_y,\Phi_x),
\end{equation}
where $\rho(x,y,t)$ is the charge density and $\mathbf{E}$ is the electric field depends on $\rho$ via the Poisson's equation \eqref{eq:GC_2}.

The 2-D incompressible Euler equation in vorticity-stream function formulation reads
\begin{equation}\label{eq:IE_1}
	\omega_t + \nabla\cdot(\mathbf{u}\omega)=0,
\end{equation}
\begin{equation}\label{eq:IE_2}
	\Delta\psi=\omega,~~\mathbf{u}=(-\psi_y,\psi_x),
\end{equation}
where $\omega(x,y,t)$ is the vorticity of the fulid, $\mathbf{u}:=(u_1,u_2)$ is the velocity field, and $\psi$ is the stream-function determined by Poisson’s equation \eqref{eq:IE_2}.

Notice that these three models can be written in the form of
\begin{equation}
	u_t + \nabla\cdot\left(\mathbf{V}(u(\mathbf{x},t))u\right)=0,
\end{equation}
where $\mathbf{V}(u(\mathbf{x},t))$ represents the velocity field. We briefly summarize the SL HWENO scheme coupled with the fourth-order RKEI as follows.

 \begin{equation}\label{eq:CF4_SL_FV_WENO}
 \begin{split}
 \overline{u}^{(1)} &= \overline{u}^n\\
 \overline{u}^{(2)} &= SLHWENO\left(\frac12\mathbf{V}\left(\overline{u}^{(1)}\right),\Delta t\right)\overline{u}^n\\
 \overline{u}^{(3)} &= SLHWENO\left(\frac12\mathbf{V}\left(\overline{u}^{(2)}\right),\Delta t\right)\overline{u}^n\\
 \overline{u}^{(4)} &= SLHWENO\left(-\frac12\mathbf{V}\left(\overline{u}^{(1)}\right)+\mathbf{V}\left(\overline{u}^{(3)}\right),\Delta t\right)\overline{u}^{(2)}\\
 \overline{u}^{n+1} &= SLHWENO\left(-\frac1{12}\mathbf{V}\left(\overline{u}^{(1)}\right)+\frac16\mathbf{V}\left(\overline{u}^{(2)}\right)+\frac16\mathbf{V}\left(\overline{u}^{(3)}\right)+\frac14\mathbf{V}\left(\overline{u}^{(4)}\right),\Delta t\right)\\
 &SLHWENO\left(\frac1{4}\mathbf{V}\left(\overline{u}^{(1)}\right)+\frac16\mathbf{V}\left(\overline{u}^{(2)}\right)+\frac16\mathbf{V}\left(\overline{u}^{(3)}\right)-\frac1{12}\mathbf{V}\left(\overline{u}^{(4)}\right),\Delta t\right)\overline{u}^n,
 \end{split}
 \end{equation}
 where $SLHWENO\left(\mathbf{V}\left(\overline{u}^{(k)}\right),\frac12\Delta t\right)\overline{u}^{(l)}$ represents the solution evolved from $\overline{u}^{(l)}$ with time step $\frac12\Delta t$ and velocity field $\mathbf{V}\left(\overline{u}^{(k)}\right)$ by a non-splitting SL HWENO scheme. For approximating the velocity field, we use Fast Fourier transform   to solve the Poisson's equations.

\section{Numerical tests}\label{sec:numerical_tests}

\subsection{Linear transport equations}

In this subsection, we test two benchmark problems: transport equation with constant coefficients and the swirling deformation flow.
As a comparison, we test the HWENO reconstruction in \cite{zhao2020hermite} with their $p_0(x,y)\in P^4$ replaced by the $q_0(x,y) \in P^3$ in this paper. Then, the corresponding HWENO reconstruction is denoted by HWENO-old. Through out the tests below, we aim to present the performance of the SL HWENO schemes and compare the differences between different HWENO reconstructions.

Unless specified, we set $\Delta t = \frac{\text{CFL}}{\frac{\text{max}\{|a(x,y,t)|\}}{\Delta x}+\frac{\text{max}\{|b(x,y,t)|\}}{\Delta y}}$ and CFL = 10.2. The PP limiter is applied for the problems with non-negative initial conditions.

\begin{exa} (Transport equation with constant coefficients).
 Consider
\begin{equation}\label{2-D_linear_transport}
  u_t+u_x+u_y=0,\quad x\in[-\pi,\pi],\quad y\in[-\pi,\pi],
\end{equation}
with a smooth initial condition, $u(x,y,0) = \text{sin}(10(x+y))$, and the periodic boundary condition. The exact solution for this problem is $u(x,y,t) = \text{sin}(10(x+y-2t))$. In \Cref{table:order_const}, we present the $L^2$ erors and corresponding orders of accuracy for the SL HWENO-1, HWENO-2, and HWENO-old schemes. As shown, all the schemes converge to fourth-order accuracy. However, the SL HWENO-old scheme requires a denser mesh to obtain a normal accuracy. In \Cref{fig:slide_const}, we show the cross-section of the numerical solution at $x=y$ of the three schemes at $T=20$ with a fixed mesh of $80\times80$. With four points per wavelength at $x=y$, we observe that HWENO-old tends to be more smeared than the other two.
\begin{table}[!htbp]
\centering
\caption{ (Transport equation with constant coefficients). $L^2$ errors and corresponding orders of accuracy of the SL HWENO schemes for (\ref{2-D_linear_transport}) with $u(x,y,0) = \text{sin}(10(x+y))$ at $T = 20$.}\label{table:order_const}
  \centering
\begin{tabular}{|c|cc|cc|cc|}
\hline
&\multicolumn{2}{|l|}{\textbf{HWENO-1}}&\multicolumn{2}{|l|}{\textbf{HWENO-2}}&\multicolumn{2}{|l|}{\textbf{HWENO-old}}\\
\cline{2-7}
mesh&$L^2$ error&order&$L^2$ error&order&$L^2$ error&order\\
\hline
  20$\times$  20    &   7.07E-01&    ---    &   7.07E-01&    ---    &   7.07E-01&   --- \\
  40$\times$  40	&	1.59E-01&   2.15    &	2.35E-01&   1.59    &	6.99E-01&   0.02 \\
  80$\times$  80	&	1.71E-02&   3.78    &	1.71E-02&   3.78	&	1.54E-01&   2.18 \\
  160$\times$  160	&	1.13E-03&   3.93    &	1.13E-03&   3.93 	&	1.61E-03&   6.58 \\
  320$\times$  320	&	7.07E-05&   4.00    &	7.07E-05&   4.00 	&	7.30E-05&   4.47 \\
\hline													
\end{tabular}
\end{table}	

\begin{figure}[h]
\centering
\subfloat{
\includegraphics[width=0.5\textwidth]{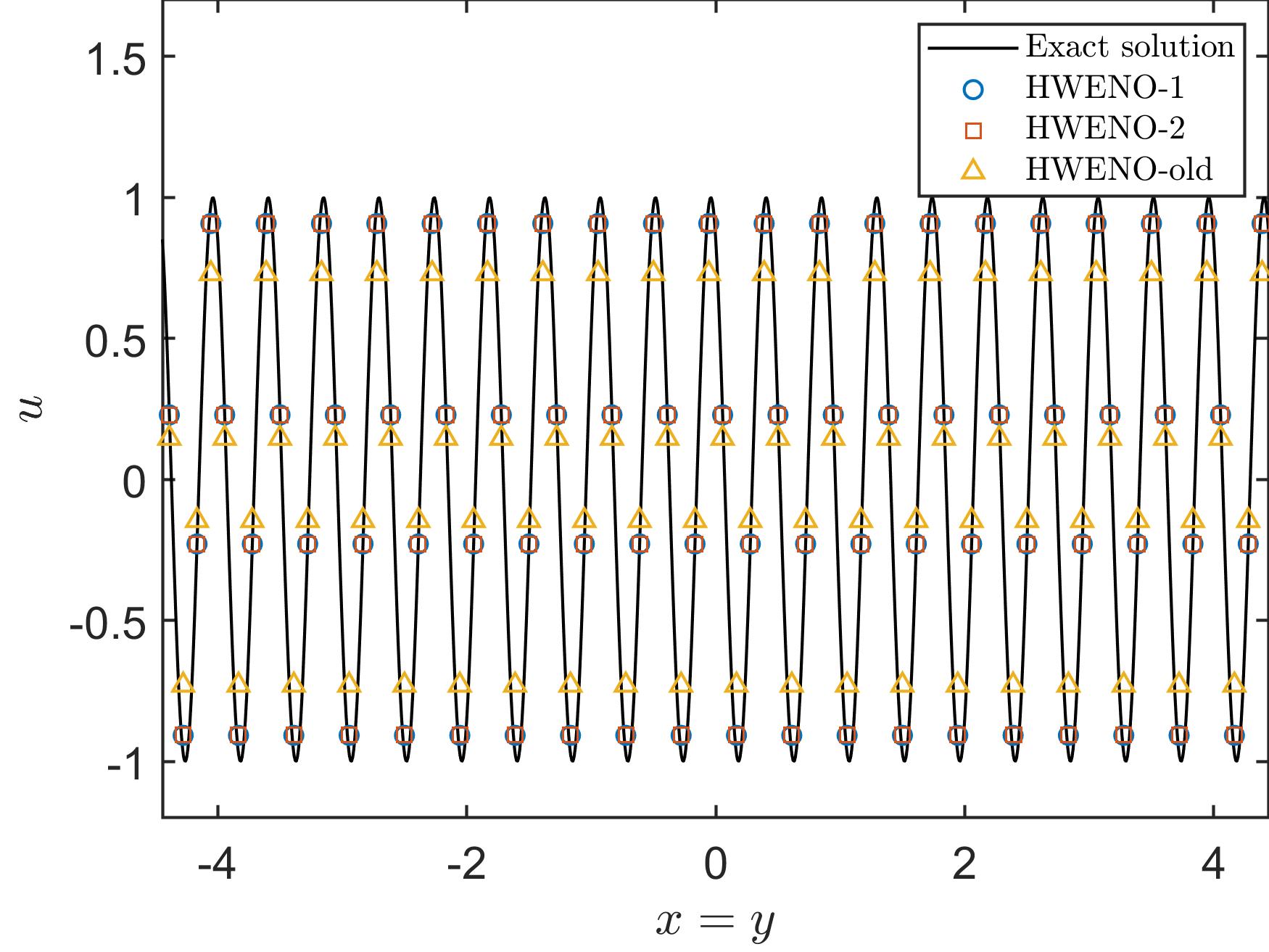}
}
\caption{(Transport equation with constant coefficients). Cross-sections of SL HWENO solutions at $T = 20$ at $x=y$ for (\ref{2-D_linear_transport}) with $u(x,y,0) = \text{sin}(x+y)$ at $T = 20$ with a mesh of $80\times80$.}\label{fig:slide_const}
\end{figure}								

\end{exa}

\begin{exa}(Swirling deformation flow). Consider
\begin{equation}\label{2_D_SDF}
\begin{split}
  u_t-(2\pi\text{cos}^2(\frac x2)\text{sin}(y)g(t)u)_x + (2\pi\text{sin}(x)\text{cos}^2(\frac y2)g(t)u)_y=0,\\
 x\in[-\pi,\pi],~~ y\in[-\pi,\pi],
 \end{split}
\end{equation}
where $g(t) = \text{cos}(\pi t/T)$ with $T = 1.5$. We consider \eqref{2_D_SDF} with the following smooth initial condition
\begin{equation}\label{eq:SDF_smooth_initial_condtion}
  u(x,y,0) = \begin{cases}
  r^b_0\text{cos}(\frac{r^b(\mathbf{x})\pi}{2r^b_0})^6~~&\text{if}~ r^b(\mathbf{x})<r^b_0,\\
  0,~~&\text{otherwise},
  \end{cases}
\end{equation}
where $r^b_0=0.3\pi$, $r^b(\mathbf{x})=\sqrt{ (x-x_0^b)^2+(y-y_0^b)^2 }$ and the center of the cosine bell $(x_0^b,y_0^b) = (0.3\pi,0)$. Zero boundary condition is equipped for this test. In \Cref{tab_2_D_SDF}, we present the $L^2$ errors and corresponding orders of accuracy of the SL HWENO-1 , HWENO-2, and HWENO-old schemes. As shown, the SL HWENO schemes are of fourth-order accuracy. The HWENO-1 shows the best convergence rate comparing with the other two HWENO reconstructions. The convergence rate of the HWENO-old is worst. This leads to the same conclusion that HWENO-old reconstruction smears more information. In \Cref{fig:Temporal_order_SDF}, we show the temporal order of accuracy of the SL HWENO-1 scheme by fixing the spatial mesh while varying the time step $\Delta t$. We find that the temporal order reaches fifth-order, which is one order higher than the expected fourth-order, for the swirling deformation flow.

\begin{table}[!htbp]
\centering
\caption{ (Swirling deformation flow).  $L^2$ errors and corresponding orders of accuracy of the SL HWENO schemes for \eqref{2_D_SDF} with initial condition \eqref{eq:SDF_smooth_initial_condtion} at $t = 1.5$.}\label{tab_2_D_SDF}
  \centering
\begin{tabular}{|c|cc|cc|cc|}
\hline
&\multicolumn{2}{|l|}{\textbf{HWENO-1}}&\multicolumn{2}{|l|}{\textbf{HWENO-2}}&\multicolumn{2}{|l|}{\textbf{HWENO-old}}\\
\cline{2-7}
mesh&$L^2$ error&order&$L^2$ error&order&$L^2$ error&order\\
\hline
  20$\times$  20    &   1.57E-02&   ---         &   1.59E-02&   ---         &   2.76E-02&   ---     \\
  40$\times$  40	&   9.57E-04&   4.04	 	&   1.70E-03&   3.22	 	&	7.35E-03&	1.91	\\
  80$\times$  80	&	1.32E-04&   2.85 	 	&	2.79E-04&   2.61 	 	&	7.84E-04&   3.23 	\\
  160$\times$  160	&	5.25E-06&   4.66 	 	&	2.75E-05&   3.34 	 	&	5.54E-05&   3.82 	\\
  320$\times$  320	&	2.23E-07&   4.56 	 	&	1.52E-06&   4.18 	 	&	8.14E-06&   2.77 	\\
  640$\times$  640	&	9.03E-09&   4.63 	 	&	6.45E-08&   4.56 	 	&	6.32E-07&   3.69 	\\
  1280$\times$ 1280 &	4.12E-10&   4.45 	 	&	4.12E-10&   7.29 	 	&	8.07E-09&   6.29 	\\
  2560$\times$ 2560 &   2.25E-11&   4.20        &   2.25E-11&   4.20        &   2.52E-11&   8.32    \\
\hline														
\end{tabular}
\end{table}

\begin{figure}[!htbp]
\centering
\subfloat{
\includegraphics[width=0.35\textwidth]{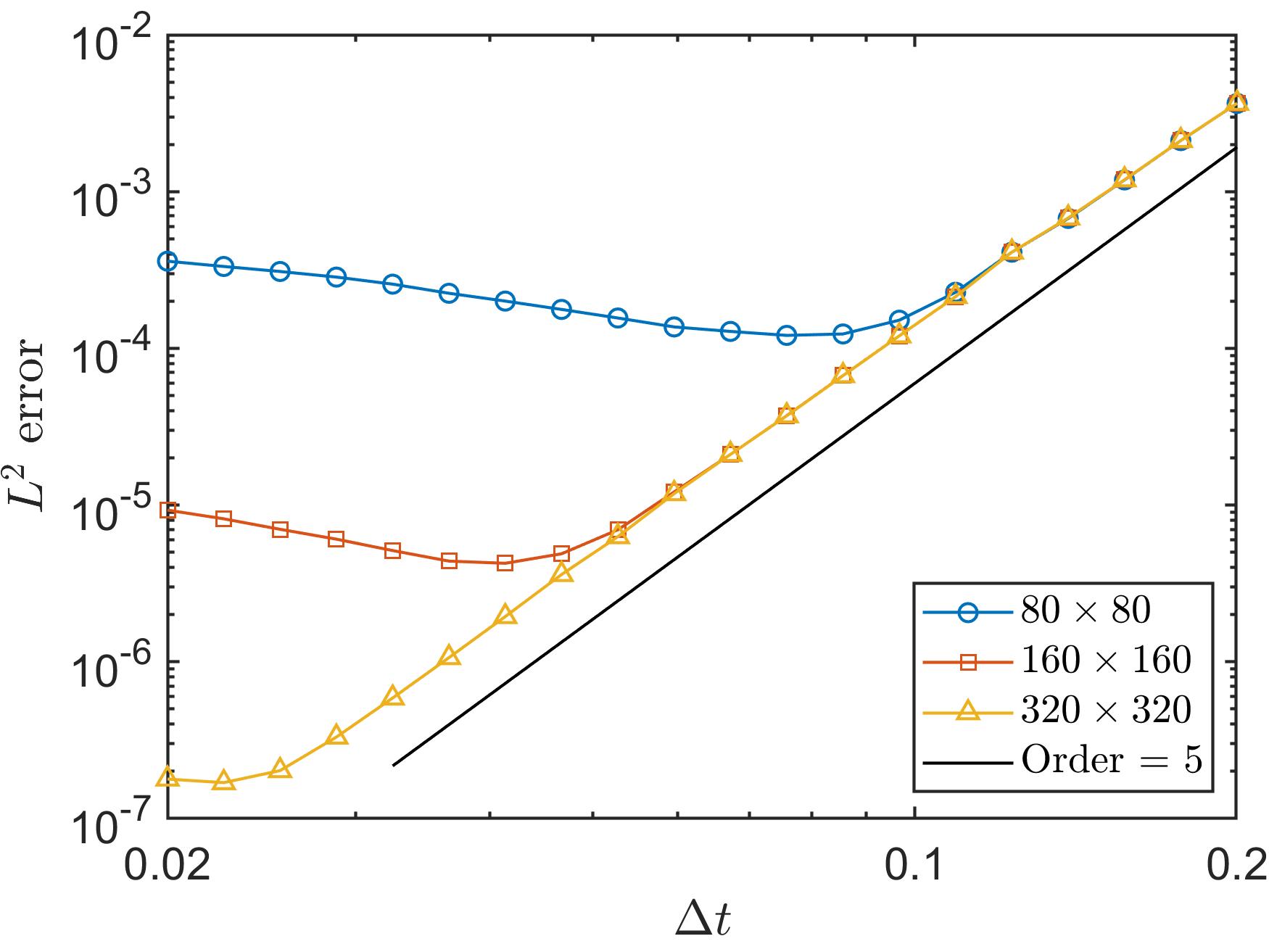}
}
\caption{(Swirling deformation flow). Temporal order of accuracy of the SL HWENO-1 scheme. The three colored lines use different mesh of $80\times80$, $160\times160$, and $320\times320$. The programs stop at $t=1.5$.}\label{fig:Temporal_order_SDF}
\end{figure}	

To investigate the performances of the SL HWENO schemes for discontinuous solutions, we test \eqref{2_D_SDF} with the initial condition shown in \Cref{fig:SDF_ini}. We show the mesh plots of the numerical solutions of the SL HWENO schemes at $t = 0.75$ and at $t = 1.5$ in \Cref{fig:SDF_0_75_1_5} with a mesh of $100\times100$. In \Cref{fig:Cross_sections}, we present two cross-sections of the numerical solutions of the SL HWENO schemes at $t=1.5$ with a mesh of $100\times100$. As shown, the SL HWENO schemes preserve the geometry of the solution. The SL HWENO-1 scheme and HWENO-2 schemes have better resolution but produces small numerical oscillation at the same time. In \Cref{fig:Cross_sections}, we observe that the SL HWENO-old scheme is more dissipative than the other two.

\begin{figure}[h]
\centering
\subfloat{
\includegraphics[width=0.35\textwidth]{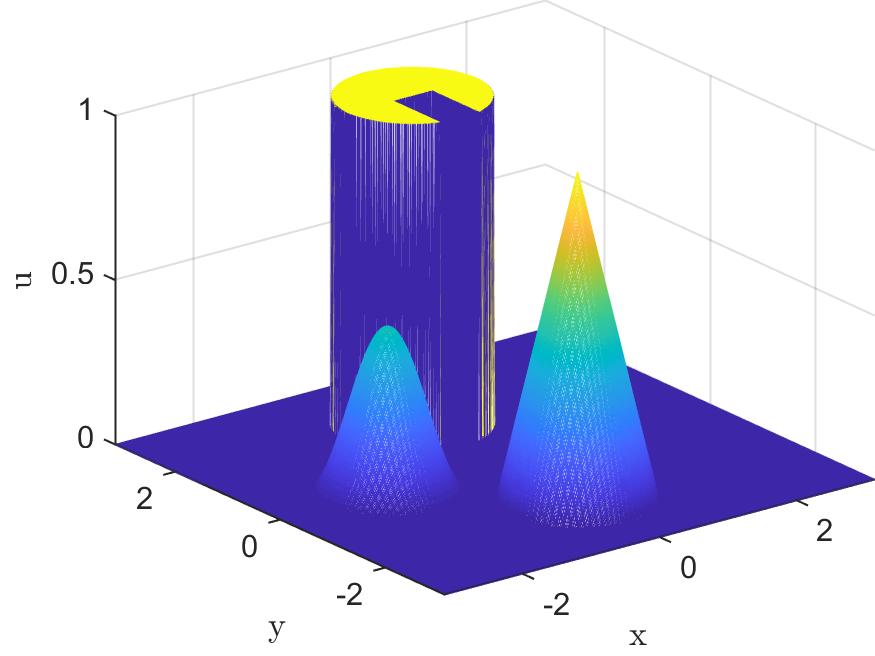}
}
\subfloat{
\includegraphics[width=0.35\textwidth]{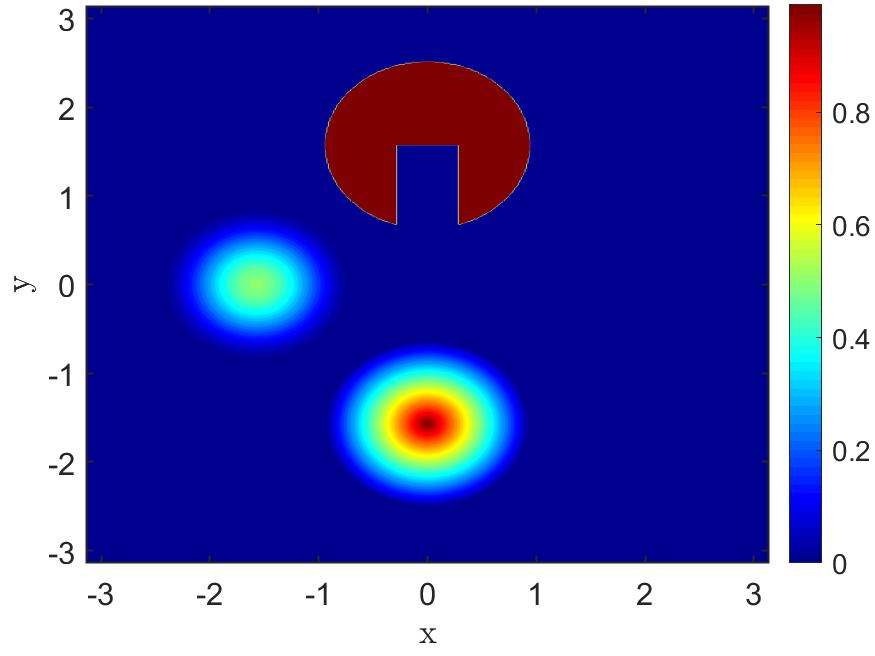}
}
\caption{(Swirling deformation flow). The mesh plot (left) and the contour plot (right) of the discontinuous initial data for \eqref{2_D_SDF}.}\label{fig:SDF_ini}
\end{figure}				

\begin{figure}[!htbp]
\centering
\subfloat{
\includegraphics[width=0.3\textwidth]{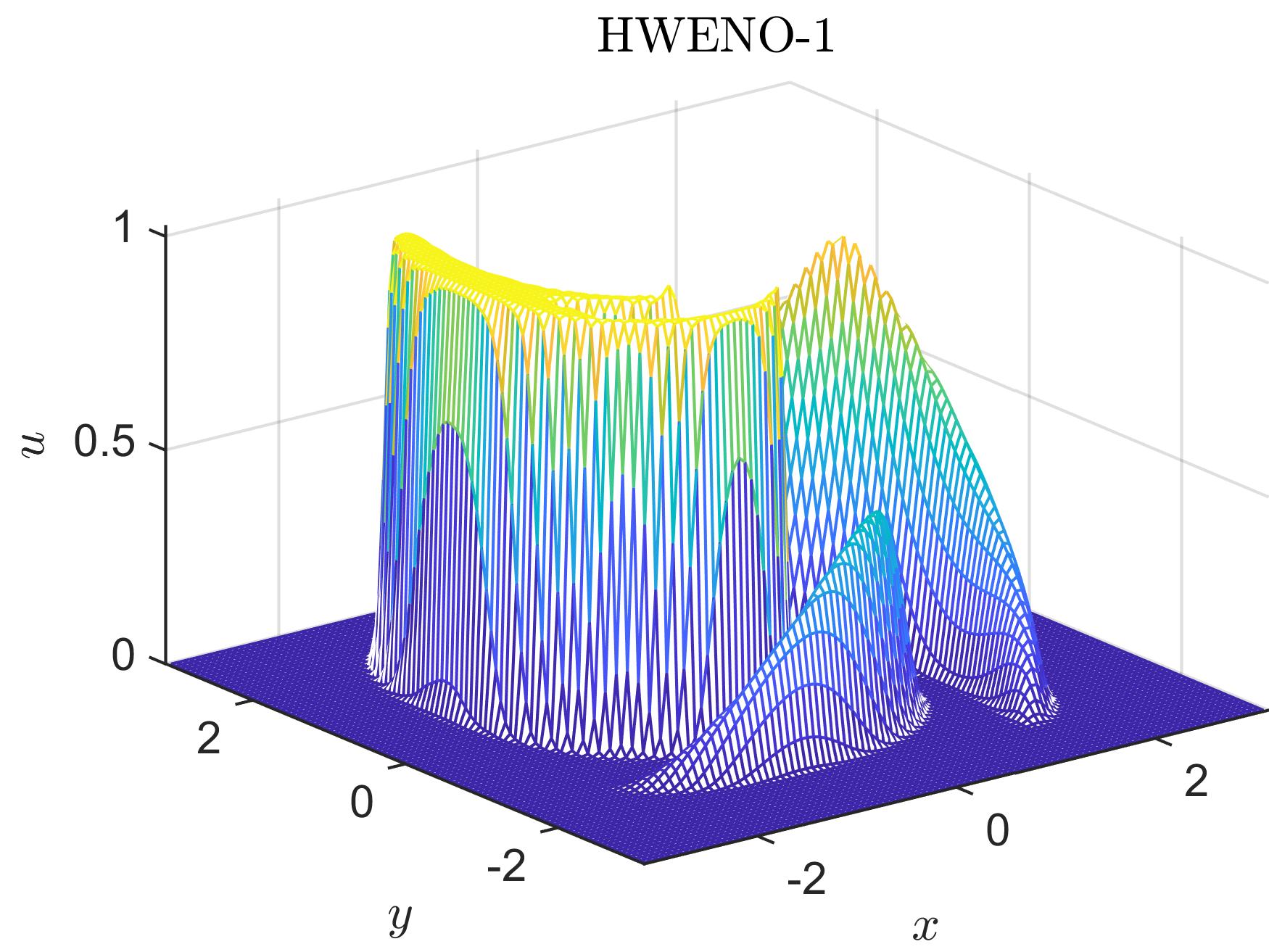}
}
\subfloat{
\includegraphics[width=0.3\textwidth]{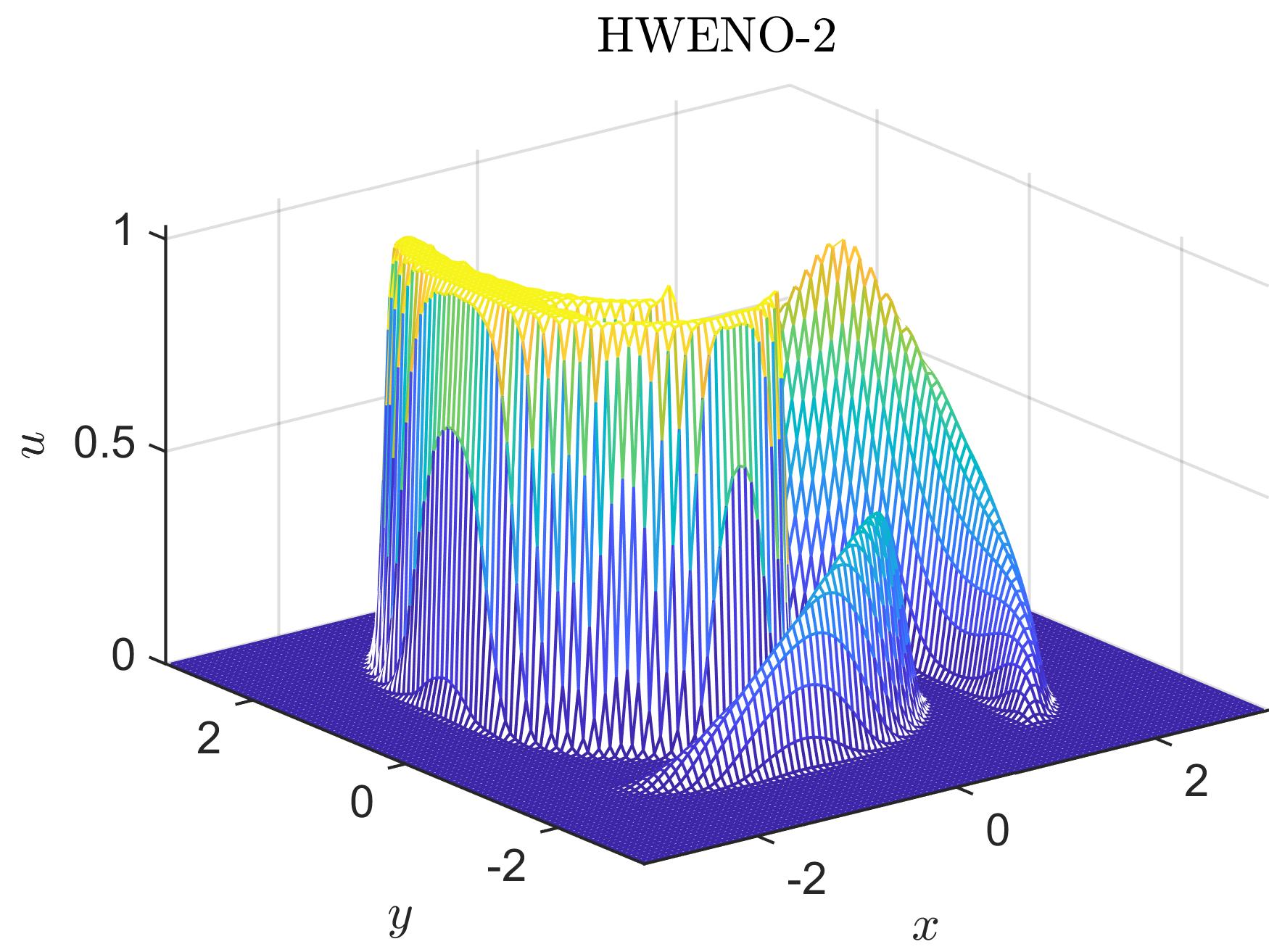}
}
\subfloat{
	\includegraphics[width=0.3\textwidth]{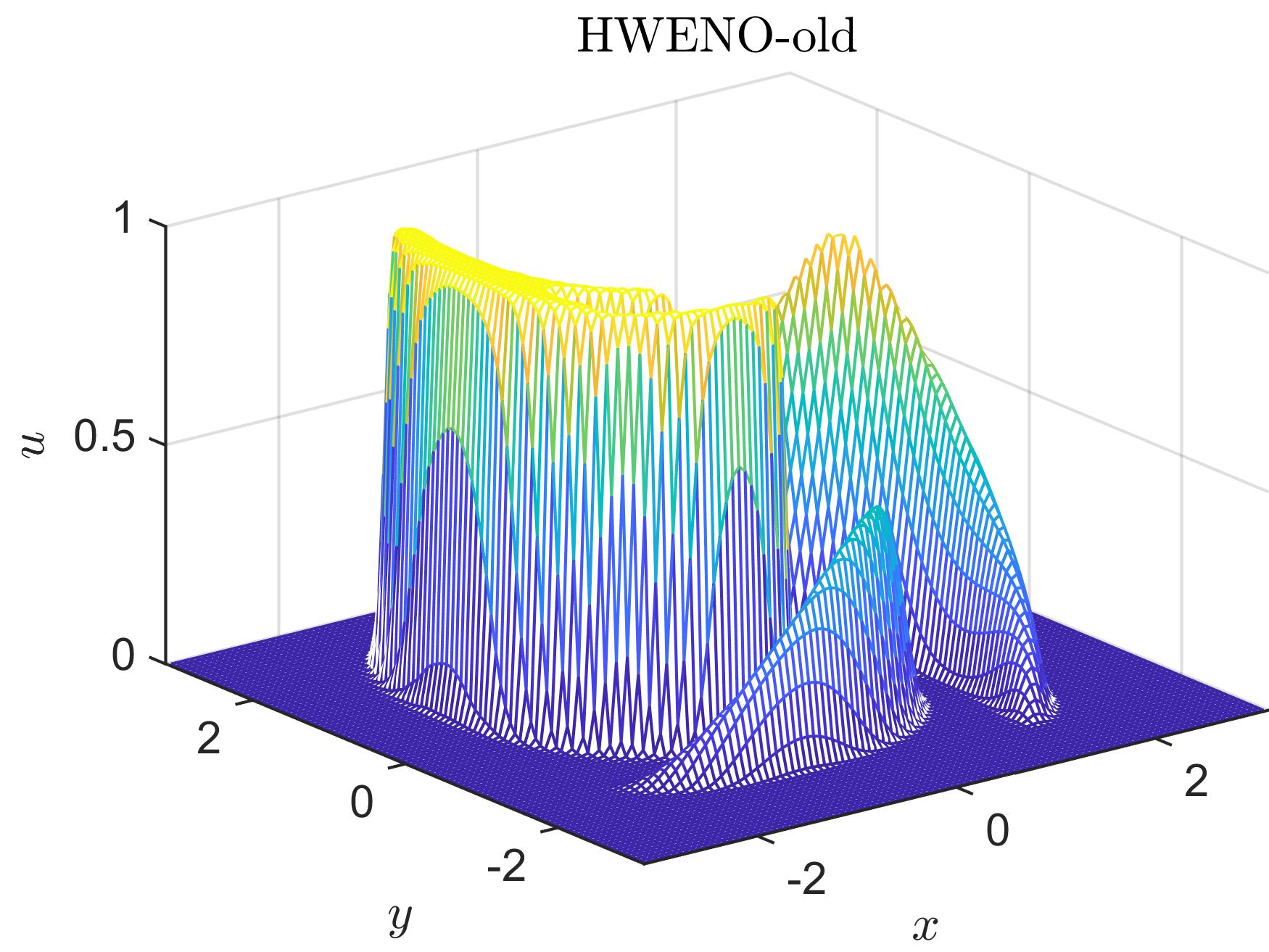}
}

\subfloat{
	\includegraphics[width=0.3\textwidth]{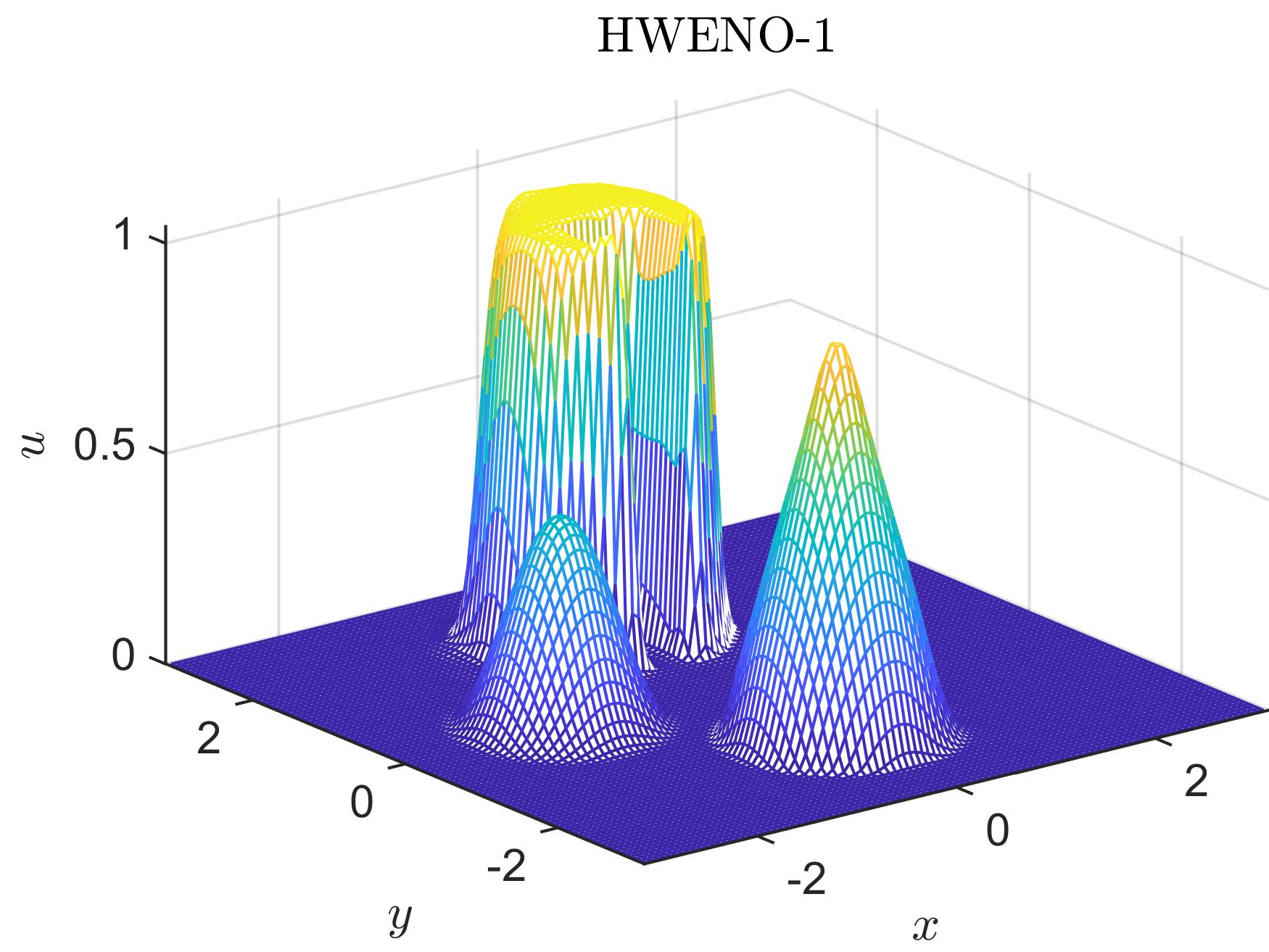}
}
\subfloat{
	\includegraphics[width=0.3\textwidth]{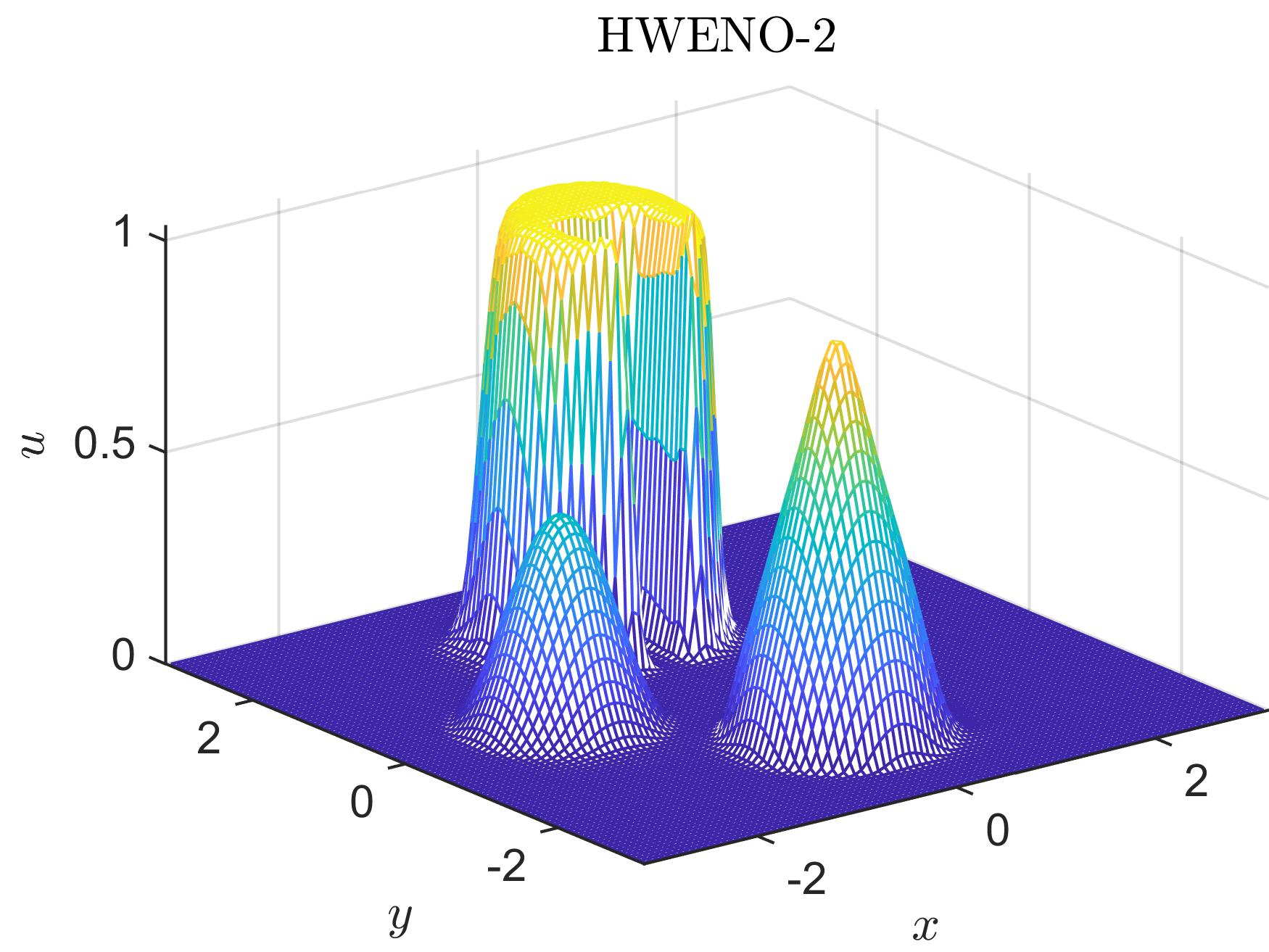}
}
\subfloat{
	\includegraphics[width=0.3\textwidth]{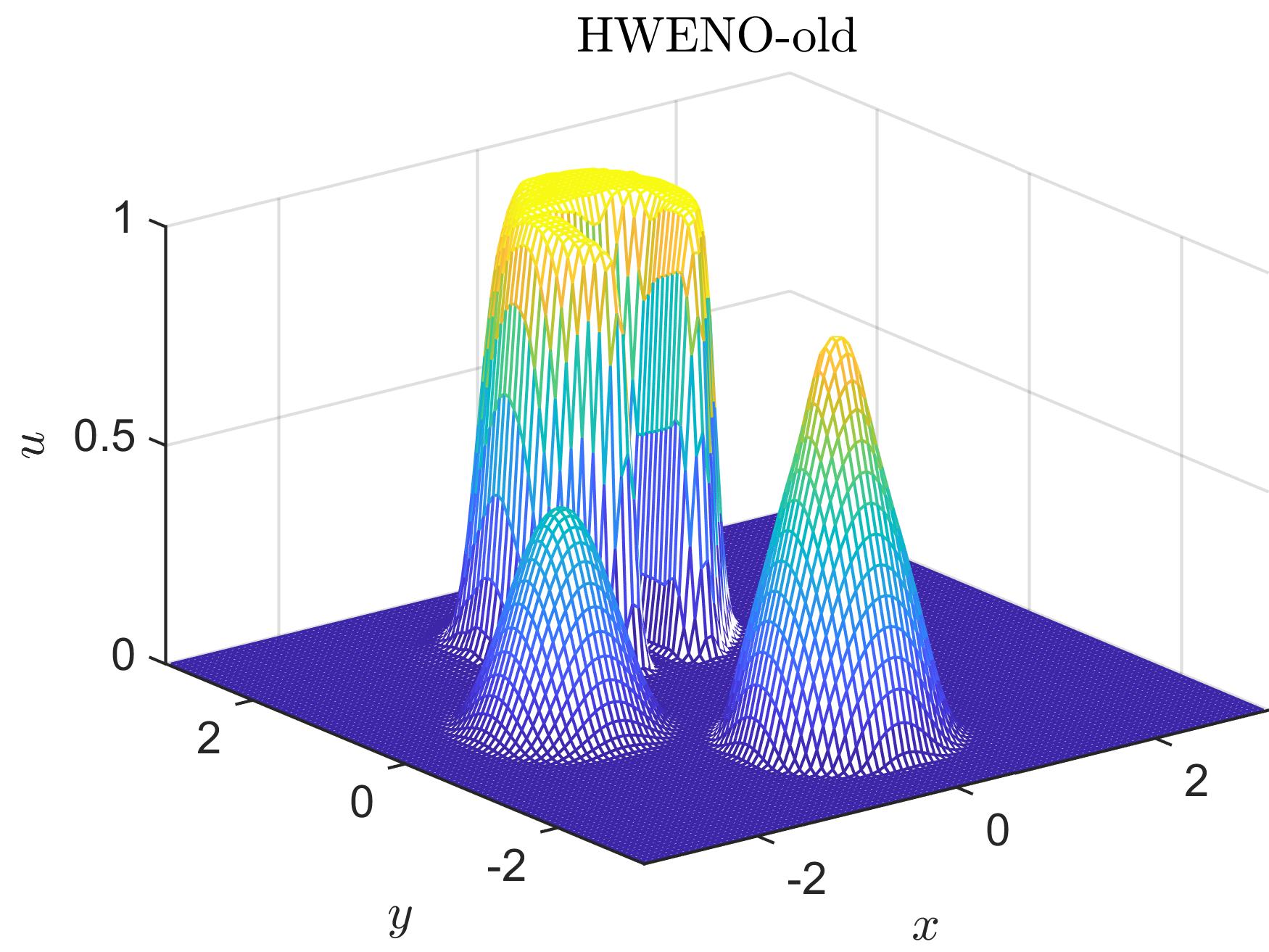}
}
\caption{(Swirling deformation flow). The mesh plots of the numerical solutions of the SL HWENO-1 (left), HWENO-2 (middle), and HWENO-old (right) schemes for \eqref{2_D_SDF} with initial condition \Cref{fig:SDF_ini} at $t = 0.75$ (top) and at $t = 1.5$ (bottom) with a mesh of $100\times100$.}\label{fig:SDF_0_75_1_5}
\end{figure}	
\begin{figure}[!htbp]
\centering
\subfloat{
\includegraphics[width=0.35\textwidth]{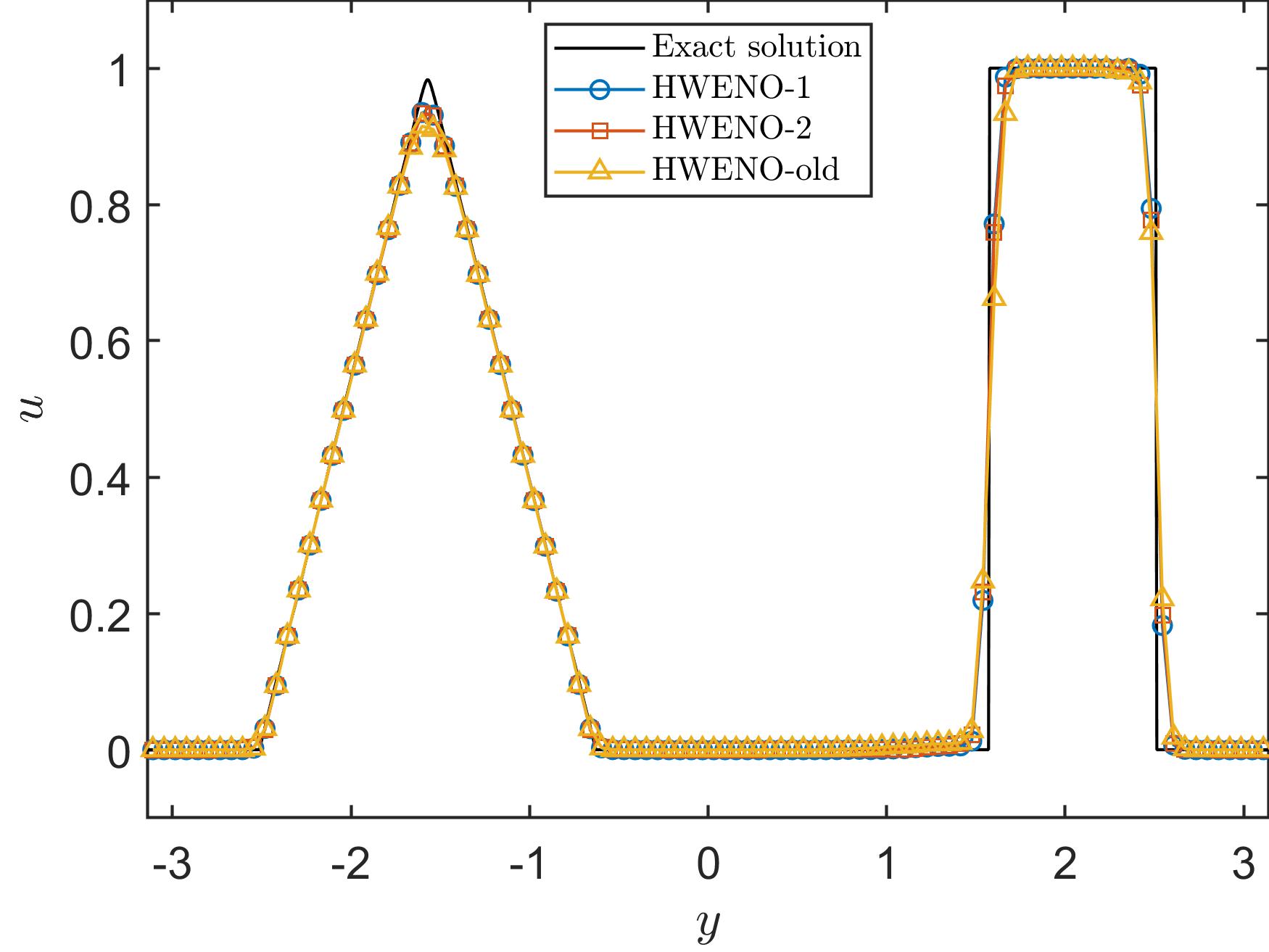}
}
\subfloat{
\includegraphics[width=0.35\textwidth]{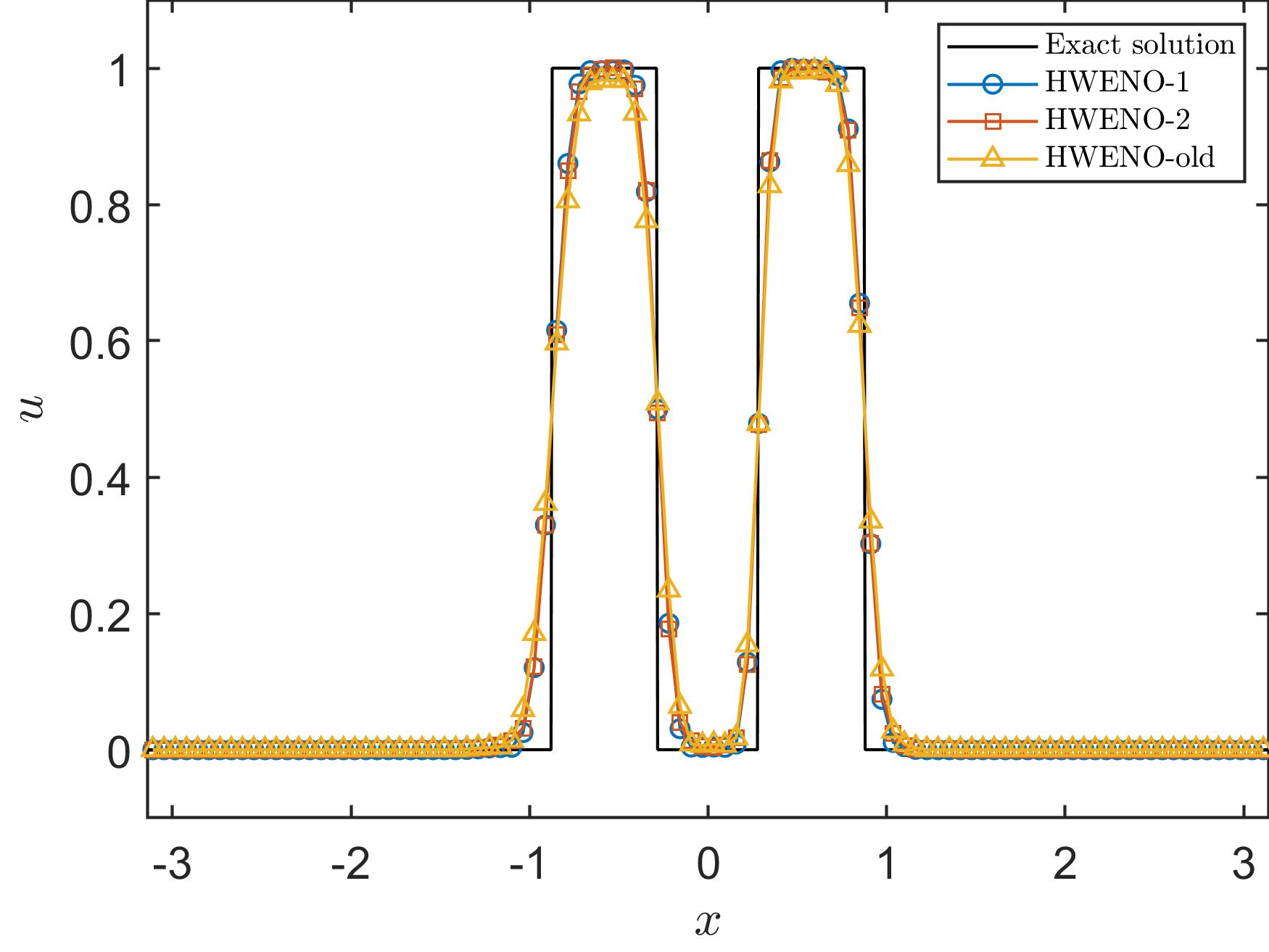}
}
\caption{(Swirling deformation flow). Cross-sections of the SL HWENO solutions at $x=0$ (left) and $y=1.2$ (right) at $t=1.5$ with a mesh of $100\times100$.}\label{fig:Cross_sections}
\end{figure}									
\end{exa}

\subsection{Nonlinear Vlasov-Poisson system}

In this subsection, we test two benchmark tests: the strong Landau damping and the bump-on-tail instability.  Unless specified, we set  $N_x=128$, $N_v = 256$, CFL = $10.2$,  $\Delta t = \text{CFL}/\left(v_{\text{max}}/\Delta x+\text{max}\{|E|\}/\Delta v\right)$ with $v_{\text{max}}$ represents the positive boundary of $v$-direction. The PP limiter is equipped for both tests.

\begin{exa} (Strong Landau damping).
 Consider the VP system with the initial condition
\begin{equation}\label{landau_damping}
f(x,v,t=0) = \frac1{\sqrt{2\pi}}\left(1+\alpha\text{cos}(kx)\right)\text{exp}\left(-\frac{v^2}{2}\right),\quad x\in[0,4\pi],\quad v\in[-2\pi,2\pi],
\end{equation}
where $k=0.5$, $\alpha = 0.5$.
 In \Cref{tab_2-D_SLD}, we show the $L^2$ errors and corresponding orders of accuracy of the SL HWENO schemes for the strong Landau damping at $T = 2$. The errors are obtained by comparing the solution to a reference solution with mesh refinement. We observe the expected order of accuracy and the three reconstructions perform similarly for this test. In \Cref{fig:Temporal_order_SLD}, we present the temporal order of accuracy of the SL HWENO-1 scheme by fixing the spatial mesh while varying $\Delta t$. We observe a clear fourth-order temporal accuracy.

\begin{table}[!htbp]
\centering
\caption{ (Strong Landau damping). $L^2$ errors and corresponding orders of accuracy of the SL HWENO schemes at $T = 2$.}\label{tab_2-D_SLD}
  \centering
\begin{tabular}{|c|cc|cc|cc|}
\hline
&\multicolumn{2}{|l|}{\textbf{HWENO-1}}&\multicolumn{2}{|l|}{\textbf{HWENO-2}}&\multicolumn{2}{|l|}{\textbf{HWENO-old}}\\
\cline{2-7}
mesh&$L^2$ error&order&$L^2$ error&order\\
\hline
  16$\times$  16	&	4.29E-03	&	---	 	&	4.22E-03	&	---	 	&	5.07E-03&	---	 	\\
  32$\times$  32	&	4.52E-04&   3.25 	 	&	4.47E-04&   3.24 	 	&	4.70E-04&   3.43 	 	\\
  64$\times$  64	&	2.10E-05&   4.43 	 	&	2.10E-05&   4.42 	 	&	1.91E-05&   4.62 	 	\\
  128$\times$  128	&	6.15E-07&   5.10 	 	&	6.15E-07&   5.09 	 	&	6.12E-07&   4.97 	 	\\
  256$\times$  256	&	2.69E-08&   4.51 	 	&	2.70E-08&   4.51 	 	&	2.69E-08&   4.51 	 	\\
\hline									
\end{tabular}
\end{table}

\begin{figure}[!htbp]
\centering
\subfloat{
\includegraphics[width=0.35\textwidth]{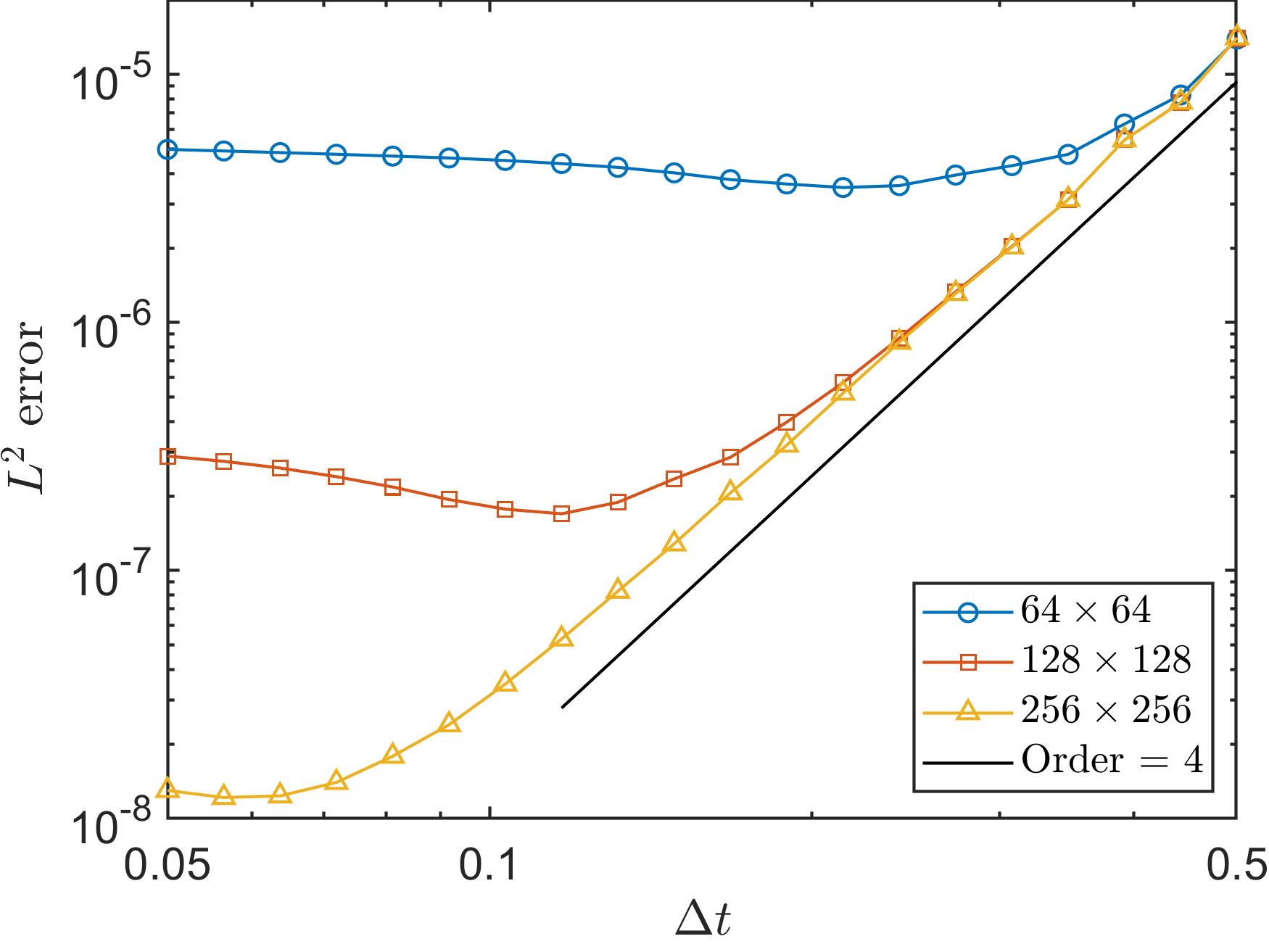}
}
\caption{(Strong Landau damping). Temporal order of accuracy of the SL HWENO-1 scheme. The three colored lines use different mesh of $64\times64$, $128\times128$, and $256\times256$. The programs stop at $T=2$.}\label{fig:Temporal_order_SLD}
\end{figure}	

In \Cref{fig:SLD_40}, we present the mesh plots and contour plots of numerical solutions of the SL HWENO schemes at $T = 40$. We observe that the numerical solution maintains the filamentation structure well for the strong Landau damping. In \Cref{fig:SLD_slide}, the cross-sections of the SL HWENO schemes at $T = 40$ and at $x=2\pi$ are provided. We observe slightly better resolution for the SL HWENO-1 and HWENO-2 schemes.

\begin{figure}[!htbp]
\centering
\subfloat{
\includegraphics[width=0.3\textwidth]{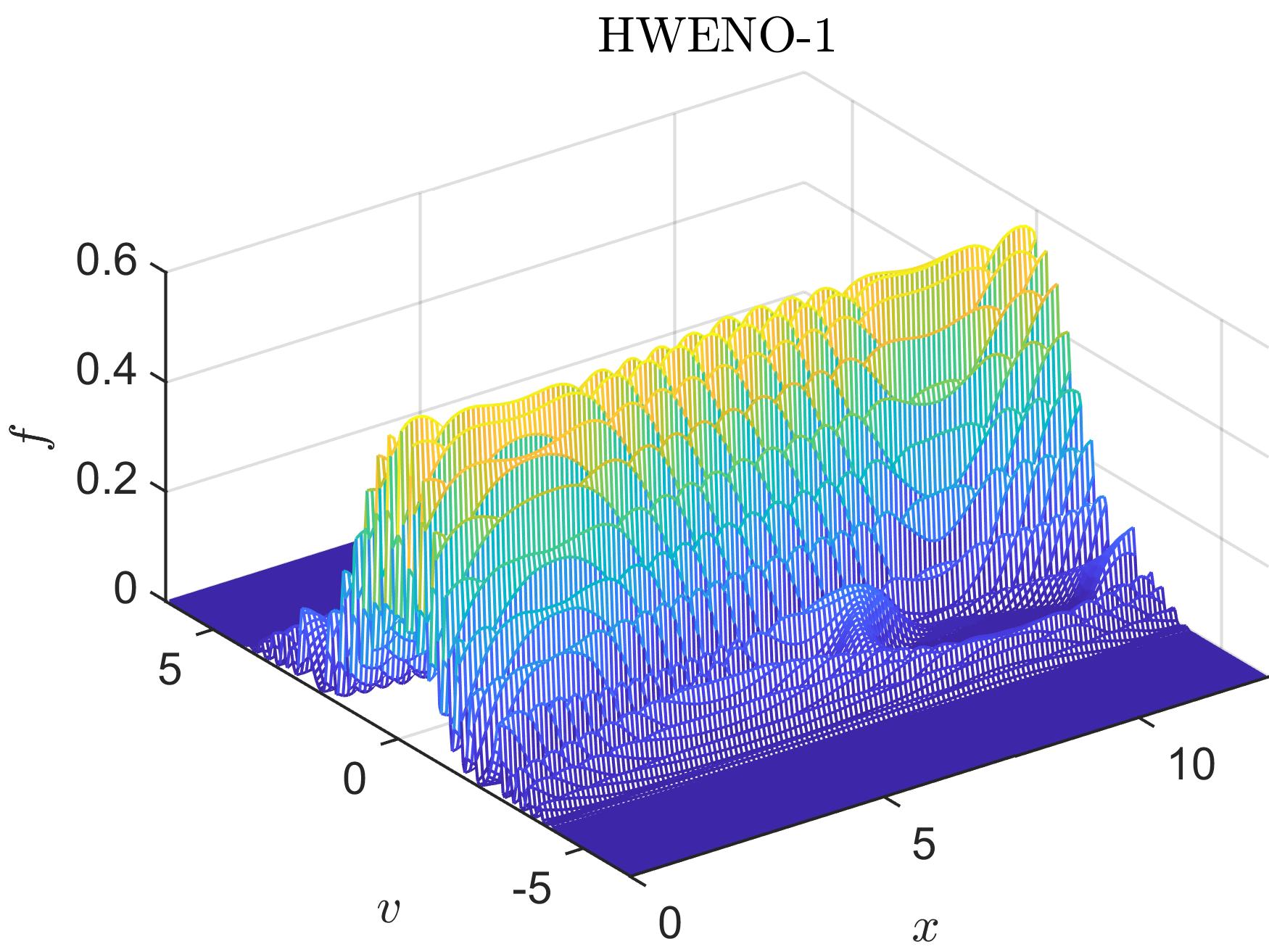}
}
\subfloat{
	\includegraphics[width=0.3\textwidth]{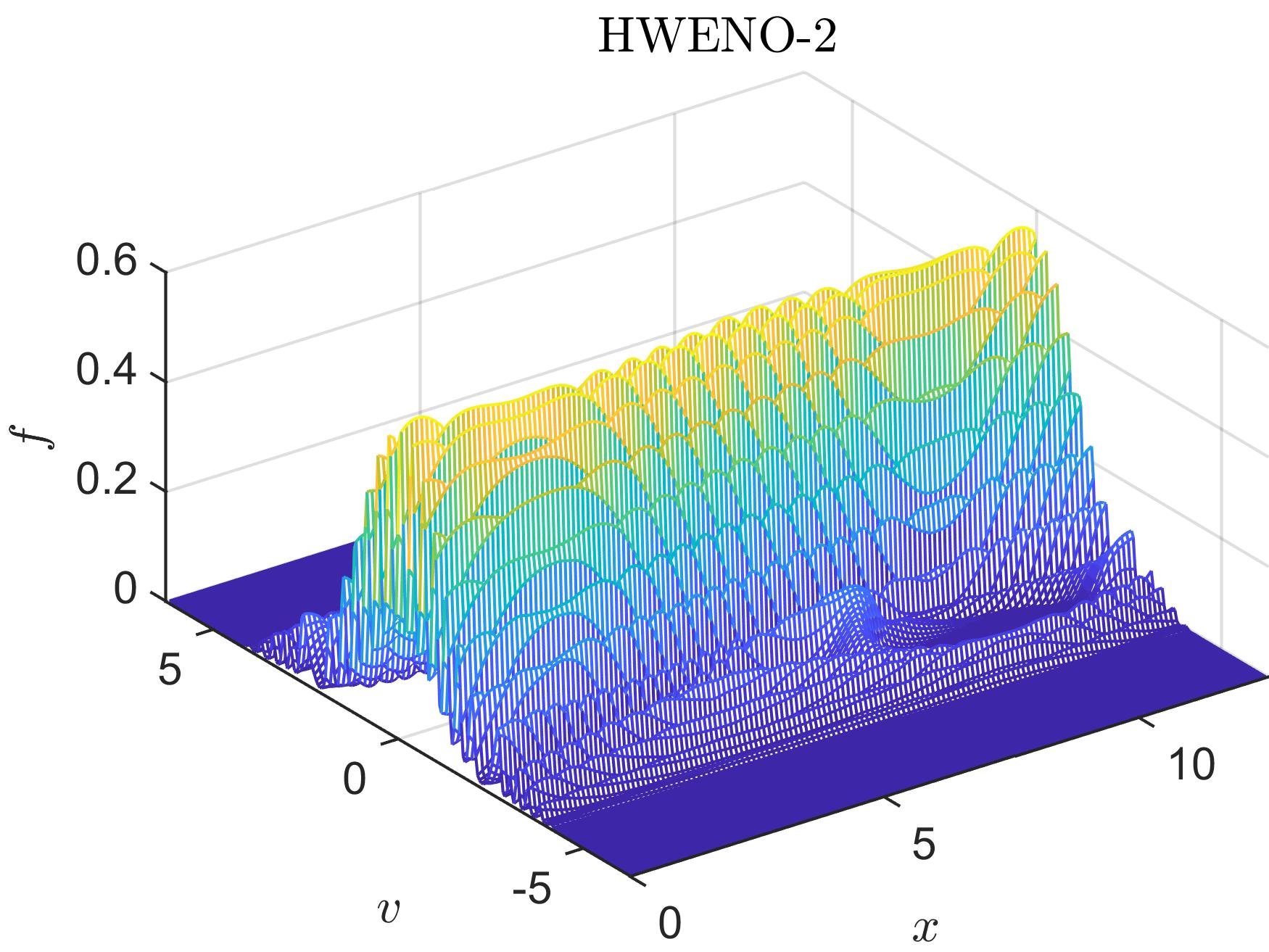}
}
\subfloat{
	\includegraphics[width=0.3\textwidth]{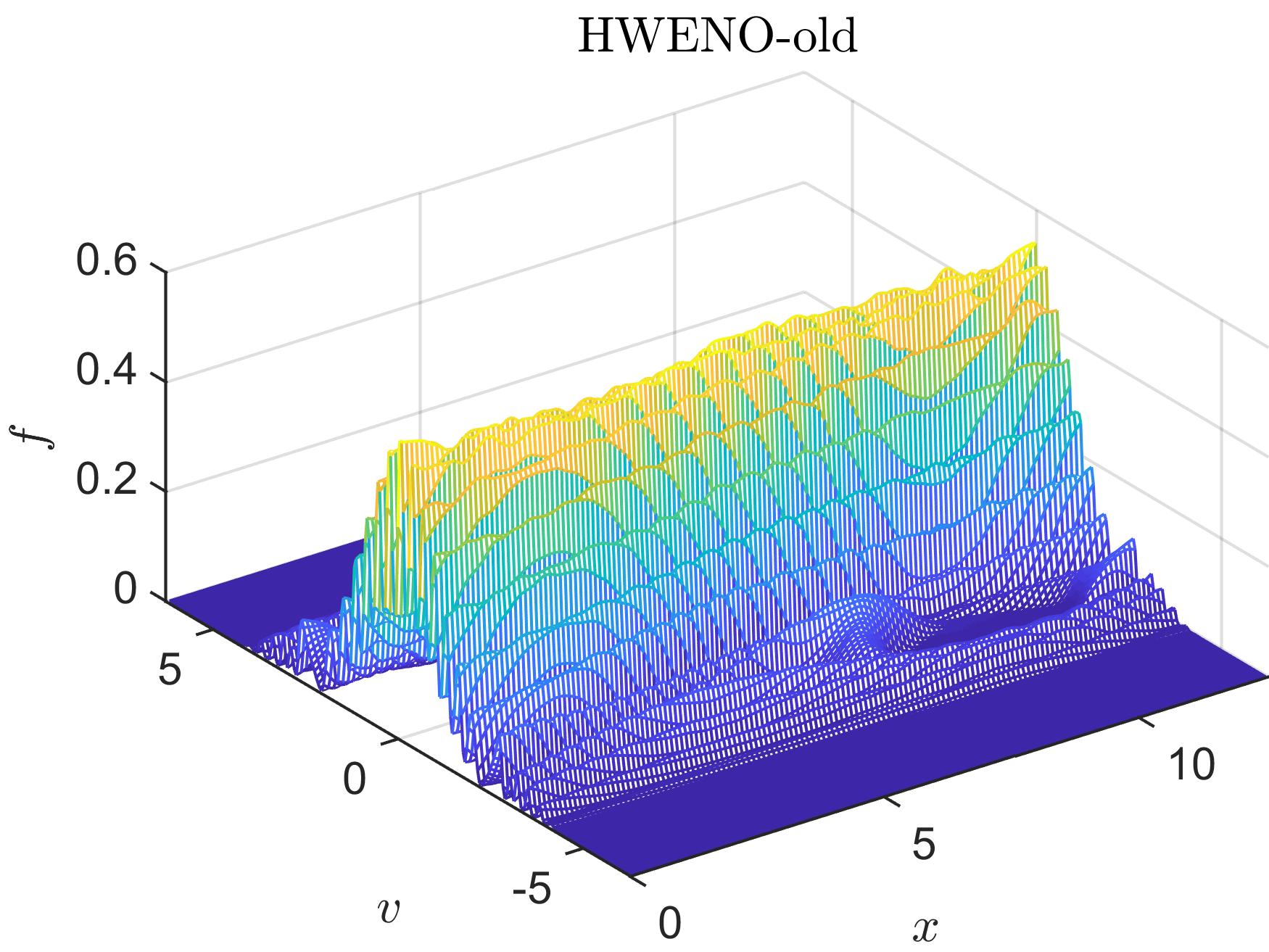}
}

\subfloat{
	\includegraphics[width=0.3\textwidth]{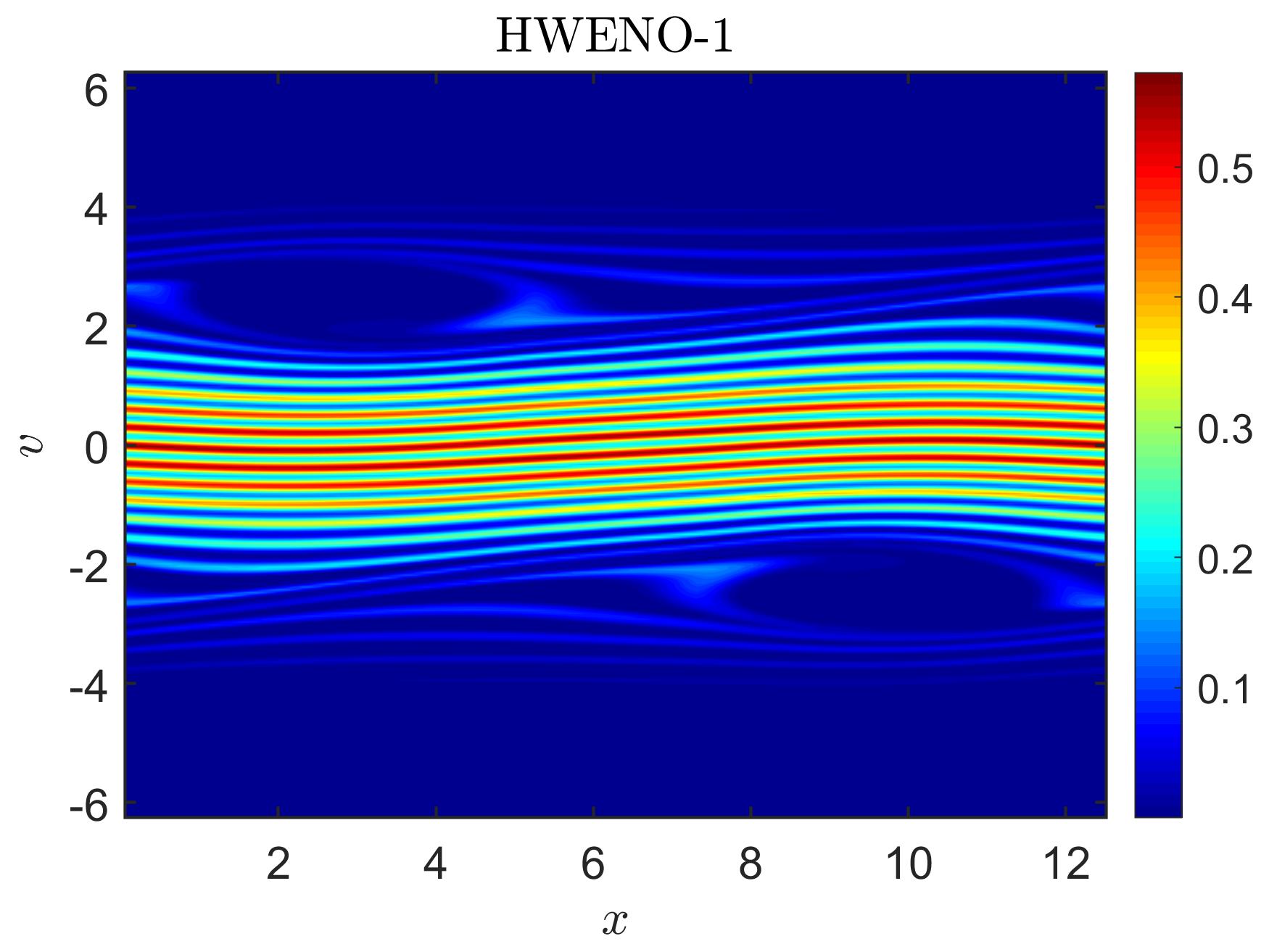}
}
\subfloat{
	\includegraphics[width=0.3\textwidth]{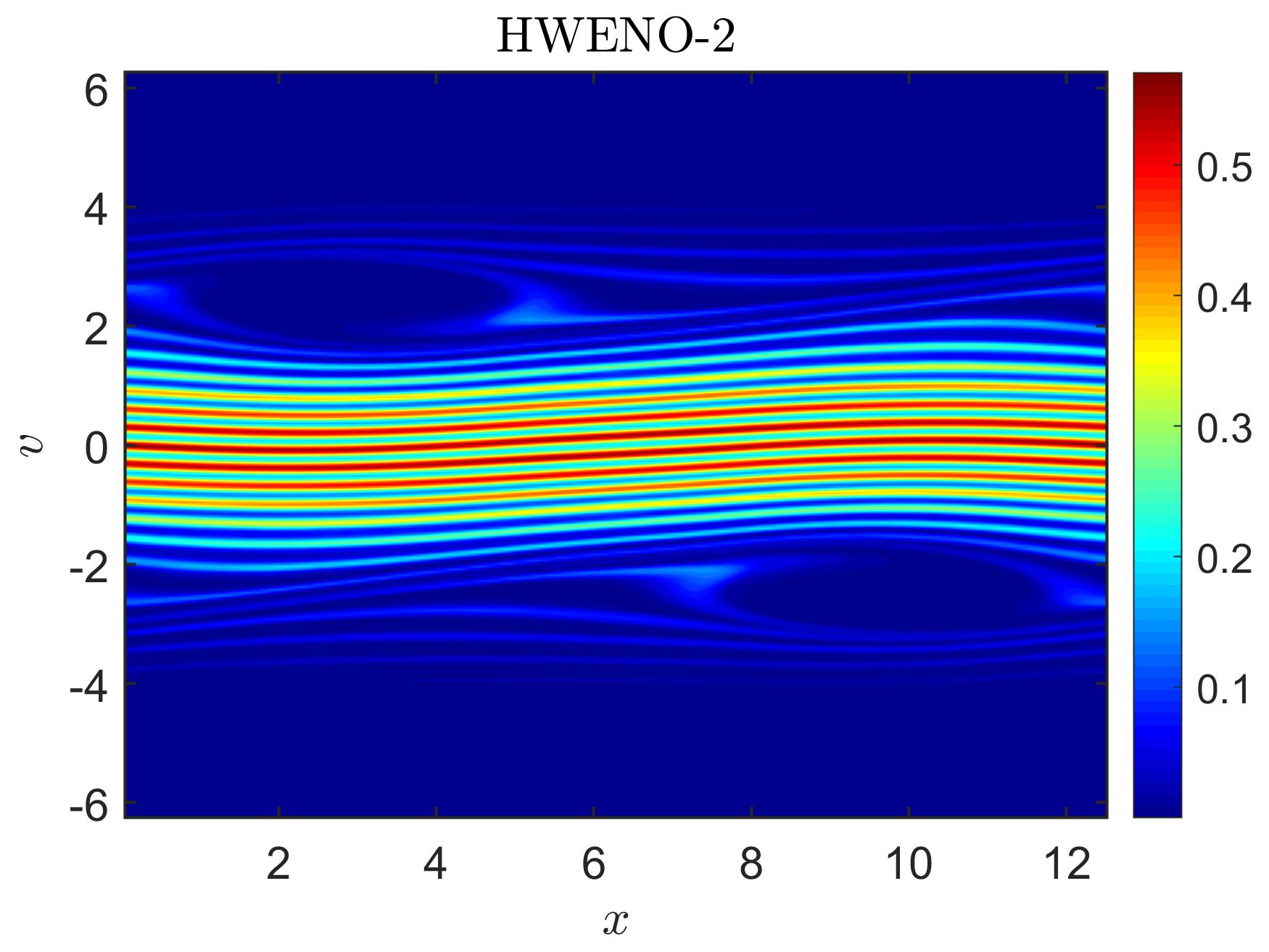}
}
\subfloat{
	\includegraphics[width=0.3\textwidth]{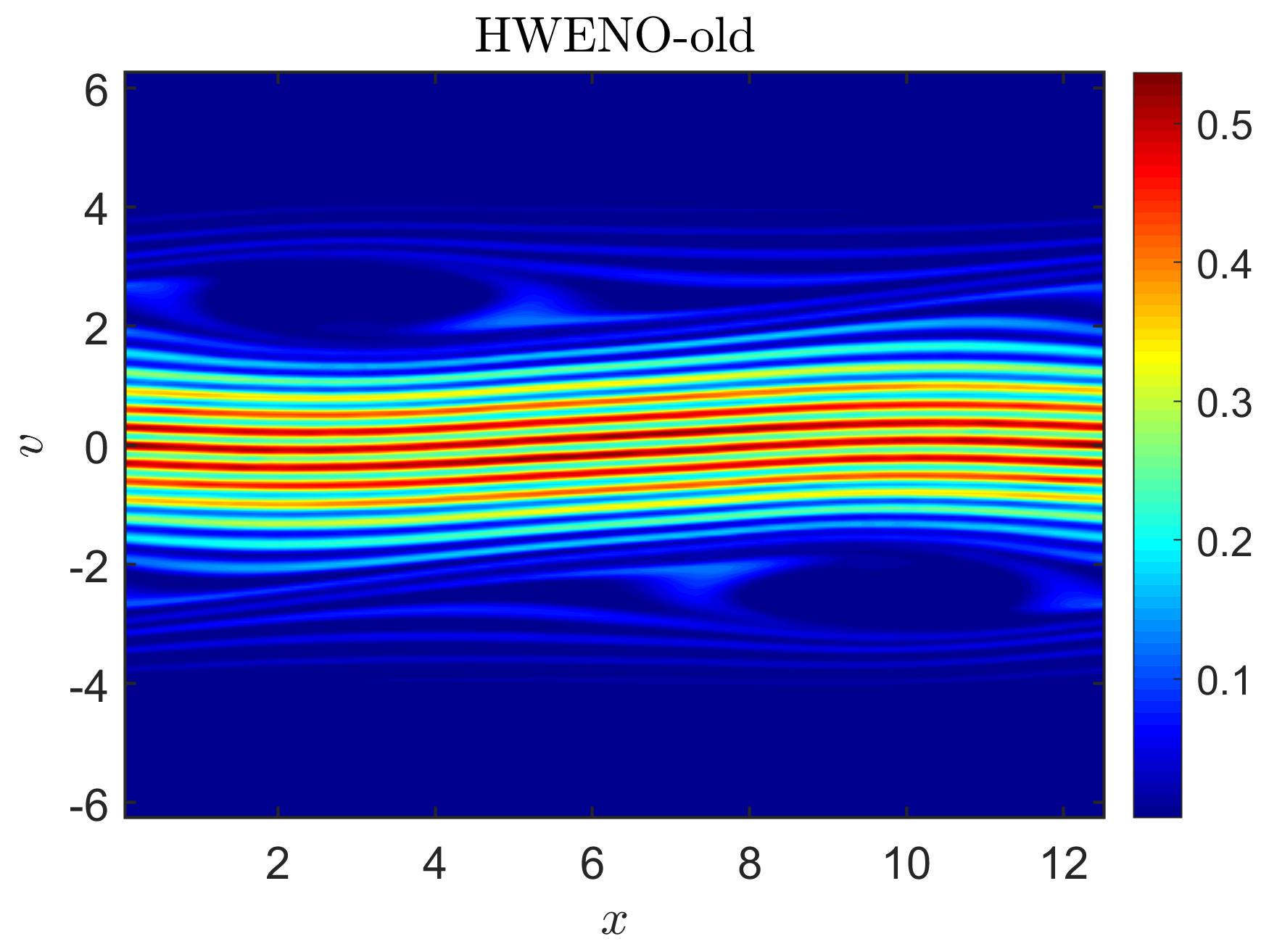}
}
\caption{(Strong Landau damping). The mesh plots (top) and contour plots (bottom) of the numerical solutions of the SL HWENO-1 (left), HWENO-2 (middle), and HWENO-old (right) schemes at $T = 40$.}\label{fig:SLD_40}
\end{figure}

\begin{figure}[!htbp]
	\centering
	\subfloat{
		\includegraphics[width=0.4\textwidth]{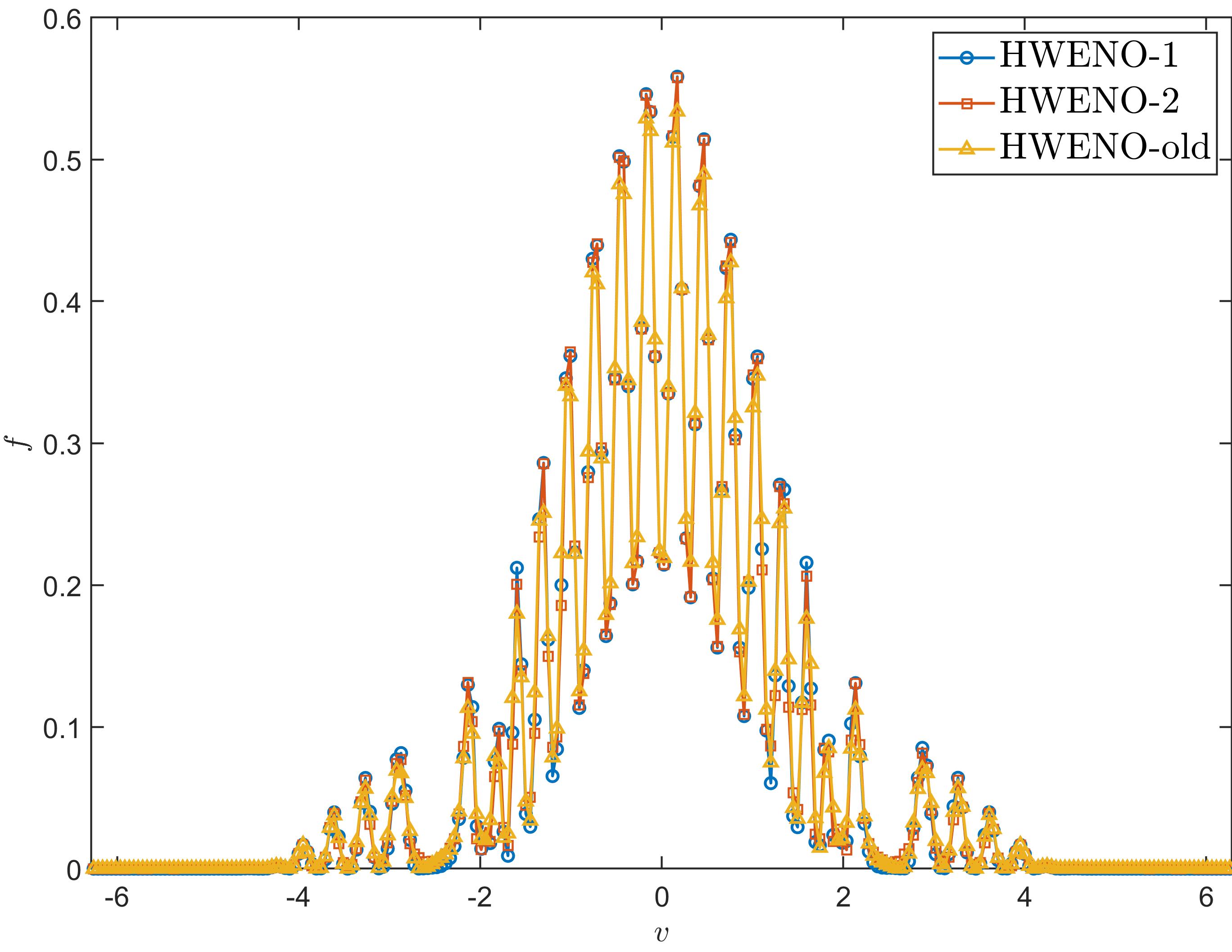}
	}
	\caption{(Strong Landau damping). Cross-sections of the SL HWENO solutions at $T = 40$ and at $x = 2\pi$.}\label{fig:SLD_slide}
\end{figure}

The VP system has several conservative quantities such as mass, $L^p$ norms, energy, and entropy \cite{ROSSMANITH20116203}. We can observe a $O(10^{-12})$ level of relative deviation of mass and $L^1$ norm with $v_{\text{max}}=10$ for strong Landau damping as shown in \Cref{fig:SLD_mass_L1}. For $L^2$ norm, energy, and entropy, we present the relative deviations of them for the SL HWENO schemes in \Cref{fig:SLD_conserve}. As shown, the magnitude of deviations are similar to the existing results \cite{ROSSMANITH20116203,QIU20118386}.

\begin{figure}[!htbp]
	\centering
	\subfloat{
		\includegraphics[width=0.35\textwidth]{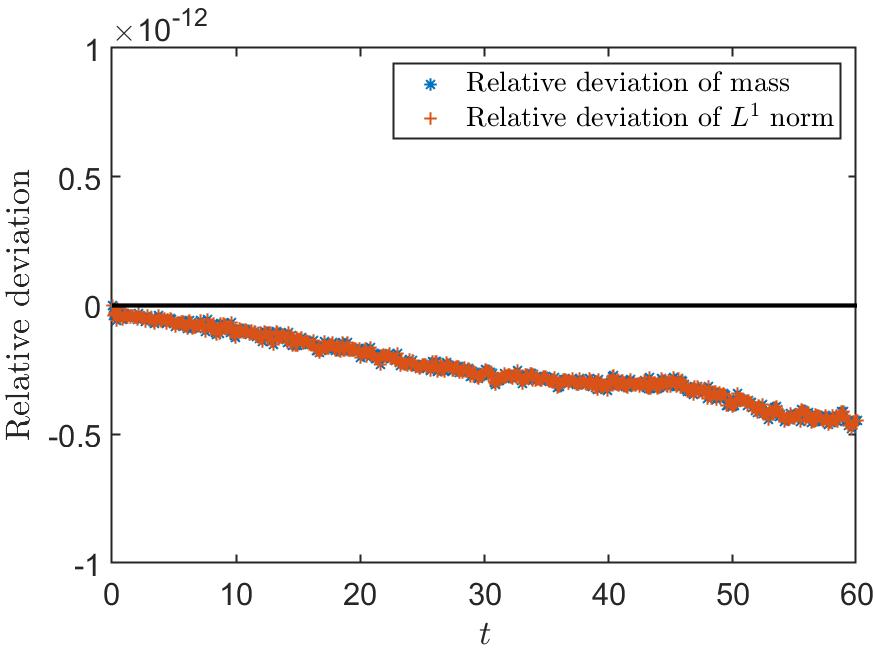}
	}
	\caption{(Strong Landau damping). Performance of mass conservation and PP properties of the SL HWENO-1 scheme for the strong Landau damping with $v_{\text{max}}=10$.}\label{fig:SLD_mass_L1}
\end{figure}	

\begin{figure}[!htbp]
	\centering
	\subfloat{
		\includegraphics[width=0.3\textwidth]{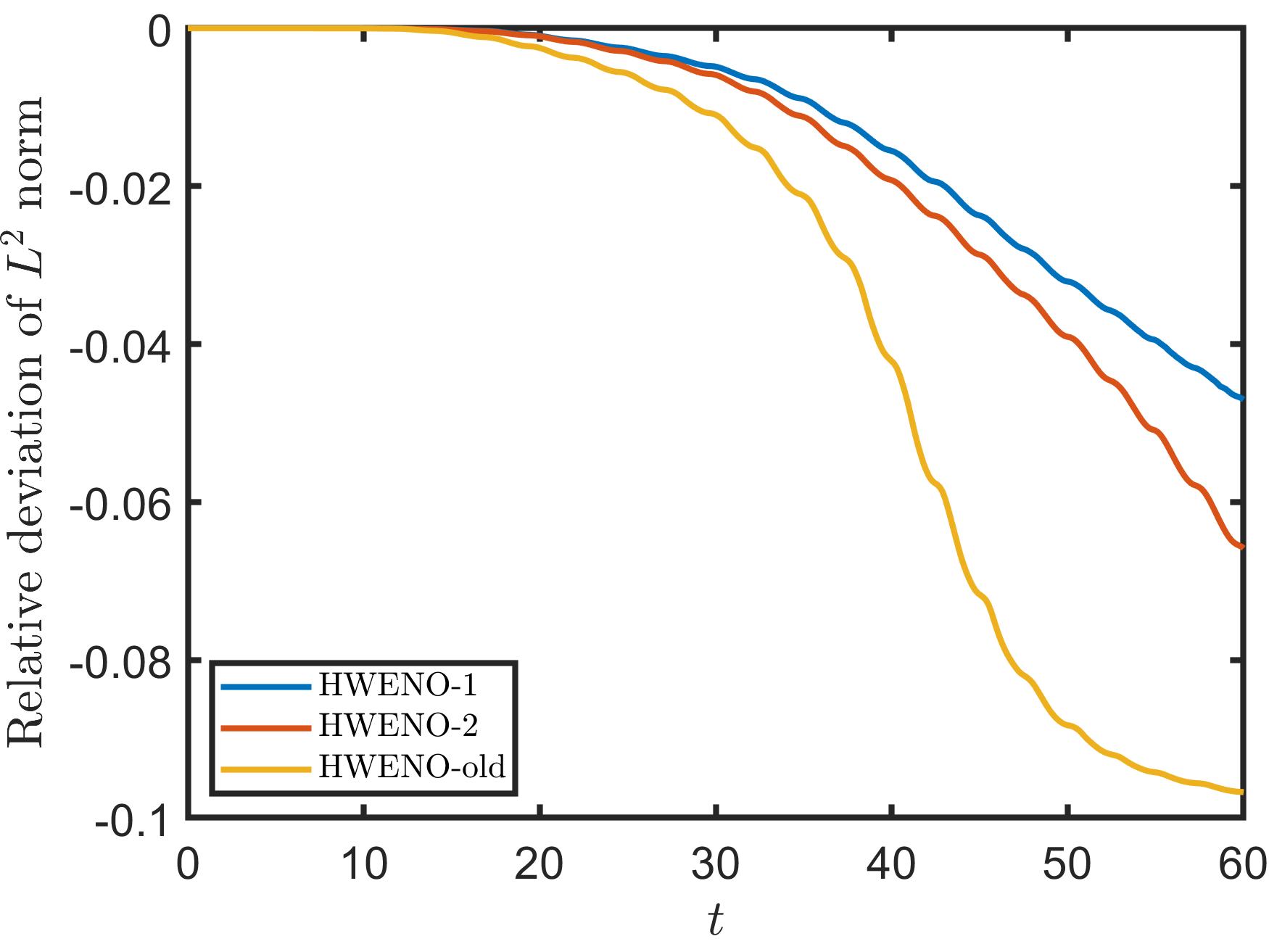}
	}
	\subfloat{
	\includegraphics[width=0.3\textwidth]{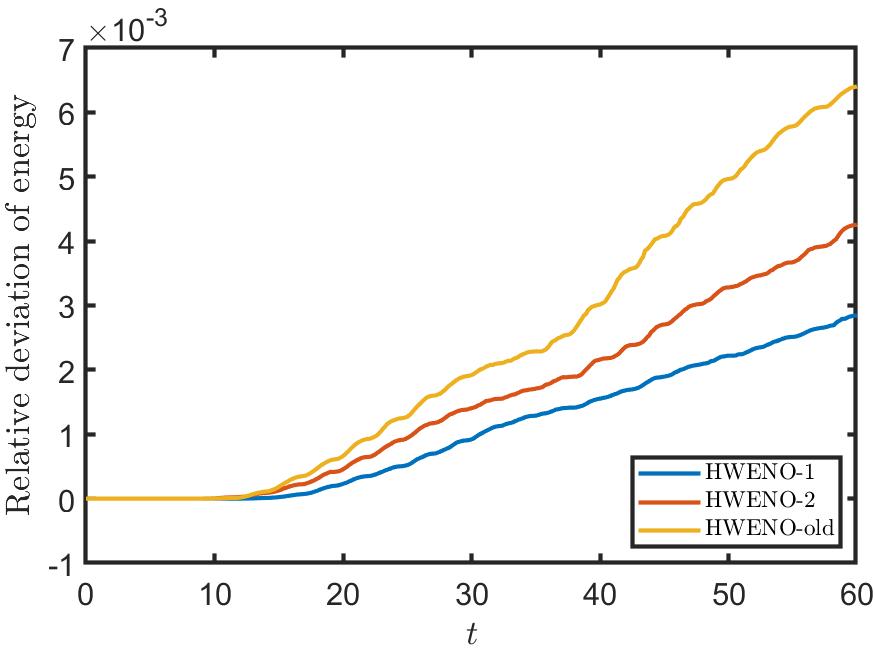}
}
	\subfloat{
	\includegraphics[width=0.3\textwidth]{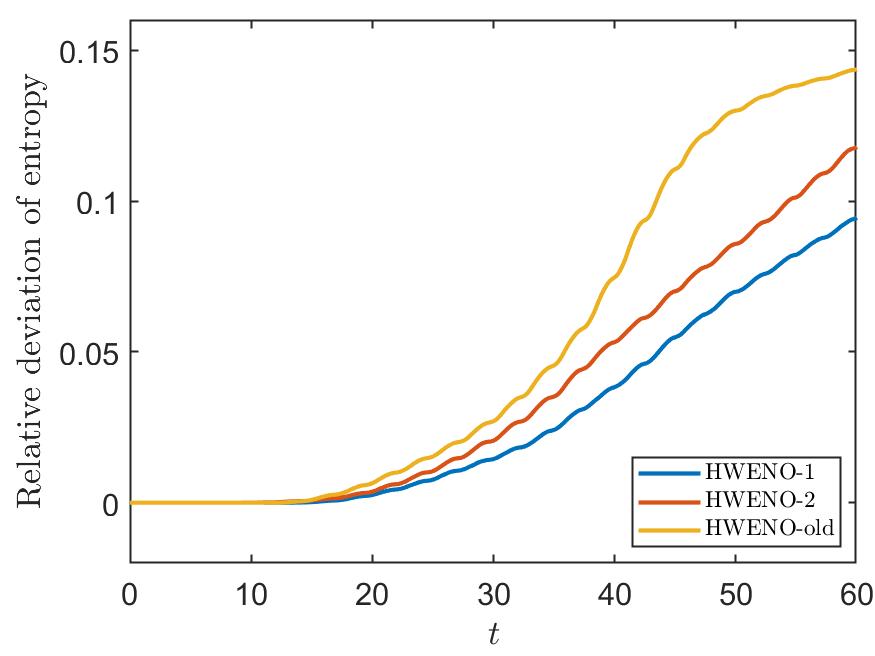}
}
	\caption{(Strong Landau damping). Relative deviations of the $L^2$ norm (left), energy (middle), and entropy (right).}\label{fig:SLD_conserve}
\end{figure}


\end{exa}

\begin{figure}[!htbp]
\centering
\subfloat{
	\includegraphics[width=0.3\textwidth]{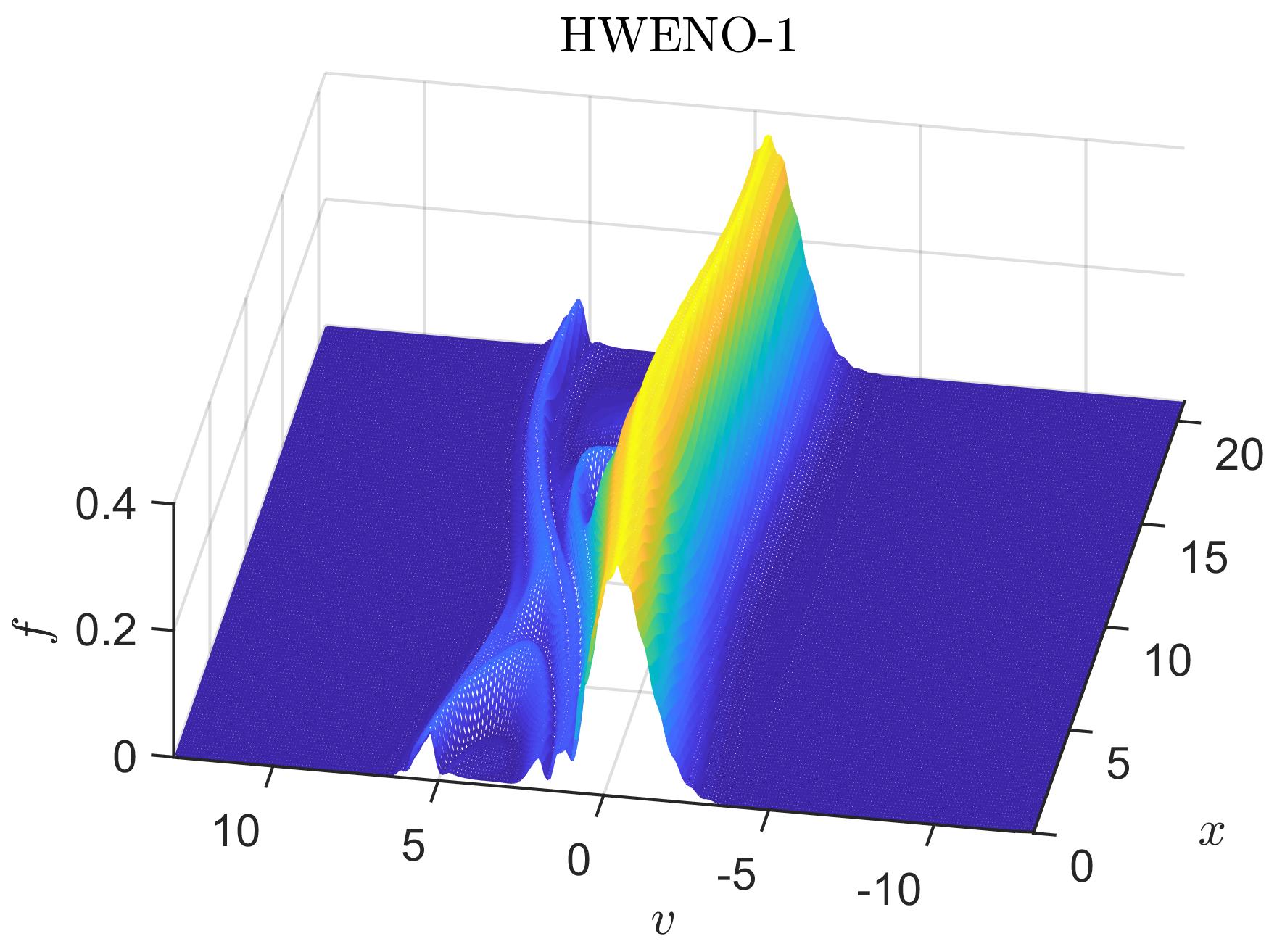}
}
\subfloat{
	\includegraphics[width=0.3\textwidth]{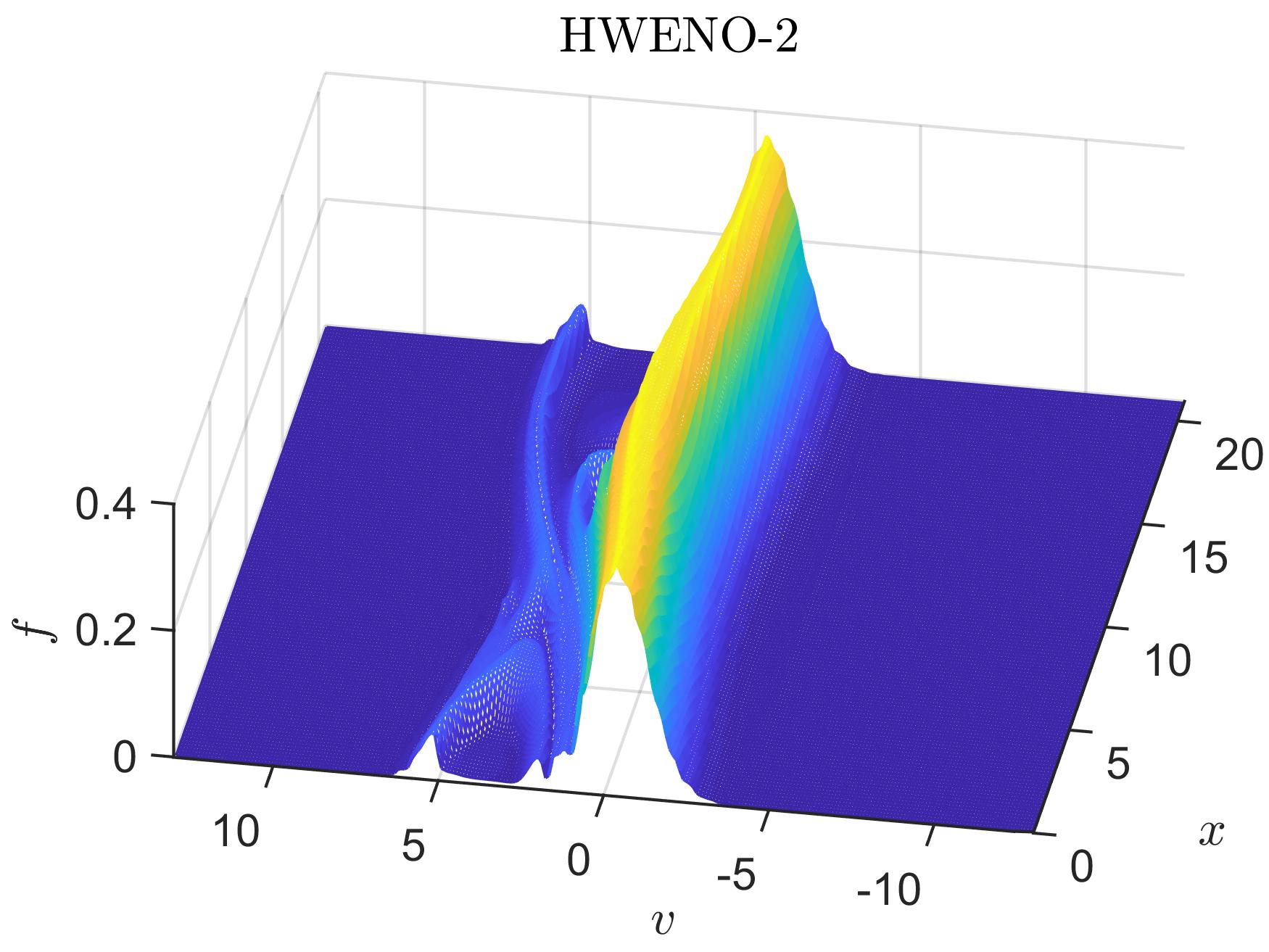}
}
\subfloat{
	\includegraphics[width=0.3\textwidth]{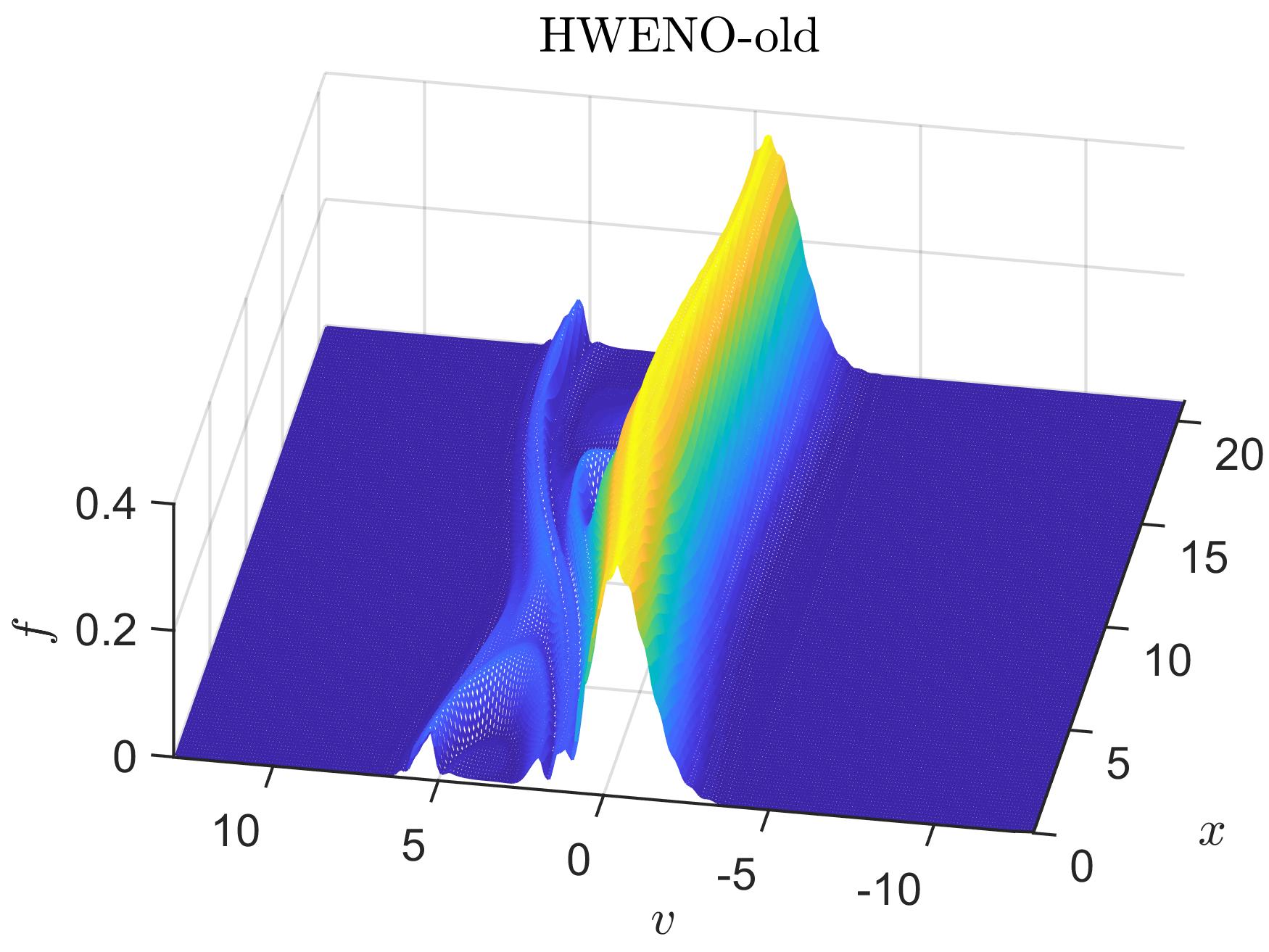}
}

\subfloat{
	\includegraphics[width=0.3\textwidth]{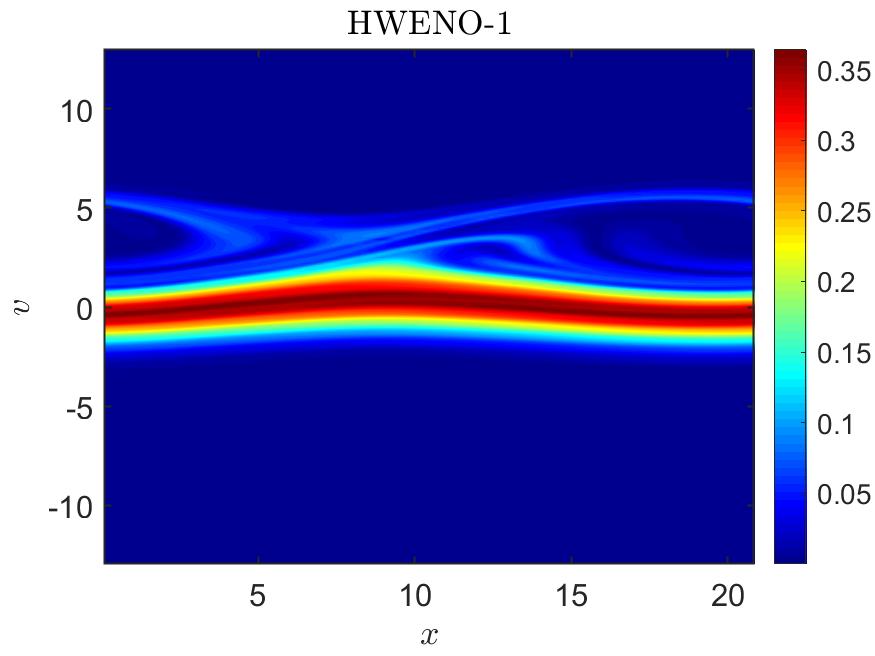}
}
\subfloat{
	\includegraphics[width=0.3\textwidth]{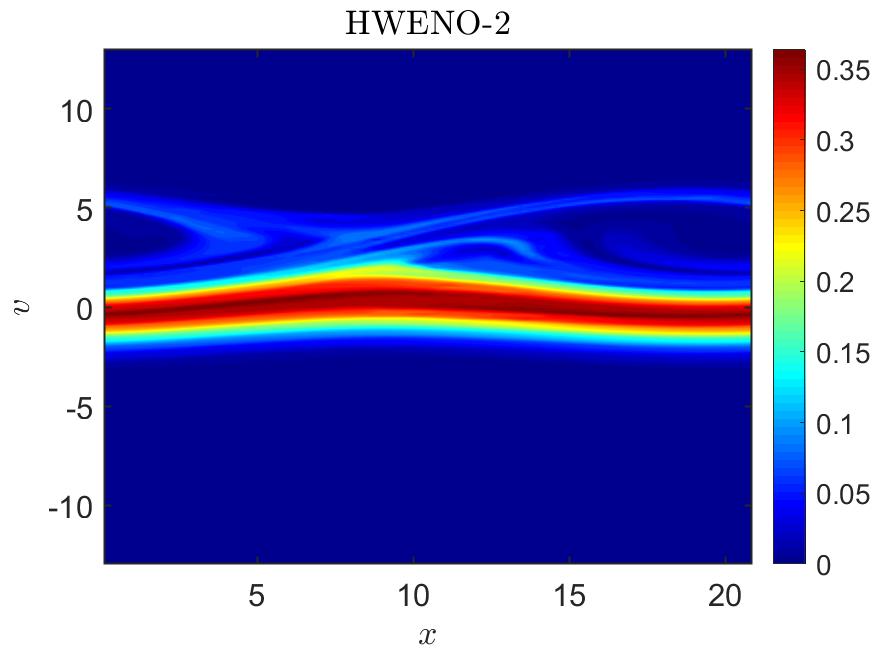}
}
\subfloat{
	\includegraphics[width=0.3\textwidth]{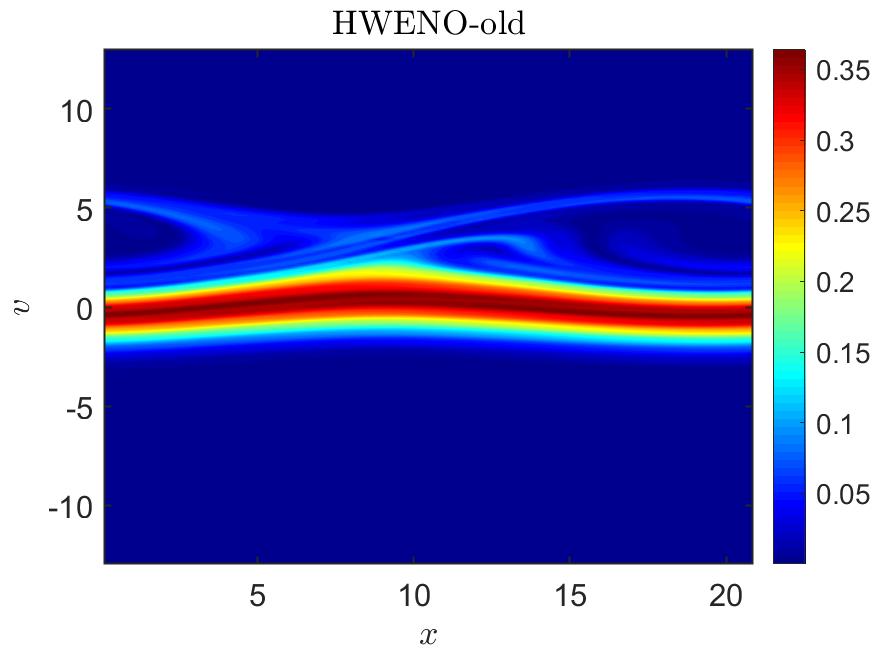}
}
\caption{(Bump-on-tail instability). The mesh plots (top) and contour plots (bottom) of the numerical solutions of the SL HWENO-1 (left), HWENO-2 (middle), and HWENO-old (right) schemes at $T = 40$.}\label{fig:BOT_plot}
\end{figure}

\begin{exa}(Bump-on-tail instability \cite{arber2002critical,CAI2018529}). Consider the bump-on-tail instability with the initial condition
	\begin{equation}\label{BOT}
		\begin{split}
		f(x,v,t=0) = \left(n_p\text{exp}\left(-\frac{v^2}{2}\right) + n_b\text{exp}\left(-\frac{(v-u)^2}{2v_t^2}\right)\right)\left(1+0.04\text{cos}(kx)\right),\\
		x\in[0,\frac{20}{3}\pi],\quad v\in [-13,13],
		\end{split}
	\end{equation}
	where $n_p=\frac{9}{10\sqrt{2\pi}}$, $n_b=\frac{2}{10\sqrt{2\pi}}$ , $u=4.5$, $v_t=0.5$ and $k=0.3$. In \Cref{fig:BOT_plot}, the meth plots and contour plots of the numerical solutions of the SL HWENO schemes at $T = 40$ are presented. The results are consistent with the existing ones in \cite{CAI2018529}.
	
\end{exa}

\subsection{Guiding center Vlasov model}

For the guiding center Vlasov model, unless specified, we set $N_x = N_y = 256$, CFL = $10.2$, and $$\Delta t = \text{CFL}/\left(\text{max}\{|E_1|\}/\Delta x+\text{max}\{|E_2|\}/\Delta y\right).$$

\begin{exa} (Kelvin-Helmholtz instability problem). Consider the guiding center Vlasov model with the periodic boundary condition and with the following initial condition
\begin{equation}\label{GC_KHI}
\begin{split}
\rho(x,y,0) = \text{sin}(y) + 0.015\text{cos}(kx),~~x\in[0,4\pi],~~y\in[0,2\pi],
\end{split}
\end{equation}
where $k=0.5$. In \Cref{tab_2-D_KHI}, we show the $L^2$ errors and corresponding orders of accuracy of the SL HWENO schemes at $T = 5$. The reference solutions are computed in the same way for the strong Landau damping. The SL HWENO schemes have similar results and the orders of accuracy are fourth. In \Cref{fig:KHI_CFL_vs_error}, we verify the fourth-order temporal accuracy for this problem.
In \Cref{fig:KHI_nonsplitting}, we show the mesh plots and contour plots of the numerical solutions of the SL HWENO schemes at $T = 40$. The numerical results are comparable to the existing ones in the literature \cite{2019ConservativeXiong,CAI2021110036}. We also provide the cross-sections of the numerical solutions of the SL HWENO schemes at $y = \pi$ with the same settings as in \Cref{fig:KHI_nonsplitting}. As shown, the SL HWENO-1, and HWENO-2 schemes have better resolution than the SL HWENO-old scheme. However, the SL HWENO-1 scheme has distinct upward and downward overshooting comparing with the other two.

The guiding center Vlasov model conserves many quantities including mass, energy, and enstrophy \cite{2019ConservativeXiong,CAI2021110036}. Similar to the strong Landau damping, we observe a $O(10^{-12})$ level of deviation of mass for this problem. We omit this result for conciseness. For energy and enstrophy, we provide the relative deviations of them for the SL HWENO schemes in \Cref{fig:KHI_conserve}. The performance for all the schemes are comparable to the existing results in the literature \cite{2019ConservativeXiong,CAI2021110036}.

\begin{table}[!htbp]
\centering
\caption{ (Kelvin-Helmholtz instability problem). $L^2$ errors and corresponding orders of accuracy of the SL HWENO schemes at $T = 5$.}\label{tab_2-D_KHI}
  \centering
\begin{tabular}{|c|cc|cc|cc|}
\hline
&\multicolumn{2}{|l|}{\textbf{HWENO-1}}&\multicolumn{2}{|l|}{\textbf{HWENO-2}}&\multicolumn{2}{|l|}{\textbf{HWENO-old}}\\
\cline{2-7}
mesh&$L^2$ error&order&$L^2$ error&order&$L^2$ error&order\\
\hline
  16$\times$  16	&	2.27E-02&       ---	 	&	2.27E-02&       ---	    	&    2.28E-02&   ---\\
  32$\times$  32	&	5.08E-03&   2.16 	 	&	5.08E-03&   2.16 	 	&    5.08E-03&   2.17\\
  64$\times$  64	&	6.54E-04&   2.96 	 	&	6.54E-04&   2.96 	 	&    6.54E-04&   2.96\\
  128$\times$  128	&	5.78E-05&   3.50 	 	&	5.78E-05&   3.50 	 	&    5.78E-05&   3.50\\
  256$\times$  256	&	4.42E-06&   3.71 	 	&	4.42E-06&   3.71 	 	&    4.42E-06&   3.71\\
\hline											
\end{tabular}
\end{table}

\begin{figure}[!htbp]
\centering
\subfloat{
\includegraphics[width=0.35\textwidth]{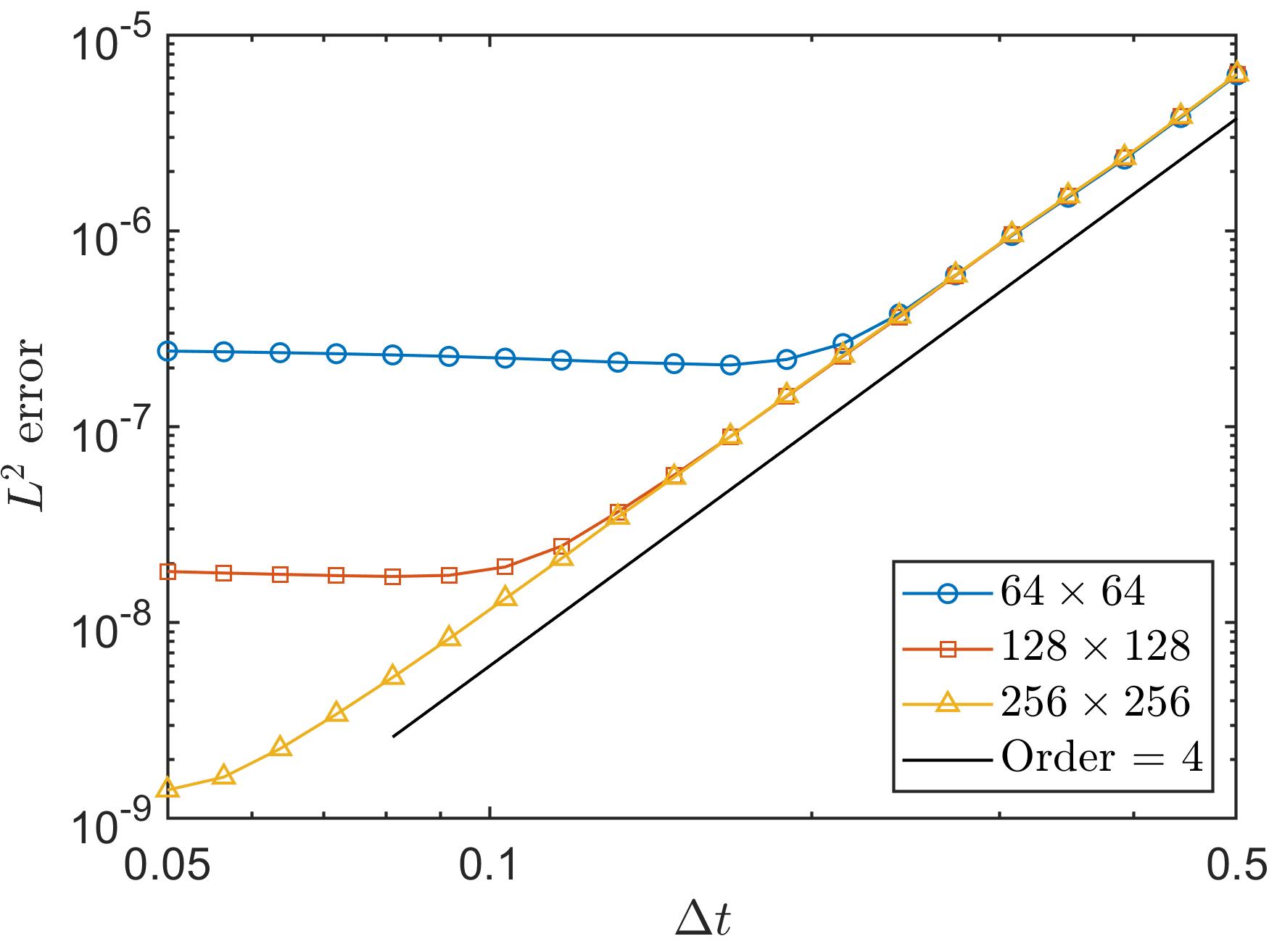}
}
\caption{(Kelvin-Helmholtz instability problem). Temporal order of accuracy of the SL HWENO-1 scheme. The three colored lines use different mesh of $64\times64$, $128\times128$, and $256\times256$. The programs stop at $T=5$.}\label{fig:KHI_CFL_vs_error}
\end{figure}	

\begin{figure}[!htbp]
\centering
\subfloat{
	\includegraphics[width=0.3\textwidth]{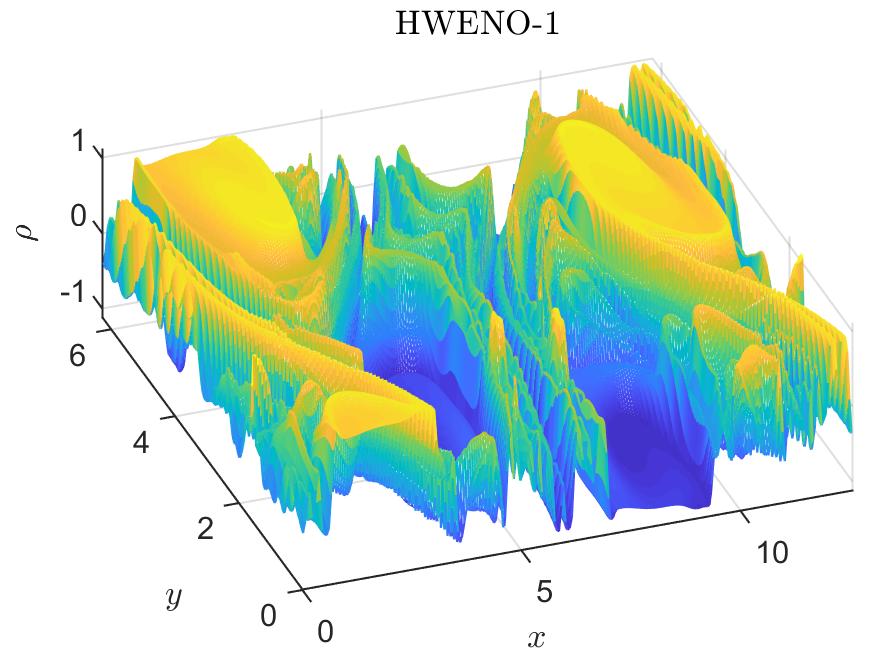}
}
\subfloat{
	\includegraphics[width=0.3\textwidth]{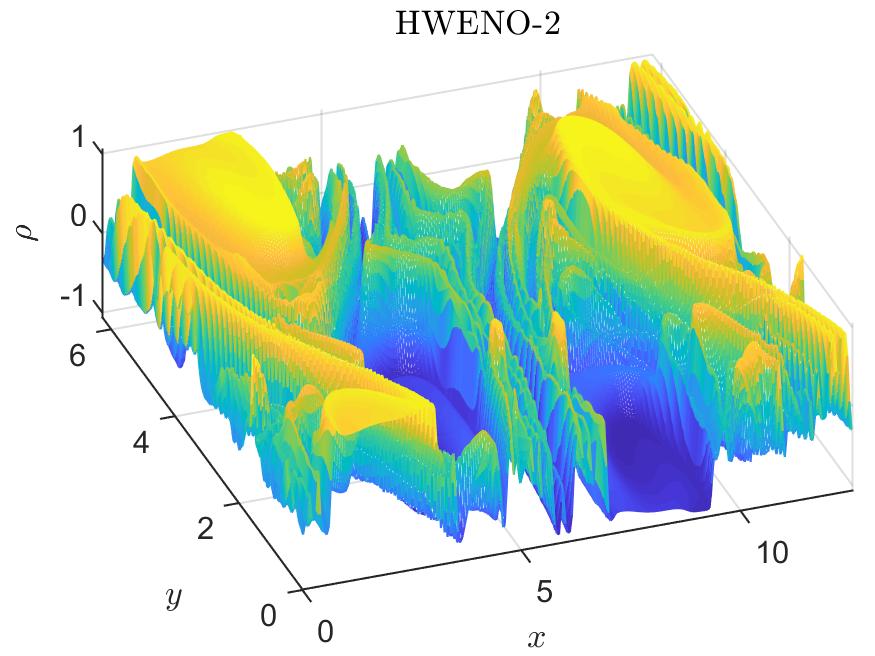}
}
\subfloat{
	\includegraphics[width=0.3\textwidth]{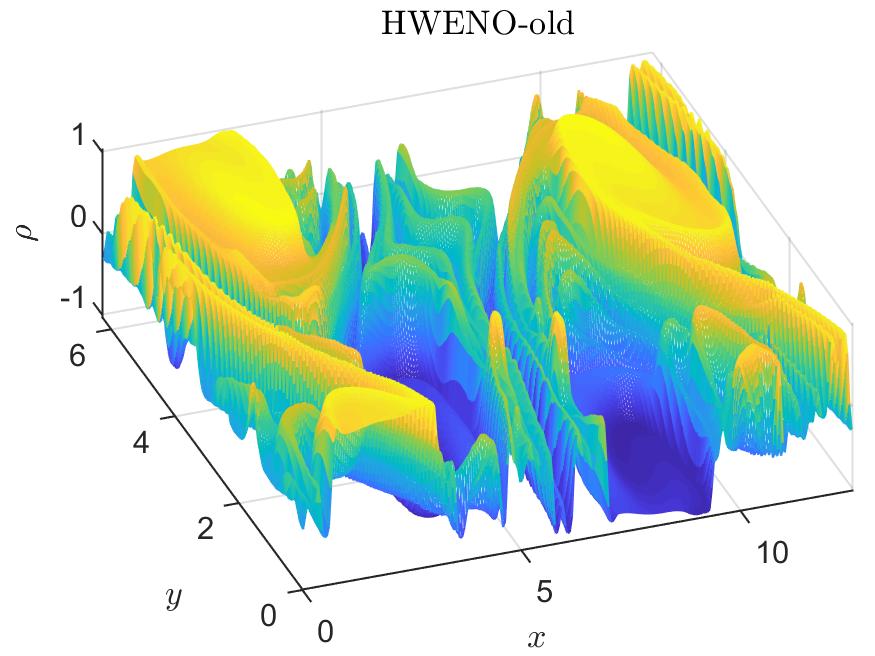}
}

\subfloat{
\includegraphics[width=0.3\textwidth]{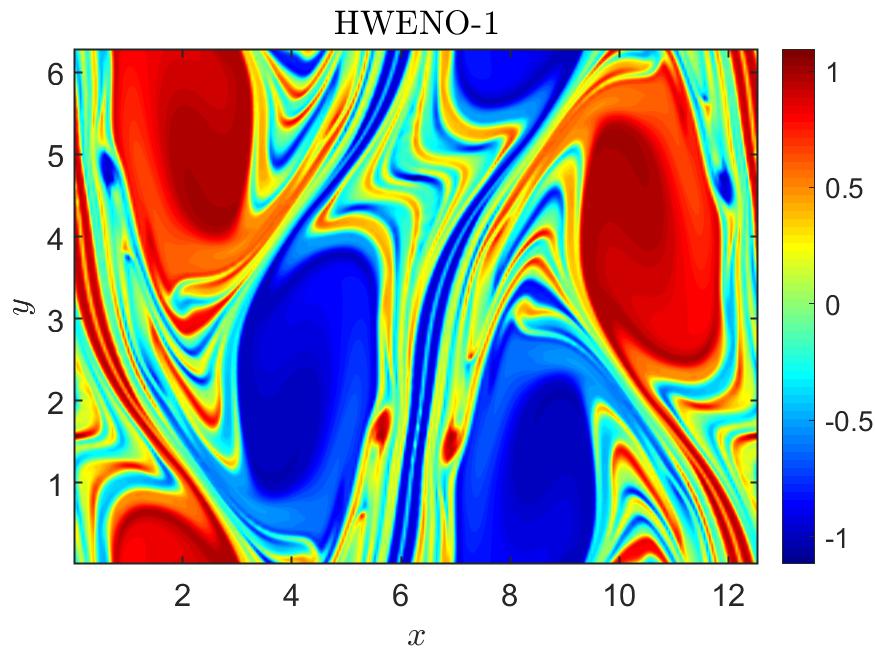}
}
\subfloat{
	\includegraphics[width=0.3\textwidth]{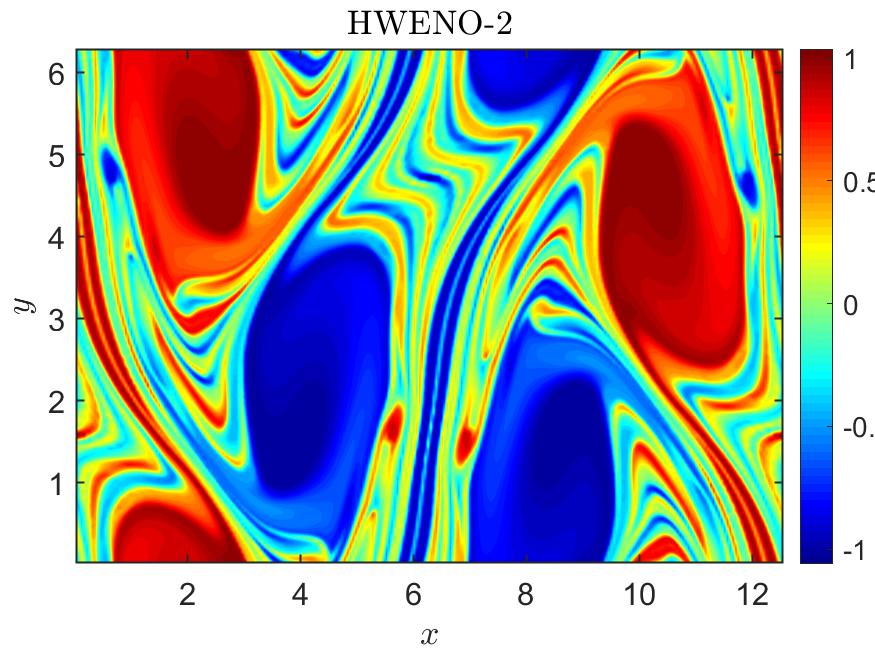}
}
\subfloat{
	\includegraphics[width=0.3\textwidth]{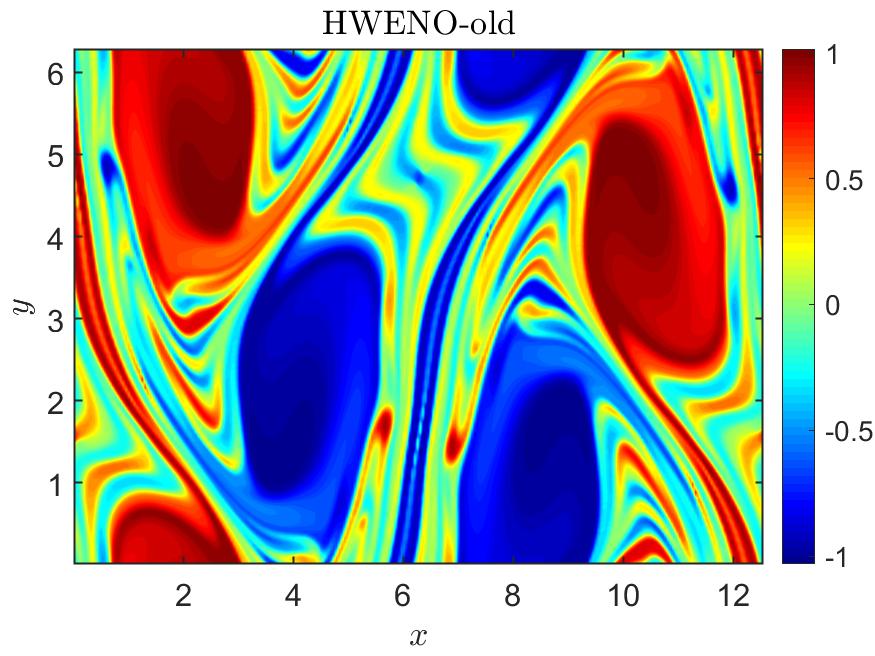}
}
\caption{(Kelvin-Helmholtz instability problem). The mesh plots (left) and contour plots of the numerical solutions of the SL HWENO-1 (left), HWENO-2 (middle), and HWENO-old (right) schemes at $T = 40$.}\label{fig:KHI_nonsplitting}
\end{figure}

\begin{figure}[!htbp]
	\centering
	\subfloat{
		\includegraphics[width=0.4\textwidth]{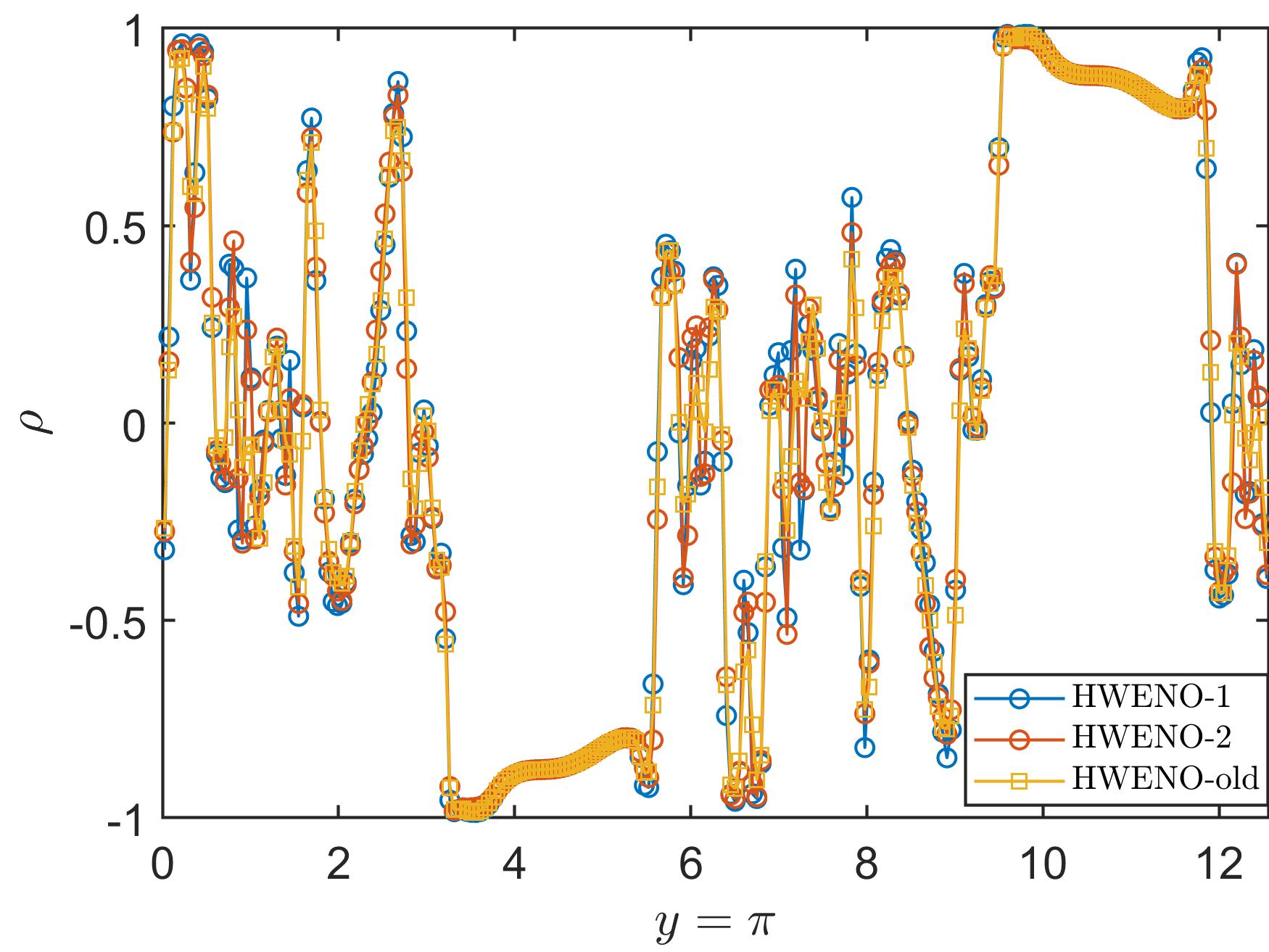}
	}
	\caption{(Kelvin-Helmholtz instability problem). Cross-sections of SL HWENO solutions at $T = 40$ and at $y = \pi$.}\label{fig:KHI_slide}
\end{figure}

\begin{figure}[!htbp]
\centering
\subfloat{
\includegraphics[width=0.35\textwidth]{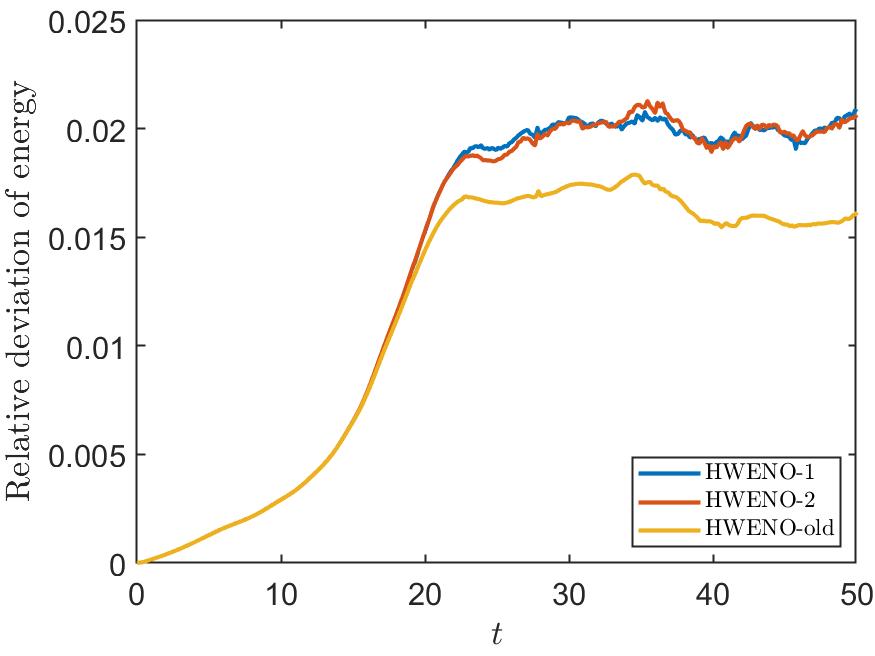}
}
\subfloat{
\includegraphics[width=0.35\textwidth]{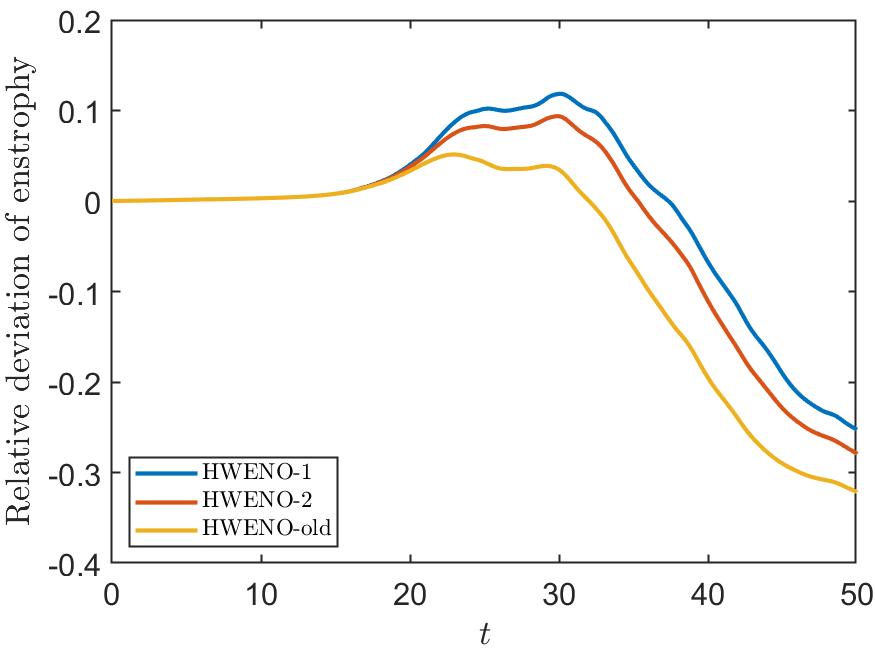}
}
\caption{(Kelvin-Helmholtz instability problem). Relative deviations of energy (left) and enstrophy (right).}\label{fig:KHI_conserve}
\end{figure}

\end{exa}

\subsection{Incompressible Euler equation}
For the incompressible Euler equation, unless specified, we set $N_x=N_y=256$, CFL = 10.2, and $$\Delta t = \text{CFL}/\left(\text{max}\{u_1\}/\Delta x+\text{max}\{u_2\}/\Delta y\right).$$
\begin{exa} (Shear flow problem). Consider the incompressible Euler equations in the domain $[0,2\pi]\times[0,2\pi]$ with the following initial condition
\begin{equation}\label{IE_VPP}
\omega(x,y,0) = \begin{cases}
\delta \text{cos}(x)-\frac1{\rho}\text{sech}^2\left(\frac{y-\pi/2}{\rho}\right),~~&\text{if}~~y\leq\pi,\\
\delta \text{cos}(x)+\frac1{\rho}\text{sech}^2\left(\frac{3\pi/2-y}{\rho}\right),~~&\text{if}~~y>\pi,
\end{cases}
\end{equation}
where $\delta = 0.05$ and $\rho = \pi/15$, and the periodic boundary condition. In \Cref{fig:SFP}, we show the mesh plots and contour plots of the numerical solutions of the SL HWENO schemes at $T = 8$. In \Cref{fig:SFP_slide}, the cross-sections of the numerical solutions of the SL HWENO schemes at $T = 8$ at $x=\pi$ are provided. The three schemes performs similarly for this problem.

The incompressible Euler equation in vorticity-stream function formulation conserves mass, energy, and enstrophy \cite{2019ConservativeXiong,CAI2021110036}. We observe a $O(10^{-12})$ level of deviation of mass. We skip this result for conciseness. In \Cref{fig:SFP_conserve}, we give the time history of relative deviations of energy and enstrophy of the SL HWENO schemes. The results are comparable to the existing ones \cite{2019ConservativeXiong,CAI2021110036}.

\begin{figure}[!htbp]
	\centering
	\subfloat{
		\includegraphics[width=0.3\textwidth]{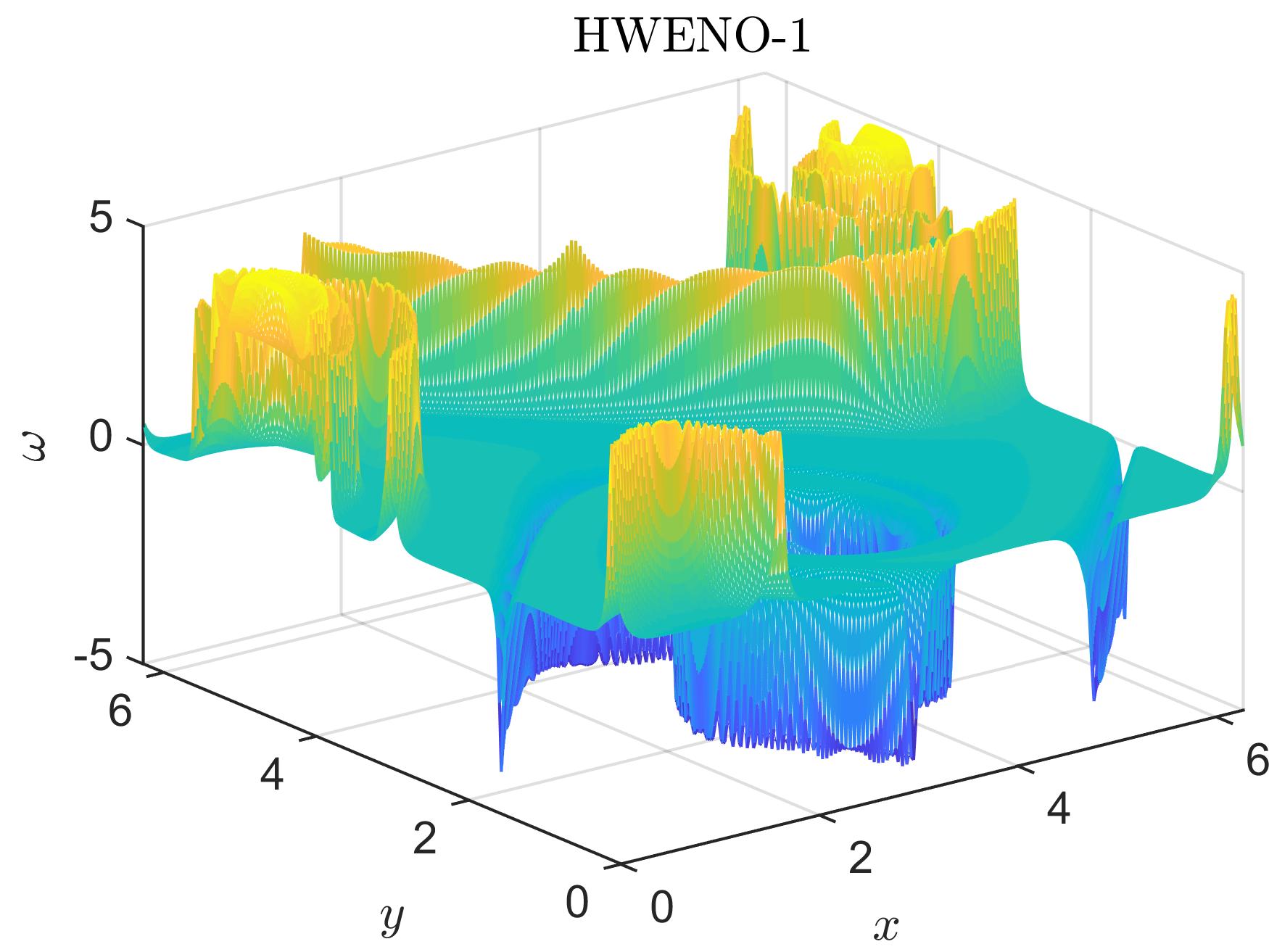}
	}
	\subfloat{
	\includegraphics[width=0.3\textwidth]{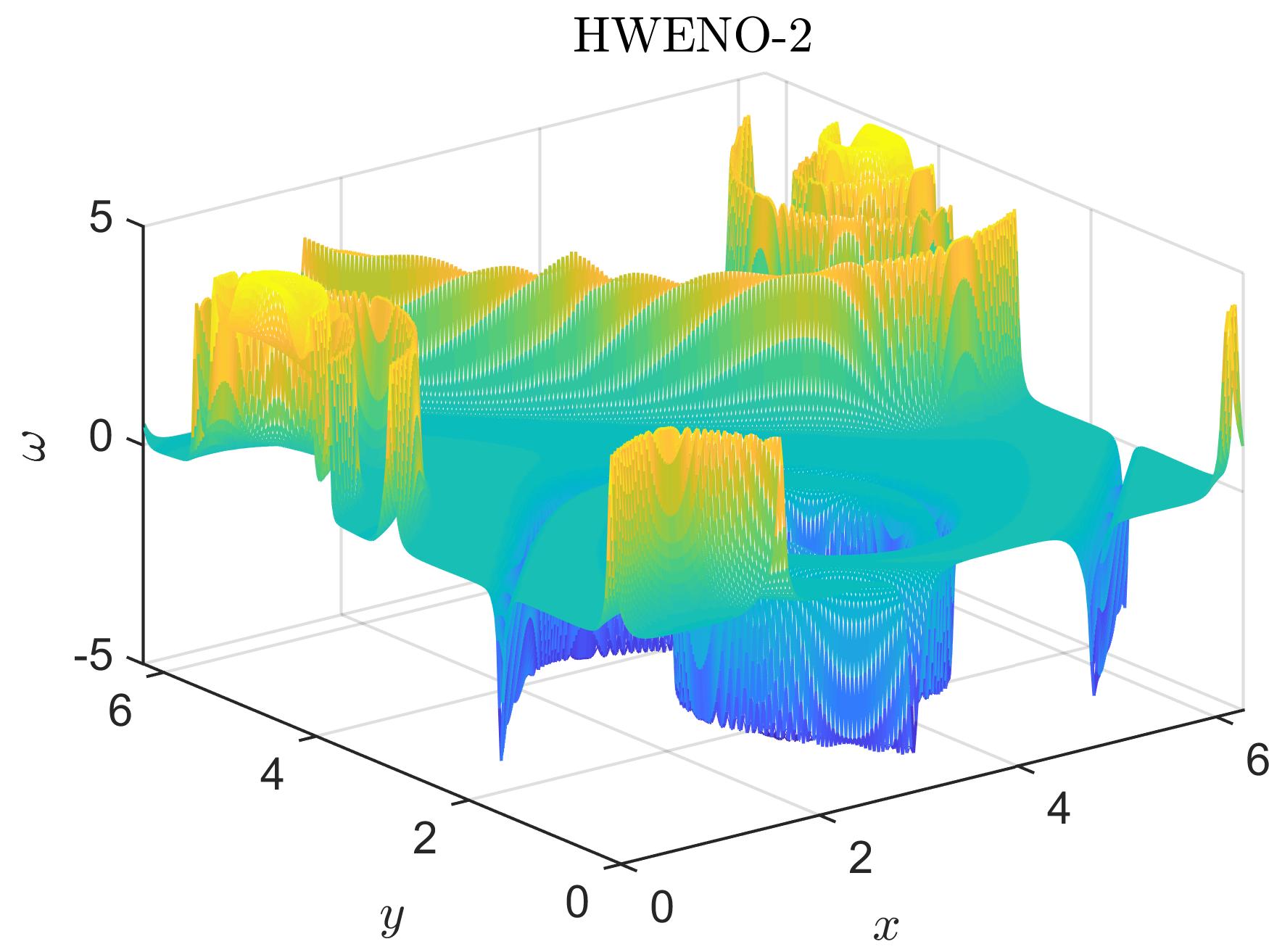}
    }
	\subfloat{
	\includegraphics[width=0.3\textwidth]{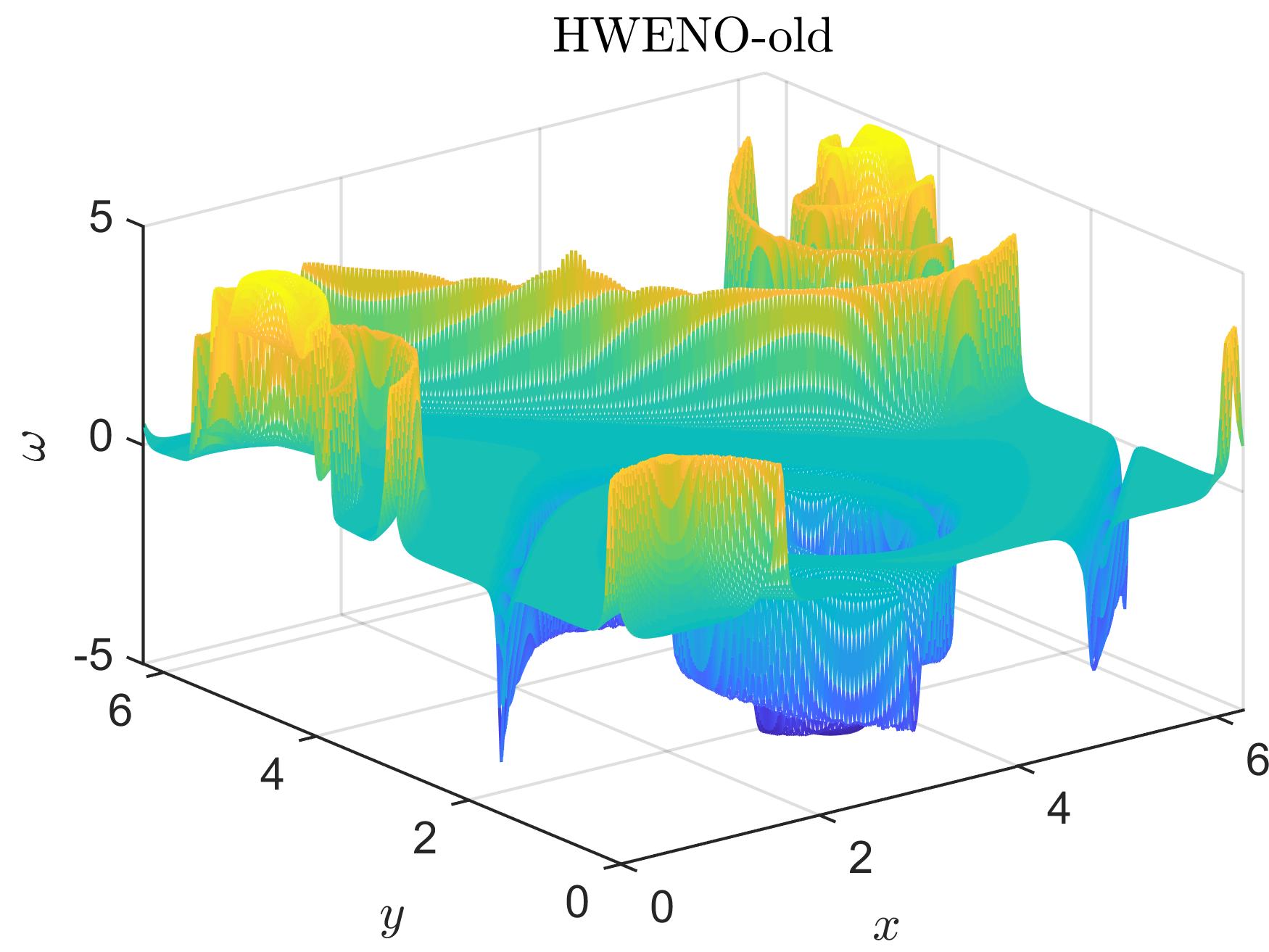}
	}
	
	\subfloat{
		\includegraphics[width=0.3\textwidth]{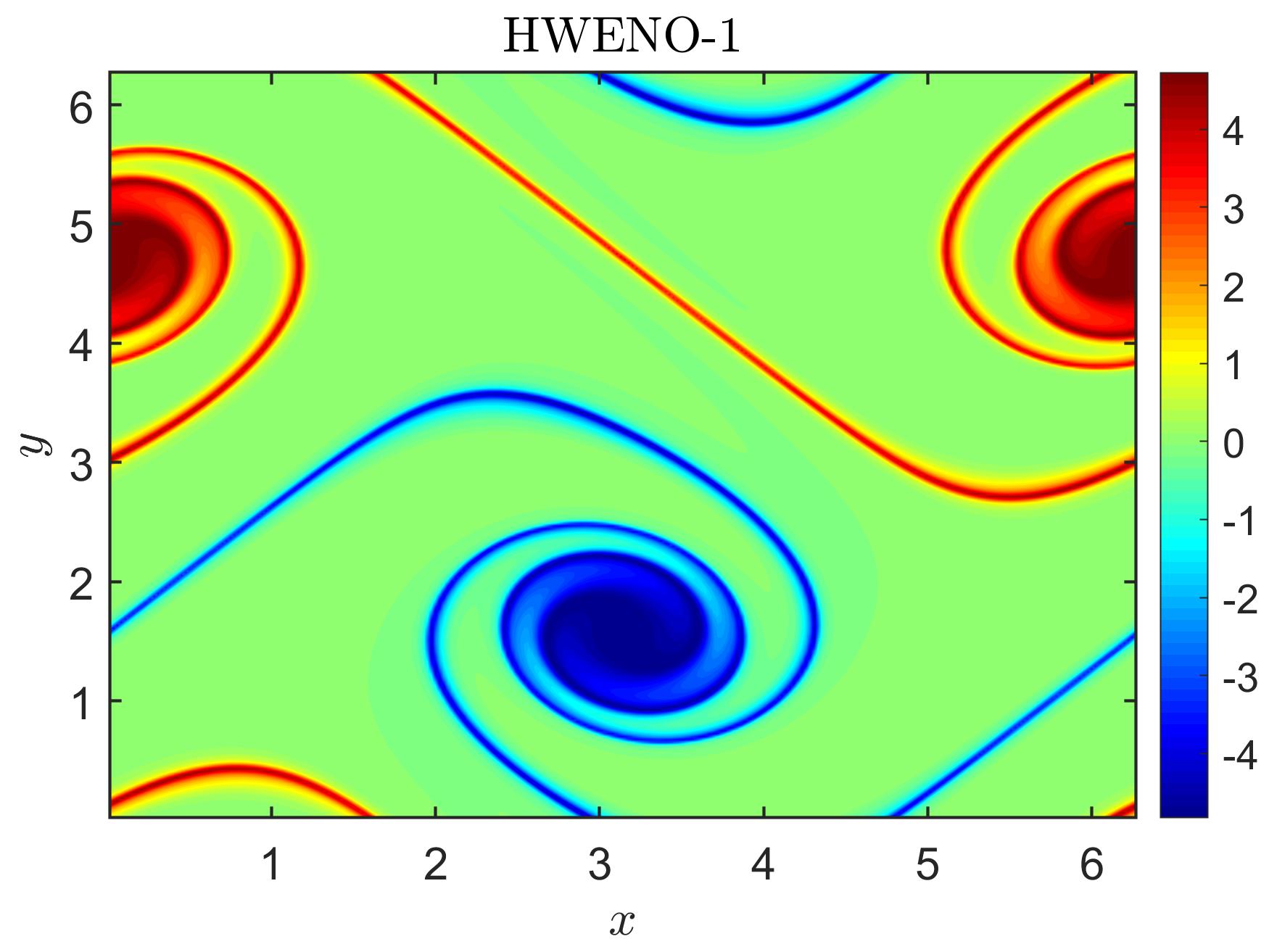}
	}
	\subfloat{
	\includegraphics[width=0.3\textwidth]{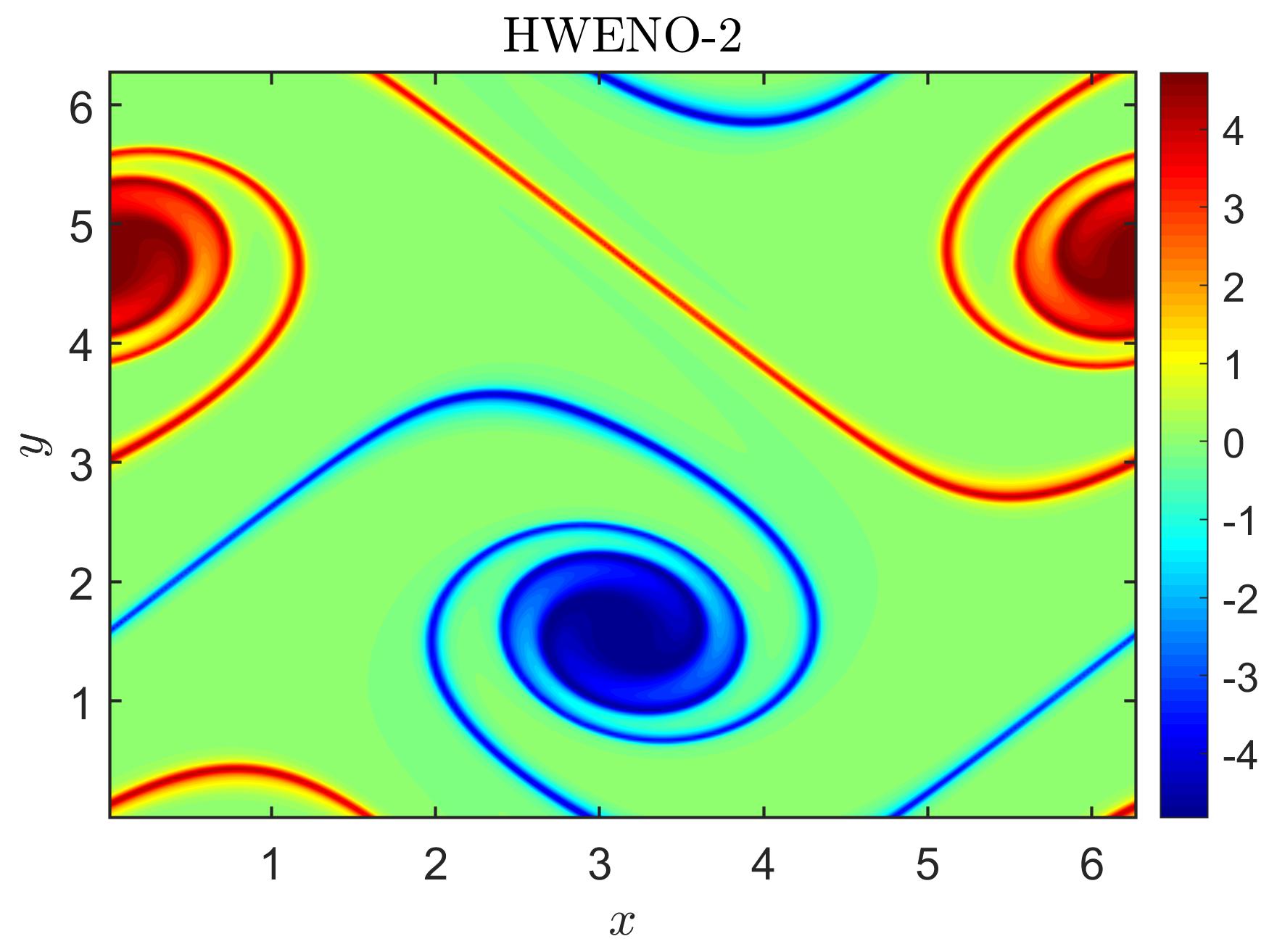}
}
	\subfloat{
	\includegraphics[width=0.3\textwidth]{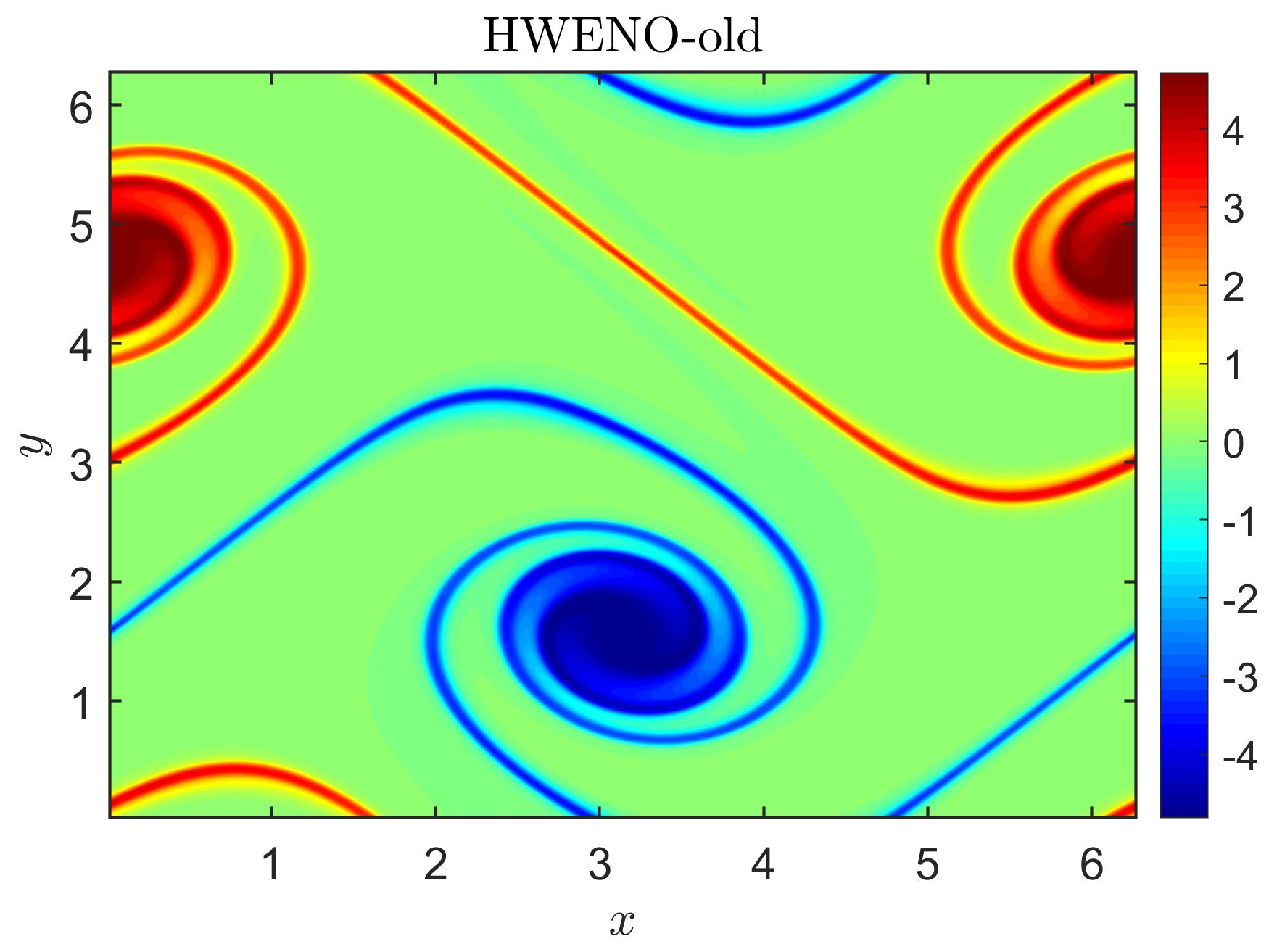}
}
	\caption{(Shear flow problem). The mesh plots (left) and contour plots of the numerical solutions of the SL HWENO-1 (left), HWENO-2 (middle), and HWENO-old (right) schemes at $T = 8$.}\label{fig:SFP}
\end{figure}

\begin{figure}[!htbp]
	\centering
	\subfloat{
		\includegraphics[width=0.4\textwidth]{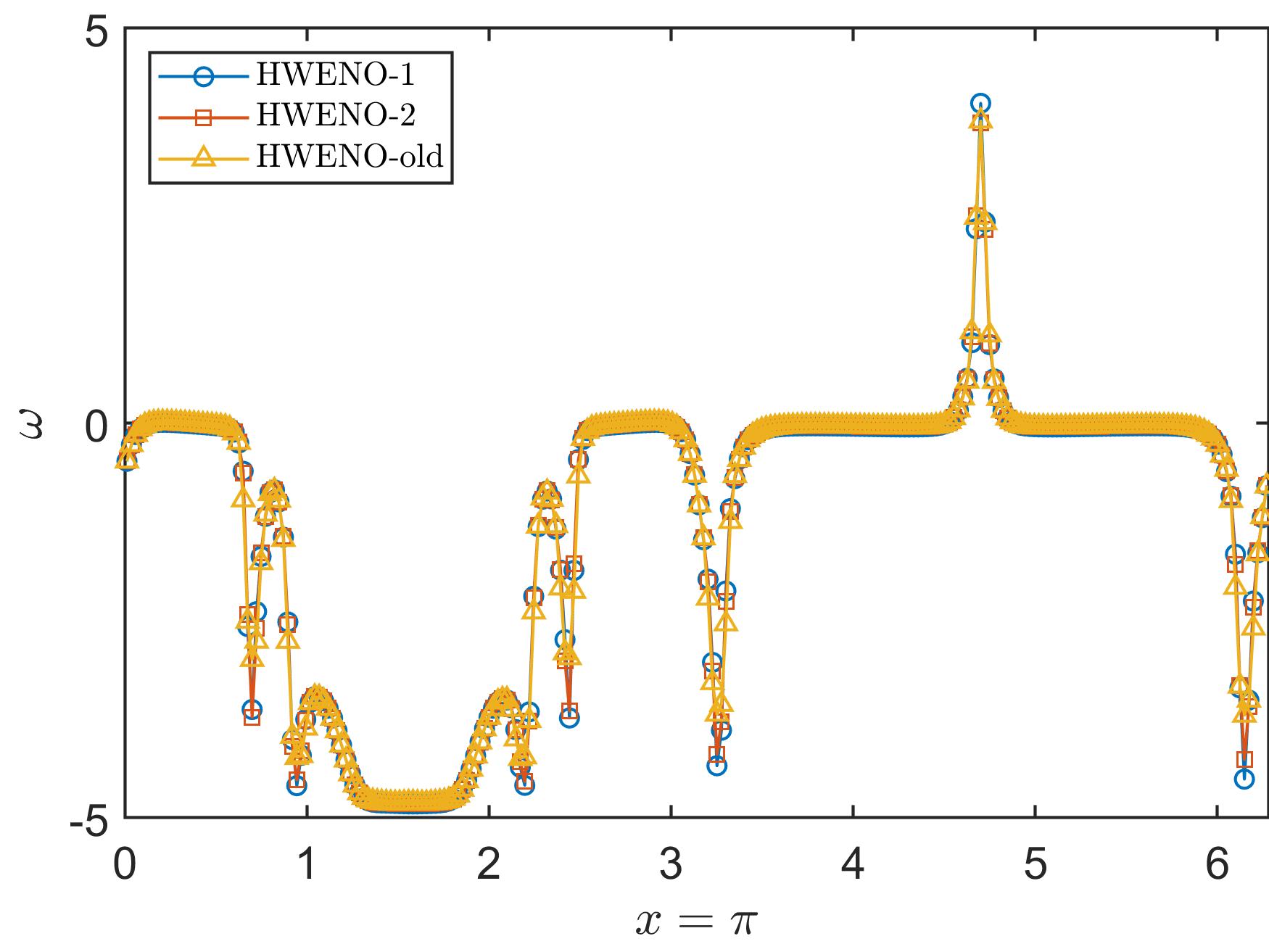}
	}
	\caption{(Shear flow problem). Cross-sections of the SL HWENO solutions at $T = 8$ and at $x = \pi$.}\label{fig:SFP_slide}
\end{figure}

\begin{figure}[!htbp]
	\centering
	\subfloat{
		\includegraphics[width=0.35\textwidth]{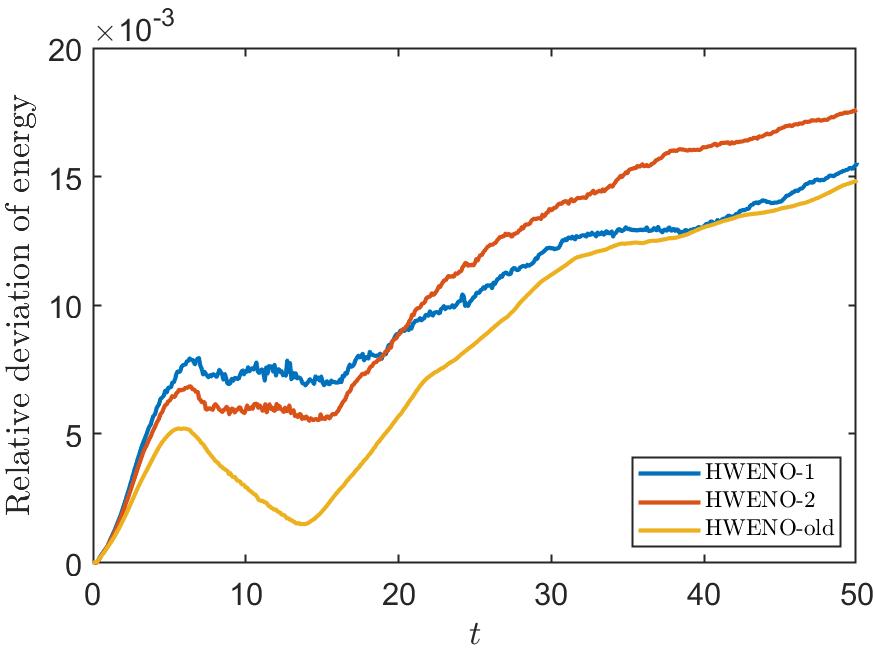}
	}
	\subfloat{
		\includegraphics[width=0.35\textwidth]{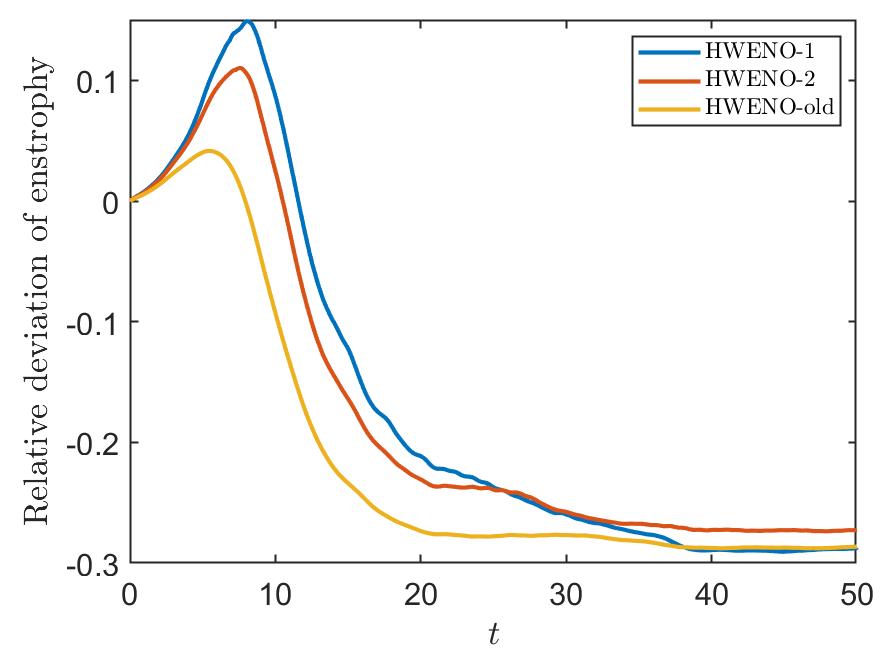}
	}
	\caption{(Shear flow problem). Relative deviations of energy (left) and enstrophy (right).}\label{fig:SFP_conserve}
\end{figure}

\end{exa}

\section{Conclusion}\label{sec:conclusion}

In this paper, we present new fourth-order SL HWENO schemes without operator splitting. Two new HWENO reconstruction methods are introduced for capturing complicated solution structures without excessive dissipation. The proposed schemes couple the weak formulation of the characteristic Galerkin method with the new HWENO reconstruction methods; the obtained SL HWENO schemes are fourth-order accurate in both space and time, mass conservative, PP, and unconditionally stable under a linearized setting.
Good numerical performance is observed through a variety of tests.

\section{Acknowledgements}
The work of the first and fourth authors was partially supported by National Natural Science Foundation (China) [Grant Number 12071392]. The work of the second author was partially support by the Foundation of Beijing Normal University [Grant Numbers 28704-310432105], UIC Research Grant with No. of UICR0700035-22 at BNU-HKBU United International College, Zhuhai, PR China,
Guangdong Higher Education Upgrading Plan (2021-2025) of "Rushing to the Top, Making Up Shortcomings and Strengthening Special Features" with No. of  UICR0400024-21.
  The work of the third author was partially support by NSF [Grant Numbers NSF-DMS-1818924 and NSF-DMS-2111253], Air Force Office of Scientific Research, United
States FA9550-18-1-0257.

\appendix
\appendixpage
\addappheadtotoc
\setcounter{section}{0}
\renewcommand\thesection{\Alph{section}}

\section{Coefficients of $\mathbf{\{q_k(x,y)\}_{k=0}^4}$}\label{appendix:coefficients}
For $q_0(x,y)$,
\begin{equation}
	\begin{split}
		&a_1^{q_0}=\overline{u}^n_{i,j},~~a_2^{q_0}=12\widetilde{\overline{v}}^n_{i,j},~~a_3^{q_0}=12\widetilde{\overline{w}}^n_{i,j},\\
		&a_4^{q_0}=\frac12\overline{u}^n_{i-1,j}-\overline{u}^n_{i,j}+\frac12\overline{u}^n_{i+1,j},\\
		&a_5^{q_0}=\frac14\overline{u}^n_{i-1,j-1}-\frac14\overline{u}^n_{i+1,j-1}-\frac14\overline{u}^n_{i-1,j+1}+\frac14\overline{u}^n_{i+1,j+1},\\
		&a_6^{q_0}=\frac12\overline{u}^n_{i,j-1}-\overline{u}^n_{i,j}+\frac12\overline{u}^n_{i,j+1},\\
		&a_7^{q_0}=-\frac5{11}\overline{u}^n_{i-1,j}+\frac5{11}\overline{u}^n_{i+1,j}-\frac{120}{11}\widetilde{\overline{v}}^n_{i,j},\\
		&a_8^{q_0}=-\frac14\overline{u}^n_{i-1,j-1}+\frac12\overline{u}^n_{i,j-1}-\frac14\overline{u}^n_{i+1,j-1}+\frac14\overline{u}^n_{i-1,j+1}-\frac12\overline{u}^n_{i,j+1}+\frac14\overline{u}^n_{i+1,j+1},\\
		&a_9^{q_0}=-\frac14\overline{u}^n_{i-1,j-1}+\frac12\overline{u}^n_{i-1,j}-\frac14\overline{u}^n_{i-1,j+1}+\frac14\overline{u}^n_{i+1,j-1}-\frac12\overline{u}^n_{i+1,j}+\frac14\overline{u}^n_{i+1,j+1},\\
		&a_{10}^{q_0}=-\frac5{11}\overline{u}^n_{i,j-1}+\frac5{11}\overline{u}^n_{i,j+1}-\frac{120}{11}\widetilde{\overline{w}}^n_{i,j}.
	\end{split}
\end{equation}
For $q_1(x,y)$,
\begin{equation}
	\begin{split}
		&a_1^{q_1}=\overline{u}^n_{i,j},~~a_2^{q_1}=12\widetilde{\overline{v}}^n_{i,j},~~a_3^{q_1}=12\widetilde{\overline{w}}^n_{i,j},\\
		&a_4^{q_1}=\overline{u}^n_{i-1,j}-\overline{u}^n_{i,j}+12\widetilde{\overline{v}}^n_{i,j},\\
		&a_5^{q_1}=\overline{u}^n_{i-1,j-1}-\overline{u}^n_{i,j-1}-\overline{u}^n_{i-1,j}+\overline{u}^n_{i,j},\\
		&a_6^{q_1}=\overline{u}^n_{i,j-1}-\overline{u}^n_{i,j}+12\widetilde{\overline{w}}^n_{i,j}.\\
	\end{split}
\end{equation}
For $q_2(x,y)$,
\begin{equation}
	\begin{split}
		&a_1^{q_2}=\overline{u}^n_{i,j},~~a_2^{q_2}=12\widetilde{\overline{v}}^n_{i,j},~~a_3^{q_2}=12\widetilde{\overline{w}}^n_{i,j},\\
		&a_4^{q_2}=-\overline{u}^n_{i,j}+\overline{u}^n_{i+1,j}-12\widetilde{\overline{v}}^n_{i,j},\\
		&a_5^{q_2}=\overline{u}^n_{i,j-1}-\overline{u}^n_{i+1,j-1}-\overline{u}^n_{i,j}+\overline{u}^n_{i+1,j},\\
		&a_6^{q_2}=\overline{u}^n_{i,j-1}-\overline{u}^n_{i,j}+12\widetilde{\overline{w}}^n_{i,j}.\\
	\end{split}
\end{equation}
For $q_3(x,y)$,
\begin{equation}
	\begin{split}
		&a_1^{q_3}=\overline{u}^n_{i,j},~~a_2^{q_3}=12\widetilde{\overline{v}}^n_{i,j},~~a_3^{q_3}=12\widetilde{\overline{w}}^n_{i,j},\\
		&a_4^{q_3}=\overline{u}^n_{i-1,j}-\overline{u}^n_{i,j}+12\widetilde{\overline{v}}^n_{i,j},\\
		&a_5^{q_3}=\overline{u}^n_{i-1,j}-\overline{u}^n_{i,j}-\overline{u}^n_{i-1,j+1}+\overline{u}^n_{i,j+1},\\
		&a_6^{q_3}=-\overline{u}^n_{i,j}+\overline{u}^n_{i,j+1}-12\widetilde{\overline{w}}^n_{i,j}.\\
	\end{split}
\end{equation}
For $q_4(x,y)$,
\begin{equation}
	\begin{split}
		&a_1^{q_4}=\overline{u}^n_{i,j},~~a_2^{q_4}=12\widetilde{\overline{v}}^n_{i,j},~~a_3^{q_4}=12\widetilde{\overline{w}}^n_{i,j},\\
		&a_4^{q_4}=-\overline{u}^n_{i,j}+\overline{u}^n_{i+1,j}-12\widetilde{\overline{v}}^n_{i,j},\\
		&a_5^{q_4}=\overline{u}^n_{i,j}-\overline{u}^n_{i+1,j}-\overline{u}^n_{i,j+1}+\overline{u}^n_{i+1,j+1},\\
		&a_6^{q_4}=-\overline{u}^n_{i,j}+\overline{u}^n_{i,j+1}-12\widetilde{\overline{w}}^n_{i,j}.\\
	\end{split}
\end{equation}

\section{Numerical proof of \Cref{prop:stability}}\label{appendix:stability}
Consider the following linear transport equation,
\begin{equation}\label{eq:const_stability}
	u_t + au_x + bu_y = 0
\end{equation}
with periodic boundary condition. Without loss of generality, we assume that $a > 0$ and $b > 0$. We define $\theta_1 = \frac{a\Delta t}{\Delta x}$ and $\theta_2 = \frac{b\Delta t}{\Delta y}$. Then \eqref{SL_formulation_numer} with $H^n(x,y)$ reconstructed by the linear reconstruction is reorganized as follows:
\begin{equation}\label{eq:SL_formulation_u_stability}
	\begin{split}
	\overline{u}^{n+1}_{i,j} &= \frac{1}{\Delta x\Delta y}\Bigg[\int_{x_{i-\frac12}-\theta_1\Delta x}^{x_{i-\frac12}-\lfloor\theta_1\rfloor\Delta x}\int_{y_{j-\frac12}-\theta_2\Delta y}^{y_{j-\frac12}-\lfloor\theta_2\rfloor\Delta y}H^{(i-1-\lfloor\theta_1\rfloor,j-1-\lfloor\theta_2\rfloor)}(x,y)dxdy\\
	&+ \int_{x_{i-\frac12}-\lfloor\theta_1\rfloor\Delta x}^{x_{i-\frac12}+(1-\theta_1)\Delta x}\int_{y_{j-\frac12}-\theta_2\Delta y}^{y_{j-\frac12}-\lfloor\theta_2\rfloor\Delta y}H^{(i-\lfloor\theta_1\rfloor,j-1-\lfloor\theta_2\rfloor)}(x,y)dxdy\\
	&+ \int_{x_{i-\frac12}-\theta_1\Delta x}^{x_{i-\frac12}-\lfloor\theta_1\rfloor\Delta x}\int_{y_{j-\frac12}-\lfloor\theta_2\rfloor\Delta y}^{y_{j-\frac12}+(1-\theta_2)\Delta y}H^{(i-1-\lfloor\theta_1\rfloor,j-\lfloor\theta_2\rfloor)}(x,y)dxdy\\
	&+ \int_{x_{i-\frac12}-\lfloor\theta_1\rfloor\Delta x}^{x_{i-\frac12}+(1-\theta_1)\Delta x}\int_{y_{j-\frac12}-\lfloor\theta_2\rfloor\Delta y}^{y_{j-\frac12}+(1-\theta_2)\Delta y}H^{(i-\lfloor\theta_1\rfloor,j-\lfloor\theta_2\rfloor)}(x,y)dxdy,
	\end{split}
\end{equation}
\begin{equation}\label{eq:SL_formulation_v_stability}
	\begin{split}
		\overline{v}^{n+1}_{i,j} &= \frac{1}{\Delta x\Delta y}\Bigg[\int_{x_{i-\frac12}-\theta_1\Delta x}^{x_{i-\frac12}-\lfloor\theta_1\rfloor\Delta x}\int_{y_{j-\frac12}-\theta_2\Delta y}^{y_{j-\frac12}-\lfloor\theta_2\rfloor\Delta y}H^{(i-1-\lfloor\theta_1\rfloor,j-1-\lfloor\theta_2\rfloor)}(x,y)\left(\frac{x+\theta_1\Delta x-x_i}{\Delta x}\right)dxdy\\
		&+ \int_{x_{i-\frac12}-\lfloor\theta_1\rfloor\Delta x}^{x_{i-\frac12}+(1-\theta_1)\Delta x}\int_{y_{j-\frac12}-\theta_2\Delta y}^{y_{j-\frac12}-\lfloor\theta_2\rfloor\Delta y}H^{(i-\lfloor\theta_1\rfloor,j-1-\lfloor\theta_2\rfloor)}(x,y)\left(\frac{x+\theta_1\Delta x-x_i}{\Delta x}\right)dxdy\\
		&+ \int_{x_{i-\frac12}-\theta_1\Delta x}^{x_{i-\frac12}-\lfloor\theta_1\rfloor\Delta x}\int_{y_{j-\frac12}-\lfloor\theta_2\rfloor\Delta y}^{y_{j-\frac12}+(1-\theta_2)\Delta y}H^{(i-1-\lfloor\theta_1\rfloor,j-\lfloor\theta_2\rfloor)}(x,y)\left(\frac{x+\theta_1\Delta x-x_i}{\Delta x}\right)dxdy\\
		&+ \int_{x_{i-\frac12}-\lfloor\theta_1\rfloor\Delta x}^{x_{i-\frac12}+(1-\theta_1)\Delta x}\int_{y_{j-\frac12}-\lfloor\theta_2\rfloor\Delta y}^{y_{j-\frac12}+(1-\theta_2)\Delta y}H^{(i-\lfloor\theta_1\rfloor,j-\lfloor\theta_2\rfloor)}(x,y)\left(\frac{x+\theta_1\Delta x-x_i}{\Delta x}\right)dxdy,
	\end{split}
\end{equation}
and
\begin{equation}\label{eq:SL_formulation_w_stability}
	\begin{split}
		\overline{w}^{n+1}_{i,j} &= \frac{1}{\Delta x\Delta y}\Bigg[\int_{x_{i-\frac12}-\theta_1\Delta x}^{x_{i-\frac12}-\lfloor\theta_1\rfloor\Delta x}\int_{y_{j-\frac12}-\theta_2\Delta y}^{y_{j-\frac12}-\lfloor\theta_2\rfloor\Delta y}H^{(i-1-\lfloor\theta_1\rfloor,j-1-\lfloor\theta_2\rfloor)}(x,y)\left(\frac{y+\theta_2\Delta y-y_j}{\Delta y}\right)dxdy\\
		&+ \int_{x_{i-\frac12}-\lfloor\theta_1\rfloor\Delta x}^{x_{i-\frac12}+(1-\theta_1)\Delta x}\int_{y_{j-\frac12}-\theta_2\Delta y}^{y_{j-\frac12}-\lfloor\theta_2\rfloor\Delta y}H^{(i-\lfloor\theta_1\rfloor,j-1-\lfloor\theta_2\rfloor)}(x,y)\left(\frac{y+\theta_2\Delta y-y_j}{\Delta y}\right)dxdy\\
		&+ \int_{x_{i-\frac12}-\theta_1\Delta x}^{x_{i-\frac12}-\lfloor\theta_1\rfloor\Delta x}\int_{y_{j-\frac12}-\lfloor\theta_2\rfloor\Delta y}^{y_{j-\frac12}+(1-\theta_2)\Delta y}H^{(i-1-\lfloor\theta_1\rfloor,j-\lfloor\theta_2\rfloor)}(x,y)\left(\frac{y+\theta_2\Delta y-y_j}{\Delta y}\right)dxdy\\
		&+ \int_{x_{i-\frac12}-\lfloor\theta_1\rfloor\Delta x}^{x_{i-\frac12}+(1-\theta_1)\Delta x}\int_{y_{j-\frac12}-\lfloor\theta_2\rfloor\Delta y}^{y_{j-\frac12}+(1-\theta_2)\Delta y}H^{(i-\lfloor\theta_1\rfloor,j-\lfloor\theta_2\rfloor)}(x,y)\left(\frac{y+\theta_2\Delta y-y_j}{\Delta y}\right)dxdy,
	\end{split}
\end{equation}
where $\lfloor\cdot\rfloor$ represents the largest integer smaller than the input number. Now, we assume that
\begin{equation}\label{eq:von_Neuemann_assume1}
	\overline{u}^n_{p,q} = \overline{u}^ne^{I\xi_1p\Delta x}e^{I\xi_2q\Delta y},\quad \overline{v}^n_{p,q} = \overline{v}^ne^{I\xi_1p\Delta x}e^{I\xi_2q\Delta y},\quad \overline{w}^n_{p,q} = \overline{w}^ne^{I\xi_1p\Delta x}e^{I\xi_2q\Delta y}\quad \forall p,q,
\end{equation}
and
\begin{equation}\label{eq:von_Neuemann_assume2}
	\overline{u}^{n+1}_{i,j}=\overline{u}^{n+1}e^{I\xi_1i\Delta x}e^{I\xi_2j\Delta y},\quad \overline{v}^{n+1}_{i,j}=\overline{v}^{n+1}e^{I\xi_1i\Delta x}e^{I\xi_2j\Delta y},\quad \overline{w}^{n+1}_{i,j}=\overline{w}^{n+1}e^{I\xi_1i\Delta x}e^{I\xi_2j\Delta y},
\end{equation}
where $I = \sqrt{-1}$. Submitting \eqref{eq:von_Neuemann_assume1} and \eqref{eq:von_Neuemann_assume2} into \eqref{eq:SL_formulation_u_stability}-\eqref{eq:SL_formulation_w_stability}, we find that
\begin{equation}
\left[\begin{array}{c}
	\overline{u}^{n+1}\\
	\overline{v}^{n+1}\\
	\overline{w}^{n+1}
\end{array}\right]=
	 A(\theta_1,\theta_2,\xi_1,\xi_2)
	 \left[\begin{array}{c}
	 	\overline{u}^{n}\\
	 	\overline{v}^{n}\\
	 	\overline{w}^{n}
	 \end{array}\right],
\end{equation}
where $A(\theta_1,\theta_2,\xi_1,\xi_2)$ is the $3\times3$ amplification matrix. We skip the explicit expression of $A(\theta_1,\theta_2,\xi_1,\xi_2)$ for conciseness since it is very complicated. We denote the spectral radius of $A(\theta_1,\theta_2,\xi_1,\xi_2)$ by $\rho(A(\theta_1,\theta_2,\xi_1,\xi_2))$. With basic algebraic manipulation, we have \begin{equation}
	A(\theta_1,\theta_2,\xi_1,\xi_2) = e^{-I\xi_1\lfloor\theta_1\rfloor\Delta x}e^{-I\xi_2\lfloor\theta_2\rfloor\Delta y}A(\theta_1-\lfloor\theta_1\rfloor,\theta_2-\lfloor\theta_2\rfloor,\xi_1,\xi_2).
\end{equation}
Hence, by von Nuemann analysis, it is sufficient to verify that $\rho(A(\theta_1,\theta_2,\xi_1,\xi_2)) \leq 1$ for any $\theta_1, \theta_2\in[0,1]$ and $\xi_1\Delta x, \xi_2\Delta y\in [0,2\pi]$. It impossible to provide the theoretical expression of $\rho(A(\theta_1,\theta_2,\xi_1,\xi_2))$. We numerically verify this relation by sampling 1000 uniform points over each $\theta_1$, $\theta_2$, $\xi_1\Delta x$, $\xi_2\Delta y$ domain. We find that all the $\rho(A(\cdot,\cdot,\cdot,\cdot))$ values computed by the sampling points are not greater than 1, which validates \Cref{prop:stability}.


  \bibliographystyle{abbrv}
  \bibliography{ref}





\end{document}